\documentclass[12pt]{article}
\usepackage{amssymb, latexsym}
\title{Rational invariant tori, phase space tunneling, and spectra for non-selfadjoint operators in dimension 2}
\author{Michael Hitrik\\Department of Mathematics \\University of California \\ Los Angeles
\\CA 90095-1555, USA\\hitrik@math.ucla.edu \and
Johannes Sj\"ostrand\\CMLS\\Ecole Polytechnique\\FR--91128 Palaiseau\\France \\and
UMR 7640 CNRS
\\johannes@math.polytechnique.fr}
\date{}

\newcommand{\real}{\mbox{\bf R}}
\def\abs#1{\left|#1\right|}
\def\begeq{\begin{equation}}
\def\endeq{\end{equation}}
\renewcommand{\Im}{\mbox{\rm Im\,}}
\def\wrtext#1{\relax\ifmmode{\leavevmode\hbox{#1}}\else{#1}\fi}

\newcommand{\eps}{\varepsilon}
\def\part#1{\frac{\partial}{\partial #1}}

\def\norm#1{||\,#1\,||}

\newcommand{\comp}{\mbox{\bf C}}
\newcommand{\z}{\mbox{\bf Z}}
\newcommand{\nat}{\mbox{\bf N}}

\renewcommand{\Re}{\mbox{\rm Re\,}}
\renewcommand{\exp}{\mbox{\rm exp\,}}
\newcommand{\supp}{\mbox{\rm supp}}

\newcommand{\neigh}{neighborhood}

\newcommand{\cantransf}{canonical transformation}

\newtheorem{dref}{Definition}[section]
\newtheorem{lemma}[dref]{Lemma}
\newtheorem{theo}[dref]{Theorem}
\newtheorem{prop}[dref]{Proposition}

\newenvironment{proof}{\vspace{.3cm}\noindent{{\em Proof:}}}{\hfill$\Box$
\vspace{.2cm}}

\begin{document}
\maketitle

\begin{abstract}
We study spectral asymptotics and resolvent bounds for non-selfadjoint
perturbations of selfadjoint $h$-pseudodifferential operators in
dimension 2, assuming that the classical flow of the unperturbed
part is completely integrable. Spectral contributions coming
from rational invariant Lagrangian tori are analyzed. Estimating the tunnel effect
between strongly irrational (Diophantine) and rational tori, we obtain an accurate description of the spectrum in a suitable
complex window, provided that the strength of the non-selfadjoint perturbation $\gg h$ (or sometimes $\gg h^2$)
is not too large.
\end{abstract}


\tableofcontents
\section{Introduction and statement of the main results}
\setcounter{equation}{0}
The present paper is the second one in a series of works dealing with the spectral analysis of small non-selfadjoint
perturbations of semiclassical selfadjoint operators in dimension 2, whose classical flow possesses invariant Lagrangian
tori. This study has been initiated in a previous work of the authors together with San V\~u Ng\d{o}c~\cite{HiSjVu},
where the distribution of eigenvalues coming from the invariant tori satisfying a Diophantine condition had been analyzed.
The original purpose of this paper was to study the contributions to the
spectrum that come from the tori that are rational --- in fact, we shall also consider more general
configurations, with both Diophantine and rational tori occurring.
In this situation, under suitable smallness assumptions on the strength of the non-selfadjoint perturbation,
we show that, in certain rectangles in the complex spectral plane,
the spectrum, up to a small error, agrees with the union of the contributions coming from each of
the invariant tori in question. The quasi-eigenvalues coming from the Diophantine tori can be computed individually,
modulo ${\cal O}(h^{\infty})$, and form a superposition of a finite number of slightly distorted lattices, while the
rational region contributes in a negligible way when compared to the Diophantine one,
already at the level of the counting function asymptotics --- see Theorem 1.1 and the remark following it,
for a precise statement.

In the semiclassical spectral analysis of non-selfadjoint operators in dimension 2, the work~\cite{MeSj2} elucidated
the r\^ole played by certain flow invariant Lagrangian tori in the complexified phase space, and showed
how the latter could be used to obtain complete asymptotic expansions for all eigenvalues in suitable domains in the
complex spectral plane, in the spirit of the classical Bohr-Sommerfeld quantization condition, well known
in dimension one~\cite{HeRo},~\cite{DiSj},~\cite{HaRo}. In the context of small
non-selfadjoint perturbations of selfadjoint operators, the work~\cite{MeSj2} has been pursued in a series of
papers~\cite{HiSj1},~\cite{HiSj2},~\cite{HiSj3a},~\cite{Sj2003}, under the additional assumption that the
classical flow of the unperturbed part should be periodic in a fixed
energy shell. The present work together with its
predecessor~\cite{HiSjVu}, as well as with the aforementioned papers, represents a part of the unified effort aiming
at obtaining a detailed information about the semiclassical behavior of the individual eigenvalues of
small perturbations of selfadjoint operators, arising in applications, such as, say,
the theory of resonances~\cite{SjZw2}, and optimal control and stabilization of linear wave equations~\cite{BurqGerard}.
In collaboration with San V\~u Ng\d{o}c,
the first author is currently investigating the effects of a
small real perturbation on the rational invariant tori in
a completely integrable system,~\cite{HiVu}. Eventually we hope to be
able to attack the problem of the distribution of the imaginary parts of scattering poles
for a strictly convex analytic obstacle in $\real^3$ --- see~\cite{SjZw1} for the results on the distribution of the
real parts of the poles. Indeed, we expect that the results established in the present paper, when combined with those
of~\cite{HiSjVu},~\cite{HiVu}, will prove to be instrumental in achieving this goal.

\subsection{General assumptions}

We shall start by describing the general assumptions on our operators, which will be the same as in~\cite{HiSjVu},
as well as in the earlier papers mentioned above. Let $M$ denote
either the space
$\real^2$ or a real analytic compact manifold of dimension 2. We shall let $\widetilde{M}$ stand for a
complexification of $M$, so that $\widetilde{M}=\comp^2$ in the Euclidean case, and in the compact case, we let
$\widetilde{M}$ be a Grauert tube of $M$ --- see~\cite{GuSt} for the definition and further references.

When $M=\real^2$, let
\begeq
\label{eq00}
P_{\eps}=P^w(x,hD_x,\eps;h),\quad 0<h \leq 1,
\endeq
be the $h$--Weyl quantization on $\real^2$ of a symbol $P(x,\xi,\eps;h)$ (i.e. the Weyl quantization of
$P(x,h \xi,\eps;h)$), depending smoothly on $\eps\in {\rm neigh}(0,\real)$ and taking values in
the space of holomorphic functions of $(x,\xi)$ in a tubular \neigh{}
of $\real^4$ in $\comp^4$, with
\begeq
\label{eq0001}
\abs{P(x,\xi,\eps;h)}\leq {\cal O}(1) m(\Re (x,\xi)),
\endeq
there. Here $m\geq 1$ is an order function on $\real^4$, in the sense that
\begeq
\label{eq0002}
m(X)\leq C_0 \langle{X-Y}\rangle^{N_0} m(Y),\quad X,\,Y\in \real^4,
\endeq
for some $C_0$, $N_0>0$. We shall assume, as we may, that $m$ belongs to its own symbol class, so that
$m\in C^{\infty}(\real^4)$ and
$\partial^{\alpha}m={\cal O}_{\alpha}(m)$ for each $\alpha\in \nat^4$. Then for $h>0$ small enough and when
equipped with the domain $H(m):=\left(m^w(x,hD)\right)^{-1}\left(L^2(\real^2)\right)$, $P_{\eps}$ becomes a closed densely
defined operator on $L^2(\real^2)$.

\medskip
Assume furthermore that
\begeq
\label{eq001}
P(x,\xi,\eps;h)\sim \sum_{j=0}^{\infty} h^j p_{j,\eps}(x,\xi)
\endeq
in the space of holomorphic functions satisfying (\ref{eq0001}) in a
fixed tubular \neigh{} of $\real^4$. We assume that $p_{0,\eps}$ is
elliptic near infinity,
\begeq
\label{eq0011}
\abs{p_{0,\eps}(x,\xi)}\geq \frac{1}{C} m(\Re (x,\xi)),\quad \abs{(x,\xi)}\geq C,
\endeq
for some $C>0$.

When $M$ is a compact manifold, for simplicity we shall take $P_{\eps}$ to be a differential operator on $M$,  such that
for every choice of local coordinates, centered at some point of $M$, it takes the form
\begeq
\label{eq0012}
P_{\eps}=\sum_{\abs{\alpha}\leq m} a_{\alpha,\eps}(x;h)(hD_x)^{\alpha},
\endeq
where $a_{\alpha}(x;h)$ is a smooth function of $\eps\in {\rm neigh}(0,\real)$ with values
in the space of bounded holomorphic functions in a complex \neigh{} of $x=0$. We further assume that
\begeq
\label{eq0013}
a_{\alpha,\eps}(x;h)\sim \sum_{j=0}^{\infty} a_{\alpha,\eps,j}(x) h^j,\quad h\rightarrow 0,
\endeq
in the space of such functions. The semiclassical principal symbol $p_{0,\eps}$, defined on $T^*M$, takes the form
\begeq
\label{eq0014}
p_{0,\eps}(x,\xi)=\sum a_{\alpha,\eps,0}(x)\xi^{\alpha},
\endeq
if $(x,\xi)$ are canonical coordinates on $T^*M$. We make the ellipticity assumption,
\begeq
\label{eq0015}
\abs{p_{0,\eps}(x,\xi)}\geq \frac{1}{C} \langle{\xi}\rangle^m,\quad (x,\xi)\in T^*M,\quad \abs{\xi}\geq C,
\endeq
for some large $C>0$. Here we assume that $M$ has been equipped with some real analytic Riemannian metric so
that $\abs{\xi}$ and $\langle{\xi}\rangle=(1+\abs{\xi}^2)^{1/2}$ are well-defined.

Sometimes, we write $p_{\eps}$ for $p_{0,\eps}$ and simply $p$ for $p_{0,0}$. We make the assumption that
$$
P_{\eps=0}\quad \wrtext{is formally selfadjoint}.
$$
In the case when $M$ is compact, we let the underlying Hilbert space be
$L^2(M, \mu(dx))$ where $\mu(dx)$ is the Riemannian volume element.

The assumptions above imply that the spectrum of $P_{\eps}$ in a fixed \neigh{} of $0\in \comp$ is discrete, when
$0<h\leq h_0$, $0\leq \eps\leq \eps_0$, with $h_0>0$, $\eps_0>0$ sufficiently small. Moreover, if
$z\in {\rm neigh}(0,\comp)$ is an eigenvalue of $P_{\eps}$ then $\Im z ={\cal O}(\eps)$.

We furthermore assume that the real energy surface
$p^{-1}(0)\cap T^*M$ is connected and that
$$
dp\neq 0\quad \wrtext{along}\quad p^{-1}(0)\cap T^*M.
$$

In what follows we shall write
\begeq
\label{eq002}
p_{\eps}=p+i\eps q+{\cal O}(\eps^2),
\endeq
in a \neigh{} of $p^{-1}(0)\cap T^*M$, and for simplicity we shall assume throughout this paper that
$q$ is real valued on the real domain. (In the general case, we should simply replace $q$ below by $\Re q$.)

Let $H_p=p'_{\xi}\cdot \partial_x-p'_x\cdot \partial_{\xi}$ be the Hamilton field of $p$.
In~\cite{HiSjVu}, it was assumed that the energy surface $p^{-1}(0)\cap T^*M$ contains
finitely many $H_p$--invariant analytic Lagrangian tori satisfying a Diophantine condition. Let us recall that according to
a classical theorem of Kolmogorov~\cite{BBGS}, the existence of such tori is assured when $p$ is a small
perturbation of a completely integrable symbol satisfying suitable non-degeneracy assumptions. Since
our primary purpose here is to examine the r\^ole of the rational tori, which are in general destroyed by perturbing
a completely integrable system, throughout this paper we shall work under the assumption that the
$H_p$--flow itself is {completely integrable}. We proceed therefore to discuss the precise assumptions on the
geometry of the energy surface $p^{-1}(0)\cap T^*M$ in this case.

\subsection{Assumptions related to the complete integrability}

As in~\cite{HiSjVu}, let us assume that there exists an analytic real
valued function $f$ on $T^*M$ such that $H_pf=0$, with the
differentials $df$ and $dp$ being linearly independent almost
everywhere. For each $E\in {\rm neigh}(0,\real)$, the level sets
$\Lambda_{a,E}=f^{-1}(a)\cap p^{-1}(E)\cap T^*M$ are invariant under the $H_p$--flow and
form a singular foliation of the 3-dimensional hypersurface
$p^{-1}(E)\cap T^*M$. At each regular point, the leaves of this
foliation are 2-dimensional Lagrangian submanifolds, and each
regular leaf is a finite union of tori. In what follows we shall use
the word ``leaf'' and notation $\Lambda$ for a connected component of
some $\Lambda_{a,E}$. Let $J$ be the set of all leaves in
$p^{-1}(0)\cap T^*M$. Then we have a disjoint union decomposition
\begeq
\label{eq0}
p^{-1}(0)\cap T^*M=\bigcup_{\Lambda \in J} \Lambda,
\endeq
where $\Lambda$ are compact connected $H_p$--invariant sets. The set $J$
has a natural structure of a graph whose edges correspond to families
of regular leaves and the set $S$ of vertices is composed of singular
leaves. The union of edges $J\backslash S$ possesses a natural real
analytic structure and the corresponding tori depend
analytically on $\Lambda\in J\backslash S$ with respect to that
structure.

As in~\cite{HiSjVu}, we shall require $J$ to be a finite connected
graph. We identify each edge of $J$ analytically with a real bounded interval
and this determines a distance on $J$ in the natural way. Assume the
continuity property
\begin{eqnarray}
\label{015}
& & \wrtext{For every}\,\,\,\Lambda_0\in J\,\,\wrtext{and
every}\,\,\eps>0,\, \exists\, \delta>0,\,\,\wrtext{such that if} \\ \nonumber
& & \Lambda\in J,\,\,{\rm
dist}_J(\Lambda,\Lambda_0)<\delta,\,\,\wrtext{then}\,\, \Lambda\subset
\{\rho\in p^{-1}(0);\, {\rm dist}(\rho,\Lambda_0)<\eps\}.
\end{eqnarray}
These assumptions are satisfied, for instance, when $f$ is a
Morse-Bott function restricted to $p^{-1}(0)\cap T^*M$, as in this
case the structure of the singular leaves is known~\cite{San}.

\medskip
Each torus $\Lambda\in J\backslash S$ carries real analytic
coordinates $x_1$, $x_2$ identifying $\Lambda$ with ${\bf T}^2 =
\real^2/2\pi \z^2$, so that along $\Lambda$, we have
\begeq
\label{eq015.5}
H_p=a_1\partial_{x_1}+a_2\partial_{x_2},
\endeq
where $a_1$, $a_2\in \real$. The rotation number is defined as the ratio
$$
\omega(\Lambda)=[a_1:a_2]\in \real {\bf P}^1,
$$
and it depends analytically on $\Lambda\in J\backslash S$. We assume that
$$
\omega(\Lambda)\,\,\wrtext{is not identically constant on any open edge}.
$$
Recall that the leading perturbation $q$ has been introduced in
(\ref{eq002}). For each torus $\Lambda\in J\backslash S$, we define the torus
average $\langle{q}\rangle(\Lambda)$ obtained by integrating $q|_{\Lambda}$ with respect to the natural smooth
measure on $\Lambda$, and assume that the analytic function $J\backslash S \ni \Lambda\mapsto \langle{q}\rangle(\Lambda)$
is not identically constant on any open edge.

\medskip
We introduce
\begeq
\label{eq016}
\langle{q}\rangle_T=\frac{1}{T} \int_{-T/2}^{T/2} q\circ
\exp(tH_p)\,dt,\quad T>0,
\endeq
and consider the compact intervals $Q_{\infty}(\Lambda)\subset \real$, $\Lambda\in J$, defined as in~\cite{HiSjVu},
\begeq
\label{eq0.1}
Q_{\infty}(\Lambda)=\left[\lim_{T \rightarrow \infty} \inf_{\Lambda}
\langle{q}\rangle_T, \lim_{T \rightarrow \infty}
\sup_{\Lambda}\langle{q}\rangle_T\right].
\endeq
Notice that when $\Lambda\in J\backslash S$ and $\omega(\Lambda)\notin {\bf Q}$ then
$Q_{\infty}(\Lambda)=\{\langle{q}\rangle(\Lambda)\}$. In the rational case,
we write $\omega(\Lambda)=\frac{m}{n}$, where $m\in \z$ and $n\in {\bf N}$ are
relatively prime, and where we may assume that $m={\cal
O}(n)$. When $k(\omega(\Lambda)):=\abs{m}+\abs{n}$ is the height of $\omega(\Lambda)$,
we recall from Proposition 7.1 in~\cite{HiSjVu} that
\begeq
\label{eq0.2}
Q_{\infty}(\Lambda)\subset \langle{q}\rangle(\Lambda)+{\cal
O}\left(\frac{1}{k(\omega(\Lambda))^{\infty}}\right)[-1,1].
\endeq

\medskip
\noindent
{\it Remark}. As $J\backslash S\ni \Lambda \rightarrow \Lambda_0\in S$,
the set of all accumulation points of $\langle{q}\rangle(\Lambda)$ is contained in the interval
$Q_{\infty}(\Lambda_0)$. Indeed, when $\Lambda\in J\backslash S$ and $T>0$, there exists
$\rho=\rho_{T,\Lambda}\in \Lambda$ such that
$\langle{q}\rangle(\Lambda)=\langle{q}\rangle_T(\rho)$. Therefore,
each accumulation point of $\langle{q}\rangle(\Lambda)$ as
$\Lambda\rightarrow \Lambda_0\in S$, belongs to
$[\inf_{\Lambda_0}\langle{q}\rangle_T,
\sup_{\Lambda_0}\langle{q}\rangle_T]$. The conclusion follows if we let
$T \rightarrow \infty$.

\medskip
Let $\Lambda_0\in J\backslash S$ be a rational invariant Lagrangian torus, so that as above,
$\omega_0:=\omega(\Lambda_0)=\frac{m}{n}\in {\bf Q}$, $m={\cal
O}(n)$. For future reference, we shall finish this subsection by considering the behavior of
the interval $Q_{\infty}(\Lambda)$ when $\Lambda\neq \Lambda_0$ is a rational torus in a \neigh{} of $\Lambda_0$.
Writing $\omega(\Lambda)=\frac{p}{q}$ where $p\in \z$ and $q\in {\bf
N}$ are relatively prime, $p={\cal O}(q)$, we get, using that $\omega(\Lambda)\neq \omega_0$,
\begeq
\label{eq0.3}
\abs{\omega(\Lambda)-\omega_0}\geq \frac{1}{n q}\geq \frac{1}{n k(\omega(\Lambda))},\
\endeq
and therefore, in view of (\ref{eq0.2}),
\begeq
\label{eq0.4}
Q_{\infty}(\Lambda)\subset \langle{q}\rangle(\Lambda)+{\cal
O}({\rm dist}(\omega(\Lambda),\omega_0)^{\infty})[-1,1].
\endeq
This estimate is uniform in $\omega_0$ provided that we have
a uniform upper bound on the height of the rotation number $\omega_0\in {\bf Q}$.


\subsection{Statement of the main result}

From Theorem 7.6 in~\cite{HiSjVu} we recall that
\begeq
\label{eq0.41}
\frac{1}{\eps}\Im \left({\rm Spec}(P_{\eps})\cap \{z; \abs{\Re z}\leq \delta\}\right)
\subset \left[\inf \bigcup_{\Lambda\in J} Q_{\infty}(\Lambda)-o(1),
\sup \bigcup_{\Lambda\in J} Q_{\infty}(\Lambda)+o(1)\right],
\endeq
as $\eps$, $h$, $\delta\rightarrow 0$. Let us also recall
from~\cite{HiSjVu} that a torus $\Lambda\in J\backslash S$ is said to be
Diophantine if representing $H_p|_{\Lambda} =
a_1\partial_{x_1}+a_2\partial_{x_2}$, as in (\ref{eq015.5}), we have
\begeq
\label{eq0.41.01}
\abs{a\cdot k}\geq \frac{1}{C_0 \abs{k}^{N_0}},\quad 0\neq k\in \z^2,
\endeq
for some fixed $C_0$, $N_0>0$.

Let $F_0\in \cup_{\Lambda\in J} Q_{\infty}(\Lambda)$ be such that
there exist finitely many Lagrangian tori
\begeq
\label{eq0.41.02}
\Lambda_{1,d},\ldots\, ,\Lambda_{L,d}\in J\backslash S
\endeq
that are uniformly Diophantine as in (\ref{eq0.41.01}), and such that
\begeq
\label{eq0.41.1}
\langle{q}\rangle(\Lambda_{j,d})=F_0\quad \wrtext{for}\,\, 1\leq j\leq L,
\endeq
 with
\begeq
\label{eq0.41.2}
d_{\Lambda=\Lambda_{j,d}}\langle{q}\rangle(\Lambda)\neq 0,\quad 1\leq j\leq L.
\endeq

Moreover, assume also that there exist tori $\Lambda_{1,r},\ldots\, ,
\Lambda_{L',r}\in J\backslash S$ with $\omega(\Lambda_{j,r})\in {\bf
Q}$, $1\leq j\leq L'$, and such that the isoenergetic condition
 \begeq
 \label{eq0.41.3}
 d_{\Lambda=\Lambda_{j,r}}\omega(\Lambda)\neq 0
 \endeq
 is satisfied for each $j$, $1\leq j\leq L'$. Assume next that the
 length $\abs{Q_{\infty}(\Lambda_{j,r})}$ of each interval
 $Q_{\infty}(\Lambda_{j,r})$ satisfies
\begeq
\label{eq0.41.3.5}
\abs{Q_{\infty}(\Lambda_{j,r})}>0,\,\,\, j=1,\ldots,\, L',
\endeq
and that
\begeq
\label{eq0.41.4}
 F_0\in Q_{\infty}(\Lambda_{j,r}),\quad 1\leq j\leq L'.
\endeq

We shall assume that
\begeq
\label{eq0.42}
\abs{\langle{q}\rangle(\Lambda_{j,r})-F_0} \geq \frac{1}{{\cal O}(1)},\quad 1\leq j \leq L'.
\endeq

Let us finally make the following global assumption:
\begeq
\label{eq0.44}
F_0\notin \bigcup_{\Lambda \in J \backslash \{\Lambda_{1,d},\ldots, \Lambda_{L,d},\Lambda_{1,r},\ldots, \Lambda_{L',r}\}}
Q_{\infty}(\Lambda).
\endeq
Here we notice that the earlier assumptions imply that $F_0 \notin
Q_{\infty}(\Lambda)$ for $\Lambda_{j,d}\neq \Lambda \in {\rm
neigh}(\Lambda_{j,d},J)$, $1\leq j\leq L$, and $\Lambda_{j,r}\neq \Lambda\in {\rm
neigh}(\Lambda_{j,r},J)$, $1\leq j \leq L'$.

\begin{theo}
Let $F_0\in \cup_{\Lambda\in J} Q_{\infty}(\Lambda)$ be such that the
assumptions {\rm (\ref{eq0.41.1})}, {\rm (\ref{eq0.41.2})}, {\rm
(\ref{eq0.41.3})}, {\rm (\ref{eq0.41.3.5})}, {\rm (\ref{eq0.41.4})},
{\rm (\ref{eq0.42})}, and {\rm (\ref{eq0.44})} are
satisfied. For $1\leq j\leq L$, we fix a basis for the first homology
group of each Diophantine torus $\Lambda_{j,d}$ given by the cycles
$\alpha_{k,j}\subset \Lambda_{j,d}$,
$k=1,2$, and let $S_j\in \real^2$ be the actions and $k_j\in \z^2$ be the Maslov
indices of $\alpha_{k,j}$. Let
\begeq
\label{eq0.44.1}
\kappa_j: {\rm neigh}(\Lambda_{j,d},T^*M)\rightarrow {\rm
neigh}(\xi=0,T^*{\bf T}^2)
\endeq
be a canonical transformation given by the action-angle variables
near $\Lambda_{j,d}$, $1\leq j\leq L$, and such that
$\kappa_j(\alpha_{k,j})=\{x\in {\bf T}^2; x_{3-k}=0\}$, $k=1,2$. Let $\delta>0$ be small and assume that
$$
h\ll \eps\leq h^{\frac{2}{3}+\delta}.
$$
Let $C>0$ be sufficiently large. Then there
exists a bijection $b$ between the spectrum of $P_{\eps}$ in
the rectangle
\begeq
\label{eq0.44.2}
R(F_0,C,\eps,\delta):=\left[-\frac{\eps}{C},\frac{\eps}{C}\right]+
i\eps \left[F_0-\frac{\eps^{\delta}}{C},F_0+\frac{\eps^{\delta}}{C}\right]
\endeq
and the union of two sets of points, $E_d$ and $E_r$, such that
$b(\mu)-\mu={\cal O}(h^{N_0})$. Here $N_0$ is fixed but can be taken
arbitrarily large. The elements of the set $E_d$, $z(j,k)$, are described by Bohr-Sommerfeld type conditions,
\begeq
\label{BS1}
z(j,k)=P_j^{(\infty)}\left(h\left(k-\frac{k_j}{4}\right)-\frac{S_j}{2\pi},\eps;h\right)+{\cal
O}(h^{\infty}),\quad k\in \z^2,\quad 1\leq j\leq L,
\endeq
with precisely one element for each $k\in \z^2$ such that the corresponding $z(j,k)$ belongs to the rectangle
{\rm (\ref{eq0.44.2})}. Here $P_j^{(\infty)}(\xi,\eps;h)$ is
smooth in $\xi\in {\rm neigh}(0,\real^2)$ and $\eps\in {\rm
neigh}(0,\real)$, real-valued for $\eps=0$. We have
\begeq
\label{BS2}
P_j^{(\infty)}(\xi,\eps;h)\sim \sum_{\ell=0}^{\infty} h^{\ell}
p^{(\infty)}_{j,\ell}(\xi,\eps),\quad 1\leq j\leq L,
\endeq
and
\begeq
\label{eqBS3}
p_{j,0}^{(\infty)}(\xi,\eps)=p(\xi)+i\eps\langle{q}\rangle(\xi)+{\cal
O}(\eps^2).
\endeq
Here $p$ and $q$ have been expressed in terms of the action-angle
variables near $\Lambda_{j,d}$ given by $\kappa_j$ in
{\rm (\ref{eq0.44.1})}, and $\langle{q}\rangle$ is the torus average
of $q$ in these coordinates. The cardinality of the set
$E_r$ is
\begeq
\label{eq0.44.3}
{\cal O}\left(\frac{\eps^{3/2}}{h^2}\right).
\endeq
\end{theo}

\medskip
\noindent
{\it Remark}. It follows from Theorem 1.1 that the total number of elements of the set $E_d$ is $\sim
\eps^{1+\delta}/h^2$. Therefore, from (\ref{eq0.44.3}) we see that for $0<\delta<1/2$, the contribution $E_r$
to the spectrum of $P_{\eps}$ in $R(F_0,C,\eps,\delta)$, coming from the rational region, is much weaker than that
of the Diophantine tori $\Lambda_{j,d}$, $1\leq j\leq L$. As will be
seen in the proof, for this result the assumption (\ref{eq0.42}) is important.

\medskip
\noindent
{\it Remark}. Assume that the subprincipal symbol of $P_{\eps=0}$ in (\ref{eq00}) and (\ref{eq0012}) vanishes. Then
it follows from the discussion in the body of the paper that Theorem
1.1 is valid in the larger range
\begeq
\label{eq0.44.4}
h^2\ll \eps={\cal O}(h^{\frac{2}{3}+\delta}),\quad 0<\delta\ll 1.
\endeq

Theorem 1.1 can be viewed as a partial generalization of one of the main results
of~\cite{HiSjVu}, where energy levels corresponding only to
Diophantine tori have been consi\-de\-red. In that paper, instead of the upper bound
$\eps\leq h^{2/3+\delta}$, it was merely required that $\eps={\cal
O}(h^{\delta})$ for some small fixed $\delta>0$. (Also the lower bounds there
were considerably weaker than $h\ll \eps$.) As will be seen in the
proof, here the strengthened upper bound on $\eps$ is
required in order to compensate for the exponential growth of the
resolvent of $P_{\eps}$ in the rational region, when considering the tunnel effect between the
Diophantine and the rational tori --- see also the discussion in the
next section.

\medskip
In the case when there are no Diophantine tori corresponding to the
energy level $(0,\eps F_0)\in \comp$,
the result of Theorem 1.1 can be improved in two ways: we can put $\delta=0$ in (\ref{eq0.44.2}), and also, the
upper bound $\eps\leq h^{2/3+\delta}$ can be replaced by $\eps={\cal
O}(h^{\widetilde{\delta}})$, $\widetilde{\delta}>0$.

\begin{theo}
Let us keep all the assumptions of Theorem {\rm 1.1}, and assume that
$L=0$ in {\rm (\ref{eq0.41.02})}. Assume furthermore that $\eps={\cal O}(h^{\widetilde{\delta}})$, $\widetilde{\delta}>0$, satisfies
$\eps \gg h$.
There exists a constant $C>0$ such that the number of eigenvalues of $P_{\eps}$ in the rectangle
\begeq
\label{eq0.45}
\abs{\Re z} <\frac{\eps}{C},\quad \abs{\frac{\Im z}{\eps}-F_0}<\frac{1}{C}
\endeq
does not exceed
\begeq
\label{eq0.46}
{\cal O}\left(\frac{\eps^{3/2}}{h^2}\right).
\endeq
\end{theo}

\medskip
\noindent
{\it Remark}. As will be explained in the beginning of section 3, the isoenergetic assumption (\ref{eq0.41.3})
implies that associated with each rational torus $\Lambda_{j,r}$,
$1\leq j\leq L'$, there is an analytic family of rational Lagrangian tori
$\Lambda_{E,j,r}\subset p^{-1}(E)\cap T^*M$ for $E\in {\rm neigh}(0,\real)$, depending analytically on $E$, and with
$\Lambda_{E=0,j,r}=\Lambda_{j,r}$, $1\leq j\leq L'$. Theorem 1.2 can therefore be interpreted as saying that that only an
$\eps^{1/2}$-\neigh{} of the set
\begeq
\label{eq0.46.4}
\bigcup_{j=1}^{L'}\bigcup_{E={\cal O}(\eps)} \Lambda_{E,j,r}
\endeq
contributes to the spectrum in the region (\ref{eq0.45}).

\medskip
\noindent
For notational simplicity only, when proving Theorem 1.1 and Theorem 1.2, we shall assume
that $L=2$, $L'=1$, and that $\Lambda_{1,d}$, $\Lambda_{2,d}$, and
$\Lambda_{1,r}$ all belong to the same open edge of $J$, so that, when identifying the edge with a real bounded
interval, we have
\begeq
\label{eq0.46.44}
\Lambda_{1,d}<\Lambda_{1,r}<\Lambda_{2,d}.
\endeq

The structure of the paper is as follows. In section 2 we present a
general outline of the proof of Theorem 1.1. Section 3 is devoted
to a formal microlocal Birkhoff normal form construction for $P_{\eps}$ near $\Lambda_{1,r}$, and in section 4 the
formal argument of the previous section is justified by constructing a microlocal Hilbert space in a full \neigh{} of
$\Lambda_{1,r}$, realizing the normal form reduction there. In the beginning of section 5 we construct the
global Hilbert space where we study our operator $P_{\eps}$, and introduce two reference operators, associated with the
Diophantine and the rational regions, respectively. Section 5 is concluded by constructing the
resolvent for $P_{\eps}$ globally, and we obtain Theorem 1.1 by comparing the spectral projections of $P_{\eps}$ and
of the reference operators. In section 6, we apply Theorem 1.1 to a small complex perturbation of
the semiclassical Laplacian on a convex analytic surface of
revolution, and give a partial ge\-ne\-ra\-lization of the corresponding discussion in~\cite{HiSjVu}.
The appendix contains a proof of a simple trace class estimate for the Toeplitz ope\-rator with a
compactly supported smooth symbol, acting on a weighted $L^2$--space of holomorphic functions.
This estimate, which seems to be of an independent interest, is used in section 5 in the main text.

\medskip
\noindent
{\bf Acknowledgment}. This project began when the first author was visiting \'Ecole Polytechnique in September of
2005. It is a pleasure for him to thank its Centre de Math\'ematiques for a generous hospitality. We are also
grateful to San V\~u Ng\d{o}c for many interesting discussions around
this work and for making a written contribution, which is
planned to be used in a future work of S. V\~u
Ng\d{o}c and the first author. The research of the first author is supported in part by
the National Science Foundation under grant DMS--0304970 and the
Alfred P. Sloan Research Fellowship.

\section{Outline of the proof}
\setcounter{equation}{0}

The purpose of this section is to provide a broad outline of the proof
of Theorem 1.1. Compared with the previous work~\cite{HiSjVu},
addressing only the case of Diophantine tori, here the essential new difficulties
will be concerned with the analysis in the rational region. We shall begin by presenting an outline of the argument
in this case.

Working microlocally in the rational region and introducing action-angle variables in a \neigh{} of $\Lambda_{1,r}\simeq {\bf T}^2$, we are led to consider an
operator, defined microlocally near $\xi=0$ in $T^*{\bf T}^2_{x,\xi}$, with the leading
symbol given by
\begeq
\label{out1}
p_{\eps}(x,\xi)=p(\xi)+i\eps q(x,\xi)+{\cal O}(\eps^2),\quad p(\xi)=\omega\cdot \xi+{\cal O}(\xi^2).
\endeq
Here $\omega=(k,l)\in \z^2$ and to fix the ideas, let us restrict the attention to the model case where $\omega=(0,1)$
and the ${\cal O}(\xi^2)$--term in (\ref{out1}) reduces to $\xi_1^2$ --- this choice of the nonlinearity in $p$
is in agreement with the isoenergetic condition (\ref{eq0.41.3}). Following the general ideas of a Birkhoff normal form
construction, we would like to eliminate, as much as possible, the $x$-dependence in the symbol in (\ref{out1}).
Performing first successive averaging procedures
along the closed orbits of the $H_p$--flow comprising the rational tori $\Lambda_{E,1,r}$,
$E\in {\rm neigh}(0,\real)$, we achieve that the leading symbol in (\ref{out1}) becomes
\begeq
\label{out2}
\widetilde{p}_{\eps}(x,\xi)=\xi_2+\xi_1^2+{\cal O}(\eps)+{\cal O}((\eps,\xi_1)^{\infty}),
\endeq
where the ${\cal O}(\eps)$--term is independent of $x_2$. In the terminology of classical mechanics, this initial
reduction is based on a secular perturbation theory --- see~\cite{LcLb}. Carrying out the reduction on the operator level,
we obtain an operator
of the form $\widetilde{P}_{\eps}=\widetilde{P}_{\eps}(x_1, hD_{x_1},hD_{x_2};h)$, which may also be viewed
as a family of one-dimensional non-selfadjoint operators acting in $x_1$, with a leading symbol of the
form
$$
hk+\xi_1^2+{\cal O}(\eps),\quad k\in \z.
$$
For this family, we cannot ex\-clude the occurrence of a pseudo\-spectral
phe\-no\-me\-non \cite{DeSjZw}, leading to the exponential growth of the resolvent norms in the
spectral regions of interest. This makes it difficult to exploit the secular perturbation theory and to
simplify the operator further.

\medskip
Nevertheless, in section 3 we show that working in a region where
\begeq
\label{out3}
\abs{\xi_1}\gg \eps^{1/2},
\endeq
and so away from an $\eps^{1/2}$--\neigh{} of the set
\begeq
\label{out3.5}
\bigcup_{\abs{E}< \delta_0} \Lambda_{E,1,r},\quad 0 < \delta \ll 1,
\endeq
the $x_1$-dependence in the symbol (\ref{out2}) can be eliminated completely, and in particular, here
the leading perturbation $q$ in (\ref{out1}) becomes replaced by its torus average. When approaching the region where
$\xi_1={\cal O}(\eps^{1/2})$, the normal form construction breaks down and no additional simplification of
the operator $P_{\eps}$ is obtained.

\medskip
To implement the complete reduction in the region (\ref{out3}) requires an introduction of a microlocal
Hilbert space of functions in a sufficiently small but fixed \neigh{}
of $\Lambda_{1,r}$. Because of the degeneration of the normal form construction very close to the
rational torus, when defining the Hilbert space
in a full \neigh{} of $\Lambda_{1,r}$, it becomes convenient and indeed, natural, to perform a second
microlocalization --- in this case, it amounts to considering our operators in $\widetilde{h}=h/\sqrt{\eps}$--quantization with respect to the
$x_1$--variable and performing an $\widetilde{h}$--Bargmann transformation in $x_1$. In section 4 we show that,
on the transform side, the microlocal Hilbert space in question becomes a
well-defined weighted space of holomorphic functions in a region $(\Re x_1,\Re x_2)\in {\bf T}^2$,
$\abs{\Im x_1}\ll \frac{1}{\sqrt{\eps}}$, $\abs{\Im x_2}\ll 1$, with the corresponding strictly
plurisubharmonic weight being uniformly well behaved and close to the standard quadratic one --- see Proposition 4.1 for
the precise statement and also, the discussion in section 4.4.

\medskip
The idea now is to use the assumption (\ref{eq0.42}) to show that
$(1/\eps)\left(P_{\eps}-z\right)$ becomes elliptic (viewed as an $\widetilde{h}$-pseudodifferential operator) near
$\Lambda_{1,r}$, when away from an ${\cal O}(\eps^{1/2})$--\neigh{} of the set in (\ref{eq0.46.4}), while the
invertibility away from the tori $\Lambda_{1,r}\cup\Lambda_{1,d}\cup\Lambda_{2,d}$ should
follow from (\ref{eq0.44}). Here the spectral parameter $z$ varies
in the domain (\ref{eq0.45}).

To handle the remaining phase space region near $\Lambda_{1,r}$, in section 5
we construct a trace class perturbation $K$, whose trace class norm does not exceed
\begeq
\label{out4}
{\cal O}\left(\frac{\eps^{3/2}}{h^{2}}\right),
\endeq
such that if $P_d:=P_{\eps}+i\eps K$  then $P_d-z$ becomes invertible, when away from the Diophantine quasi-eigenvalues
$z(j,k)$ in (\ref{BS1}). Moreover, we obtain a sufficiently good control on the norm of the inverse of $P_d$, when
the latter is considered in a global Hilbert space, obtained by gluing together the microlocal Hilbert space near $\Lambda_{1,r}$ and
the space away from $\Lambda_{1,r}$, defined using the Diophantine analysis of ~\cite{HiSjVu}.
The trace class perturbation $K$ is constructed as a Toeplitz operator on
the FBI--Bargmann transform side and when deriving the trace class norm bound (\ref{out4}), we use a general estimate
of Proposition A.1 in the appendix.

\medskip
It will be fruitful to think of the Diophantine and the rational tori in question as of microlocal wells to which the main
difficulties of our problem are localized. From this point of view we may think of the operator $P_d$ as a
reference operator associated to the Diophantine region. Proceeding in the spirit of tunneling problems,
in section 5 we next define and study a reference operator associated to the rational region, $P_r$,
obtained by modifying $P_{\eps}$ away from the rational region and
such that $P_r-z$ is invertible outside of a small
\neigh{} of $\Lambda_{1,r}$. Because of the pseudospectral difficulties in the normal form construction for
$P_{\eps}$ in an ${\cal O}(\eps^{1/2})$--\neigh{} of $\Lambda_{1,r}$, when estimating the resolvent of $P_r$, we are only able to
show that it enjoys an exponential upper bound, with the exponent there being given, roughly speaking,
by the phase space volume of the region near the rational torus, not covered by the normal form, multiplied by $h^{-2}$,
or, equivalently, by the trace class norm of the perturbation $K$ in (\ref{out4}).

\medskip
Using the operators $P_d$ and $P_r$, together with an additional reference ope\-ra\-tor corresponding
to the elliptic region, we next construct and study an approximate, and then exact,
resolvent of $P_{\eps}$. To obtain the main result of Theorem 1.1 we
would like to compare the spectral
projections of $P_{\eps}$ with those of the reference operators. Due to the exponential growth of
the resolvent of $P_{\eps}$ near $\Lambda_{1,r}$, at this point
it becomes very important to estimate the tunnel effect between the
Diophantine and the rational tori and to
show that it is small enough to overrule the pseudospectral growth of the resolvent in the rational region. This tunneling
analysis is carried out at the end of section 5 and it involves an
additional modification of phase space exponential weights
near the invariant tori. Imposing the upper bound $\eps ={\cal
O}\left(h^{2/3+\delta}\right)$, $0<\delta\ll 1$, assures that
our perturbative argument goes through, and we can conclude the proof by comparing the
spectral projections, as indicated above.

\vskip 2mm
\noindent
{\it Remark}. The idea of using auxiliary trace class perturbations to create a gap
in the spectrum of a non-selfadjoint operator has a long tradition in abstract
non-selfadjoint spectral theory and seems to go back to the work of
Markus and Matsaev~\cite{MaMa}, see also~\cite{Ma}. It has been used by the second author in the
theory of resonances~\cite{Sj97},~\cite{Sj01}, and when studying
spectral asymptotics for damped wave equations on compact domains~\cite{Sj00} (see also~\cite{Hi03}).
In the present paper, in the absence of the Diophantine tori, once the trace class perturbation $K$, alluded to above,
has been constructed, we can conclude the proof of Theorem 1.2, in section 5,
by relying upon some standard Fredholm determinant estimates~\cite{GoKr}.

\section{The normal form construction near $\Lambda_{1,r}$}
\setcounter{equation}{0}
For simplicity, we shall concentrate throughout the following
discussion on the case when $M=\real^2$, the
compact real analytic case being analogous --- see also the appendix
in~\cite{HiSj1} for the basic facts about FBI transforms on manifolds.
We shall keep all the assumptions made in the introduction, and consider an operator $P_{\eps}$ in (\ref{eq00}) with a
principal symbol
\begeq
\label{eq0.47}
p_{\eps}=p+i\eps q+{\cal O}(\eps^2),
\endeq
in a \neigh{} of $p^{-1}(0)\cap \real^4$. In order to simplify the
presentation, we shall furthermore assume that the order function $m$ introduced in
(\ref{eq0001}) belongs to $L^{\infty}(\real^{4})$. It will be clear that the analysis below extends
to the case of a general order function $m\geq 1$. From the
introduction, let us also recall the simplifying assumption that $L'=1$ so that
$\Lambda_{1,r}$ is the only rational torus corresponding to the level
$(0,\eps F_0)$.

\medskip
In this section, we shall work microlocally near $\Lambda_{1,r}\subset p^{-1}(0)\cap \real^4$. Let
\begeq
\label{eq0.5}
\kappa_0: {\rm neigh}(\Lambda_{1,r},\real^4)\rightarrow {\rm neigh}(\xi=0,T^*{\bf T}^2),
\endeq
be a real and analytic canonical transformation, given by the action-angle variables, and such that
$\kappa_0(\Lambda_{1,r})$ is the zero section in $T^*{\bf T}^2$. Then $p\circ \kappa_0^{-1}$ is a function of $\xi$ only, and
to simplify the notation we shall write $p\circ \kappa^{-1}_{0}=p(\xi)$. We have $p(0)=0$ and without loss of generality
we may assume that
\begeq
\label{eq0.6}
\partial_{\xi_1}p(0)=0,\quad \partial_{\xi_2}p(0)>0.
\endeq
The isoenergetic assumption (\ref{eq0.41.3}) takes the
following form,
\begeq
\label{eq1}
\partial_{\xi_1}^2 p(0)\neq 0.
\endeq
In order to fix the ideas, we assume that $\partial_{\xi_1}^2 p(0)>0$.

By the implicit function theorem, the equation
$\partial_{\xi_1}p(\xi)=0$ has a unique analytic local solution
$\xi_1=f(\xi_2)$ with $f(0)=0$. The function $\xi_2\mapsto
p(f(\xi_2),\xi_2)$ has a positive derivative near $0$, and therefore the
equation $p(f(\xi_2),\xi_2)=E$ has a unique solution $\xi_2(E)$ close
to $0$ for $E\in {\rm neigh}(0,\real)$. We obtain a family of
rational Lagrangian tori $\Lambda_{E,1,r}\subset p^{-1}(E)$, defined by
\begeq
\label{eq1.2}
\xi_2=\xi_2(E),\quad \xi_1=f(\xi_2(E)).
\endeq
By construction,
$\partial_{\xi_1}p=0$ on $\Lambda_{E,1,r}$, and hence,
\begeq
\label{eq1.5}
\partial_{\xi_1}p(\xi_1,\xi_2)={\cal O}(\xi_1-f(\xi_2)),
\partial_{\xi_2} p(\xi_1,\xi_2)=\partial_{\xi_2}
p(f(\xi_2),\xi_2)+{\cal O}((\xi_1-f(\xi_2))).
\endeq

Implementing $\kappa_0$ in (\ref{eq0.5})
by means of a microlocally unitary Fourier integral operator with a
real phase as in Theorem 2.4 in~\cite{HiSj1},
and conjugating $P_{\eps}$ by this operator, we obtain a new $h$-pseudodifferential operator, still denoted by $P_{\eps}$, defined microlocally
near $\xi=0$ in $T^*{\bf T}^2$. The full symbol of $P_{\eps}$ is holomorphic in a fixed complex \neigh{} of $\xi=0$,
and the leading symbol is given by
\begeq
\label{eq2}
p_{\eps}(x,\xi)=p(\xi)+i\eps q(x,\xi)+{\cal O}(\eps^2),
\endeq
with
\begeq
\label{eq2.1}
p(\xi_1,\xi_2)=p(f(\xi_2),\xi_2)+g(\xi_1,\xi_2)(\xi_1-f(\xi_2))^2.
\endeq
Here $g(0)>0$, since we have assumed that $\partial_{\xi_1}^2 p(0)>0$, and the function $q$ in (\ref{eq2})
is real on the real domain. On the operator level, $P_{\eps}$ acts in the space of microlocally defined Floquet periodic functions on ${\bf T}^2$,
$L^2_{\theta}({\bf T}^2)\subset L^2_{{\rm loc}}(\real^2)$, elements $u$ of which satisfy
\begeq
\label{2.2}
u(x-\nu)=e^{i\theta\cdot \nu} u(x),\quad \theta=\frac{S}{2\pi h}+\frac{k_0}{4},\quad \nu\in 2\pi \z^2.
\endeq
Here $S=(S_1,S_2)$ is given by the classical actions,
$$
S_j=\int_{\alpha_j}\eta\,dy,\quad j=1,2,
$$
with $\alpha_j$ forming a system of fundamental cycles in
$\Lambda_{1,r}$, such that
$$
\kappa_0(\alpha_j)=\beta_j,\quad j=1,2, \quad \beta_j=\{x\in {\bf T}^2;\, x_{3-j}=0\}.
$$
The tuple $k_0=(k_0(\alpha_1),k_0(\alpha_2))\in \z^2$ stands for the Maslov indices of the cycles $\alpha_j$, $j=1,2$.

As a first step in the normal form construction for $P_{\eps}$, we shall apply the secular perturbation theory
to the principal symbol $p_{\eps}$ in (\ref{eq2}) --- see also~\cite{LcLb}.

Let
\begeq
\label{eq3}
\langle{q}\rangle_2(x_1,\xi)=\frac{1}{2\pi} \int_0^{2\pi}
q(x,\xi)\,dx_2
\endeq
denote the average of $q$ with respect to $x_2$. Using the assumption
(\ref{eq0.6}) and proceeding as in section 4 of~\cite{HiSj1}
(see also section 2 of~\cite{HiSjVu}),
it is straightforward to construct, by successive averagings in $x_2$, a symbol $G_1(x,\xi)=G_1^{(N)} (x,\xi)$, analytic
in $(x,\xi)$, such that
\begeq
\label{eq4}
H_p G_1 = q - \widetilde{\langle{q}\rangle_2}(x_1,\xi)+{\cal
O}((\xi_1-f(\xi_2))^N),
\endeq
where
$\widetilde{\langle{q}\rangle_2}(x_1,\xi)=\langle{q}\rangle_2(x_1,\xi)+{\cal
O}(\xi_1-f(\xi_2))$ is independent of $x_2$. Here $N\in \nat$ can be taken arbitrarily large but fixed. We get from
(\ref{eq4}), by a Taylor expansion,
\begin{eqnarray*}
p_{\eps}\left(\exp(i\eps H_{G_1})(x,\xi)\right) & = & p(\xi)+i\eps
\widetilde{\langle{q}\rangle_2}(x_1,\xi)+{\cal O}\left(\eps^2+\eps
\left(\xi_1-f(\xi_2)\right)^N\right)\\
& = & p(\xi)+i\eps \widetilde{\langle{q}\rangle_2}(x_1,\xi)+i\eps^2
\widetilde{q}+{\cal O}\left(\eps^3+\eps
\left(\xi_1-f(\xi_2)\right)^N\right),
\end{eqnarray*}
where $\widetilde{q}=\widetilde{q}(x,\xi)$. We next construct $G_2$, analytic in $(x,\xi)$ and such that
$$
H_p
G_2=\widetilde{q}-\widetilde{\langle{\widetilde{q}\rangle}_2}(x_1,\xi)+{\cal
O}((\xi_1-f(\xi_2))^N).
$$
Then
\begin{eqnarray*}
& & p_{\eps}\left(\exp(i\eps H_{G_1})(\exp(i\eps^2
H_{G_2})(x,\xi))\right) \\
& = & p(\xi)+i\eps\widetilde{\langle{q}\rangle_2}+i\eps^2
\widetilde{\langle{\widetilde{q}\rangle}_2}+
{\cal O}\left(\eps^3+\eps\left(\xi_1-f(\xi_2)\right)^N\right).
\end{eqnarray*}
It is clear that this procedure can be iterated, and after $N$ steps, we define
\begeq
\label{eq5}
\kappa_{\eps}:=\exp(i\eps H_{G_1})\circ \exp
(i\eps^{2}H_{G_{2}})\circ\ldots \circ \exp(i\eps^N H_{G_N}).
\endeq
It follows that
\begin{eqnarray}
\label{eq6}
& & p_{\eps}\left(\kappa_{\eps}(x,\xi)\right) =  p(\xi)+i\eps
\widetilde{\langle{q}\rangle_2}(x_1,\xi)+\eps^2 r_{\eps}(x_1,\xi) \\ \nonumber
& + & {\cal
O}\left(\eps^{N+1}+\eps\left(\xi_1-f(\xi_2)\right)^N\right)=
p'_{\eps}(x_1,\xi)+{\cal O}\left(\eps^{N+1}+\eps\left(\xi_1-f(\xi_2)\right)^N\right).
\end{eqnarray}
Here the last equality defines $p'_{\eps}(x_1,\xi)$.

Using the same averaging procedure as above also on the level of lower order symbols, as in section 4 of ~\cite{HiSj1}
and section 3 of ~\cite{HiSjVu},
we conclude that there exists an analytic elliptic Fourier integral operator $F=F^{(N)}_{\eps}$ in the complex domain,
quantizing the holomorphic canonical transformation $\kappa_{\eps}$ in
(\ref{eq5}), such that
\begeq
\label{eq7}
F^{-1} P_{\eps} F = P'_{\eps}(x_1,hD_x;h)+R_{\eps}(x,hD_x;h).
\endeq
Here the full symbol of $P'_{\eps}$ is independent of $x_2$ and
\begeq
\label{eq7.1}
R_{\eps}(x,\xi;h)={\cal
O}\left(\eps^{N+1}+\eps
\left(\xi_1-f(\xi_2)\right)^N+h^{N+1}\right).
\endeq
The leading symbol of $P'_{\eps}(x_1,hD_x;h)$ is $p'_{\eps}(x_1,\xi)$
in (\ref{eq6}). As in section 6 of~\cite{HiSj1} and section 2 of~\cite{MeSj2},
the operator $F$ is defined by working on the FBI--Bargmann transform side.

\medskip
When discussing further reductions of $P'_{\eps}$, it is natural to exploit the fact that this operator
is independent of $x_2$, and hence, at least formally, by taking a Fourier series expansion in $x_2$, we can reduce
the study of $P'_{\eps}$ to the study of a family of one-dimensional operators $P'_{\eps}(x_1,hD_{x_1},\xi_2;h)$, with
\begeq
\label{eq7.5}
\xi_2=h\left(k-\frac{k_0(\alpha_2)}{4}\right)-\frac{S_2}{2\pi}\in {\rm neigh}(0,\real),\quad k\in \z.
\endeq
The family $P'_{\eps}(x_1,hD_{x_1},\xi_2,h)$ acts on the microlocal space of Floquet periodic functions
$L^2_{\theta_1}({\bf T}^1)$, ${\bf T}^1=\real/2\pi \z$,
$\theta_1=S_1/2\pi h+k_0(\alpha_1)/4$, defined similarly to (\ref{2.2}).
We would like to eliminate the $x_1$--dependence in the symbol of
$P'_{\eps}(x_1,hD_{x_1},\xi_2;h)$ by
means of an additional conjugation by an elliptic Fourier integral
operator. Using (\ref{eq2.1}) and (\ref{eq6}) we get
\begin{eqnarray}
\label{eq8}
& & P'_{\eps}(x_1,hD_{x_1},\xi_2;h) =
p(f(\xi_2),\xi_2)+g(hD_{x_1},\xi_2)\left(hD_{x_1}-f(\xi_2)\right)^2
\\ \nonumber & & +i\eps
\widetilde{\langle{q}\rangle_2}(x_1,hD_{x_1},\xi_2)+
\eps^2 r_{\eps}(x_1,hD_{x_1},\xi_2)+{\cal O}(h)+{\cal O}(h^2).
\end{eqnarray}
Let us recall that $g(0)>0$, and the ${\cal O}(h)$--contribution in (\ref{eq8}) is the subprincipal term
in the full symbol of $P'_{\eps}$. After a conjugation by $\exp(\frac{i}{h} f(\xi_2)x_1)$, modifying the Floquet
condition on ${\bf T}^1$, we get
\begin{eqnarray}
\label{eq9}
& & e^{-\frac{i}{h} f(\xi_2)x_1}P'_{\eps}\left(x_1,hD_{x_1},\xi_2;h\right)e^{\frac{i}{h} f(\xi_2)x_1}
= p(f(\xi_2),\xi_2) \\ \nonumber
& & +g\left(f(\xi_2)+hD_{x_1},\xi_2\right)(hD_{x_1})^2 \\ \nonumber
& & + \left(i\eps \widetilde{\langle{q}\rangle_2}+\eps^2 r_{\eps}+{\cal
O}(h)+{\cal O}(h^2)\right)(x_1,f(\xi_2)+hD_{x_1},\xi_2).
\end{eqnarray}

In section 4 of~\cite{HiSj1}, it is explained how to eliminate the
$x_1$--dependence in (\ref{eq9}) by means of a Fourier integral operator conjugation. Here we shall
follow the procedure there after a suitable change of Planck's constant. Let us work microlocally in a region
\begeq
\label{eq10}
\abs{\xi_1}\sim \mu,\quad (\varepsilon+h)^{1/2} \ll \mu \ll 1.
\endeq
We write
$$
hD_{x_1}=\mu \widetilde{h}D_{x_1},\quad \widetilde{h}=\frac{h}{\mu}\ll
1.
$$
If $\xi_1$, $\widetilde{\xi_1}$ denote the cotangent variables
corresponding to $hD_{x_1}$ and $\widetilde{h}D_{x_1}$, respectively,
we have
$$
\xi_1=\mu \widetilde{\xi_1}.
$$
Then (\ref{eq9}) gives
\begin{eqnarray}
\label{eq11}
& & \mu^{-2} e^{-\frac{i}{h} f(\xi_2)x_1}P'_{\eps}(x_1, hD_{x_1}, \xi_2;h) e^{\frac{i}{h}
f(\xi_2)x_1} \\ \nonumber
& = & \frac{1}{\mu^2}
p(f(\xi_2),\xi_2)+g(f(\xi_2)+\mu\widetilde{h}D_{x_1},\xi_2)\left(\widetilde{h}D_{x_1}\right)^2
\\ \nonumber
& & +
\left(\frac{\eps}{\mu^2}i\widetilde{\langle{q}\rangle_2}+\frac{{\cal
O}(h)}{\mu^2}+\frac{\eps^2}{\mu^2}r_{\eps}\right)(x_1,f(\xi_2)+\mu
\widetilde{h}D_{x_1},\xi_2) \\ \nonumber
& & + {\cal O}(\widetilde{h}^2)\left(x_1,f(\xi_2)+\mu
\widetilde{h}D_{x_1},\xi_2\right),
\end{eqnarray}
which can be viewed as an $\widetilde{h}$--pseudodifferential operator. The symbol
associated to the second term in the right hand side of (\ref{eq11}) is then
\begeq
\label{eq11.5}
g(f(\xi_2)+\mu \widetilde{\xi_1},\xi_2)\widetilde{\xi_1}^2,
\endeq
and it follows from (\ref{eq10}) that we work in a region where
\begeq
\label{eq12}
\abs{\widetilde{\xi_1}}\sim 1.
\endeq
Notice that in this region, the $\widetilde{\xi}_1$--gradient of (\ref{eq11.5}) is of the order of magnitude 1.

\medskip
\noindent
We set next
\begeq
\label{eq12.2}
r_0\left(x_1,\widetilde{\xi_1},\frac{\eps}{\mu^2},\frac{h}{\mu^2},\xi_2\right)=g(f(\xi_2)+\mu
\widetilde{\xi_1},\xi_2)\widetilde{\xi_1}^2+{\cal
O}\left(\frac{\eps+h}{\mu^2}\right),
\endeq
where the ${\cal O}\left(\frac{\eps+h}{\mu^2}\right)$--term stands for
the third term in the right hand side of (\ref{eq11}). Following the
argument of section 4 of~\cite{HiSj1}, we shall now recall how the
$x_1$--dependence in $r_0$ can be eliminated by means of a
suitable canonical transformation.

\medskip
We look for
$\varphi_0=\varphi_0\left(x_1,\widetilde{\xi_1},\frac{\eps}{\mu^2},\frac{h}{\mu^2},\xi_2\right)$,
such that
\begeq
\label{eq12.5}
r_0\left(x_1,\widetilde{\xi_1}+\partial_{x_1}\varphi_0,\frac{\eps}{\mu^2},\frac{h}{\mu^2},\xi_2\right)=
\biggl\langle{r_0\left(\cdot,\widetilde{\xi_1},\frac{\eps}{\mu^2},\frac{h}{\mu^2},\xi_2\right)\biggr\rangle}_1.
\endeq
Here, for a smooth function $f(x,\xi)$ defined near $\xi=0$ in $T^*{\bf T}^2$,
the expression $\langle{f}\rangle_1$ stands for the average with respect to $x_1$,
$$
\langle{f}\rangle_1(x_2,\xi)=\frac{1}{2\pi}\int_0^{2\pi} f(x,\xi)\,dx_1.
$$

By the implicit function theorem, (\ref{eq12.5}) has an analytic solution with
$\partial_{x_1}\varphi_0$ single-valued and ${\cal
O}((\eps+h)/\mu^2)$. Taking a Taylor expansion of (\ref{eq12.5}) and
using (\ref{eq12.2}), we get
$$
\left(\partial_{\widetilde{\xi_1}}r_0\right)\left(x_1,\widetilde{\xi_1},\frac{\eps}{\mu^2},\frac{h}{\mu^2},\xi_2\right)
\partial_{x_1}\varphi_0 + \left(r_0-\langle{r_0}\rangle_1\right)={\cal
O}\left(\left(\frac{\eps+h}{\mu^2}\right)^2\right),
$$
and using also that the $\widetilde{\xi_1}$--gradient of
(\ref{eq11.5}) is $\sim 1$, we conclude that
$$
\varphi_0=\varphi_{\mu}+x_1\widetilde{\zeta_1},
$$
where
$$
\widetilde{\zeta_1}=\widetilde{\zeta_1}\left(\widetilde{\xi_1},\frac{\eps}{\mu^2},\frac{h}{\mu^2},\xi_2\right)=
{\cal
O}\left(\left(\frac{\eps+h}{\mu^2}\right)^2\right),
$$
and $\varphi_{\mu}={\cal O}((\eps+h)/\mu^2)$ is periodic in $x_1$. We set
$\widetilde{\eta_1}=\widetilde{\xi_1}+\widetilde{\zeta_1}$, and view
$\varphi_{\mu}$ as a function of $\widetilde{\eta_1}$ rather than
$\widetilde{\xi_1}$.

Summarizing the discussion above, we see that there exists a holomorphic phase function
\begeq
\label{eq13}
\varphi_{\mu}(x_1,\widetilde{\eta_1})=\varphi_{\mu}\left(x_1,\widetilde{\eta_1},\frac{\eps}{\mu^2},
\frac{h}{\mu^2},\xi_2\right)={\cal
O}\left(\frac{\eps+h}{\mu^2}\right)
\endeq
defined in a fixed complex \neigh{} of $x_1\in {\bf T}^1$,
$\abs{\widetilde{\eta_1}}\sim 1$, such that if
$$
\psi(x_1,\widetilde{\eta_1})=x_1\widetilde{\eta_1}+\varphi_{\mu}(x_1,\widetilde{\eta_1}),
$$
then the canonical transformation
\begeq
\label{eq13.5}
\kappa_{\mu,\eps,h}:(y_1,\widetilde{\eta_1})=(\psi'_{\widetilde{\eta_1}},\widetilde{\eta_1})\mapsto
(x_1,\psi'_{x_1})=(x_1,\widetilde{\xi_1})
\endeq
is ${\cal O}\left(\frac{\eps+h}{\mu^2}\right)$--close to the identity,
and
$$
\left(g(f(\xi_2)+\mu
\widetilde{\xi_1},\xi_2)\widetilde{\xi_1}^2+\left(\frac{\eps}{\mu^2}i\widetilde{\langle{q}\rangle_2}+
\frac{{\cal
O}(h)}{\mu^2}+\frac{\eps^2}{\mu^2}r_{\eps}\right)(x_1,f(\xi_2)+\mu\widetilde{\xi_1},\xi_2)\right)\circ
\kappa_{\mu,\eps,h}
$$
is independent of $y_1$ and is equal to
\begin{eqnarray}
\label{eq14}
& & g(f(\xi_2)+\mu \widetilde{\eta_1})\widetilde{\eta_1}^2 \\ \nonumber
& + & \left(\frac{\eps}{\mu^2}
i\langle{\widetilde{\langle{q}\rangle_2}\rangle}_1+\frac{{\cal
O}(h)}{\mu^2}+\frac{\eps^2}{\mu^2}\langle{r_{\eps}}\rangle_1\right)
(f(\xi_2)+\mu \widetilde{\eta_1},\xi_2)+ {\cal O}\left(\left(\frac{\eps+h}{\mu^2}\right)^2\right).
\end{eqnarray}

\noindent
In what follows, we shall fix the choice of $\varphi_{\mu}$ by requiring that $\left(\varphi_{\mu}\right)_{x_1=0}=0$.

\medskip
Associated to $\kappa_{\mu,\eps,h}$, we can construct an elliptic $\widetilde{h}$--Fourier integral
operator of the form
\begeq
\label{eq15}
Gu(x_1)=\frac{1}{2\pi \widetilde{h}}\int\!\!\!\!\int
e^{\frac{i}{\widetilde{h}}
(\varphi_{\mu}(x_1,\widetilde{\eta_1},\frac{\eps}{\mu^2},\frac{h}{\mu^2},\eta_2)+(x_1-y_1)\widetilde{\eta_1})}
a(x_1,\widetilde{\eta_1},\frac{\eps}{\mu^2},\frac{h}{\mu^2},\eta_2;\widetilde{h})u(y_1)\,dy_1
d\widetilde{\eta_1},
\endeq
such that the full symbol of the $\widetilde{h}$-pseudodifferential operator
\begeq
\label{eq16}
G^{-1} \mu^{-2} e^{-\frac{i}{h}f(\xi_2)x_1}P'_{\eps}(x_1,\mu \widetilde{h}D_{x_1},\xi_2;h)
e^{\frac{i}{h}f(\xi_2)x_1}G
\endeq
is independent of $x_1$ (and of $x_2$), with the principal symbol given by
(\ref{eq14}). For the amplitude in (\ref{eq15}), we shall require that
$\left(a\right)_{x_1=0}=1$.

\medskip
\noindent
{\it Remark.} Working microlocally in a region
$$
\abs{\xi_1}\sim \mu,
$$
where $\mu \ll 1$ is such that
\begeq
\label{eq16.5}
\frac{(\eps+h)^{1/2}}{\mu}\leq h^{\delta_1},\,\,\,\delta_1>0,
\endeq
and following some further arguments of section 4 of~\cite{HiSj1}, we see that
the canonical transformation $\kappa_{\mu,\eps,h}$ and the
$\widetilde{h}$--Fourier integral operator $G$ in (\ref{eq15}) can be
constructed by a formal Taylor series in the
asymptotically small parameter $(\eps+h)/\mu^2={\cal O}(h^{2\delta_1})$.

\medskip
\noindent
{\it Remark.} Assume that the subprincipal symbol of $P_{\eps=0}$ in
(\ref{eq00}) vanishes. Then it follows from some arguments
in sections 2 and 4 in~\cite{HiSj1} that the $x_1$--dependence in
$P'_{\eps}(x_1,hD_{x_1},\xi_2;h)$ in (\ref{eq8}) can be eliminated
microlocally in a region $\abs{\xi_1}\sim \mu$, where
$$
(\eps+h^2)^{1/2}\ll \mu \ll 1.
$$

\bigskip
By rescaling, we can express $G$ in (\ref{eq15}) as an $h$--Fourier integral
operator. Indeed, using that $\frac{d\widetilde{\eta_1}}{\widetilde{h}}=\frac{d\eta_1}{h}$, we get
\begin{eqnarray}
\label{eq17}
Gu(x_1)=\frac{1}{2\pi h}\int\!\!\!\!\int e^{\frac{i}{h}\left(\mu
\varphi_{\mu}(x_1,\frac{\eta_1}{\mu},\frac{\eps}{\mu^2},\frac{h}{\mu^2},\xi_2)+(x_1-y_1)\eta_1\right)}\times
\\ \nonumber
a(x_1,\frac{\eta_1}{\mu}, \frac{\eps}{\mu^2},
\frac{h}{\mu^2},\xi_2;\frac{h}{\mu})u(y_1)\,dy_1\,d\eta_1.
\end{eqnarray}
Moreover, the introduction of the small parameter $\mu$ in (\ref{eq10}) was artificial, and
the\-re\-fore we can carry out the const\-ruc\-tions in such a way that the phase function
$\mu \varphi_{\mu}(x_1,\frac{\eta_1}{\mu},\frac{\eps}{\mu^2},\frac{h}{\mu^2},\xi_2)$ and the
amplitude $a(x_1,\frac{\eta_1}{\mu},\frac{\eps}{\mu^2},\frac{h}{\mu^2},\xi_2;\frac{h}{\mu})$ in (\ref{eq17})
are independent of $\mu$. We write then
\begin{eqnarray}
\label{eq18}
Gu(x_1)=\frac{1}{2\pi h}\int\!\!\!\!\int e^{\frac{i}{h}\left(
\varphi_{{\rm new}}(x_1,\eta_1,\eps,h,\xi_2)+(x_1-y_1)\eta_1\right)}\times
\\ \nonumber
a_{{\rm new}}(x_1,\eta_1,\eps,h,\xi_2;h)u(y_1)\,dy_1\,d\eta_1,
\end{eqnarray}
with $\varphi_{{\rm new}}$, $a_{\rm new}$ defined for $\eps+h \ll \abs{\eta_1}^2 \ll 1$ and satisfying
\begeq
\label{eq19}
\varphi_{{\rm new}}={\cal O}\left(\frac{\eps+h}{\abs{\eta_1}}\right),
\endeq
and
\begeq
\label{eq20}
a_{{\rm new}}\sim \sum_{j=0}^{\infty} a_{{\rm new},j}h^j,\quad a_{{\rm
new},j}={\cal O}(\abs{\eta_1}^{-2j}).
\endeq
Here $a_{{\rm new},j}$ do not depend on $h$. Since we work in the
complex domain, we can estimate the derivatives of $\varphi_{{\rm
new}}$ and $a_{{\rm new},j}$ using the Cauchy inequalities. In
particular, when $(\alpha_1,\beta_1)\in \nat^2$, we get using
(\ref{eq19}),
\begeq
\label{eq20.5}
\partial_{x_1}^{\alpha_1}\partial_{\eta_1}^{\beta_1} \varphi_{{\rm
new}}={\cal
O}_{\alpha_1 \beta_1}\left(\frac{(\eps+h)}{\abs{\eta_1}^{1+\abs{\beta_1}}}\right).
\endeq

\medskip
Since, as we have just observed,
$$
\mu \varphi_{\mu}\left(x_1,\frac{\eta_1}{\mu},\frac{\eps}{\mu^2}, \frac{h}{\mu^2},\xi_2\right)=
\varphi_{{\rm new}}\left(x_1,\eta_1,\eps,h,\xi_2\right),
$$
where $\varphi_{{\rm new}}$ satisfies (\ref{eq19}), it follows that the phase $\varphi_{\mu}$ in (\ref{eq13})
extends to a region $1 \ll \abs{\widetilde{\eta_1}} \ll \frac{1}{\mu}$ and satisfies there
$$
\varphi_{\mu}\left(x_1,\widetilde{\eta_1},\frac{\eps}{\mu^2},\frac{h}{\mu^2},\xi_2\right)
={\cal O}\left(\frac{(\eps+h)}{\mu^2
\abs{\widetilde{\eta_1}}}\right).
$$
Similarly, the normal form (\ref{eq16}) corresponds, after a multiplication by $\mu^2$,
to an operator which is independent of $\mu$,
\begin{eqnarray}
\label{eq21}
& &P_{\eps}''(hD_{x_1},\xi_2;h)=G^{-1}
e^{-\frac{i}{h}f(\xi_2)x_1}P'_{\eps}(x_1,hD_{x_1},\xi_2;h)e^{\frac{i}{h}f(\xi_2)x_1}G
\\ \nonumber & = & p(f(\xi_2),\xi_2)+g(f(\xi_2)+hD_{x_1})(hD_{x_1})^2
\\ \nonumber
& & + \left(i\eps \langle{\widetilde{\langle{q}\rangle_2}\rangle}_1+{\cal
O}(h)+\eps^2 \langle{r_{\eps}}\rangle_1\right)(f(\xi_2)+hD_{x_1},\xi_2) \\
\nonumber
& & + {\rm Op}_h \left({\cal O}\left(\frac{(\eps+h)^2}{\xi_1^2}\right)\right)+R(hD_{x_1},\xi_2,\eps;h),
\end{eqnarray}
where
$$
R\sim \sum_{j=2}^{\infty} h^j R_j(\xi),\quad R_j(\xi)={\cal
O}\left(\frac{1}{\abs{\xi_1}^{2j-2}}\right).
$$

For future reference we remark that we can also view the operator $G$
in (\ref{eq15}) as acting on (Floquet periodic) functions
on ${\bf T}^2$. If we maintain the scaling, we get
\begeq
\label{eq23}
Gu(x)=\frac{1}{(2\pi \widetilde{h})(2\pi h)}\int\!
e^{\frac{i}{\widetilde{h}}(\varphi_{\mu}(x,\widetilde{\eta_1})
+(x_1-y_1)\widetilde{\eta_1})+\frac{i}{h}(x_2-y_2)\eta_2}
a(x_1,\widetilde{\eta_1},\eta_2;\widetilde{h},h)u(y)\,dyd\widetilde{\eta_1}\,d\eta_2,
\endeq
where $\eta_2$ is the same variable as $\xi_2$. Without the scaling,
we have a similar formula by adding a $y_2$, $\eta_2$--integration to
(\ref{eq17}) (after replacing $\xi_2$ there by $\eta_2$), and adding a
phase factor $e^{\frac{i}{h}(x_2-y_2)\eta_2}$.

\medskip
Naturally, the argument so far is formal, with the various normal forms
computed by formal stationary phase expansions. Also, let us recall that the phase
$\varphi_{{\rm new}}$ in (\ref{eq18}) is defined only for
$(\eps+h)^{1/2} \ll \abs{\eta_1} \ll 1$.

\medskip
We summarize the discussion in this section in the following proposition.

\begin{prop}
Let $P_{\eps}$ be an $h$-pseudodifferential operator defined microlocally near $\xi=0$ in $T^*{\bf T}^2$, and assume
that the principal symbol of $P_{\eps}$,
$$
p_{\eps}(x,\xi)=p(\xi)+i\eps q(x,\xi)+{\cal O}(\eps^2),
$$
is such that $p(\xi)$ satisfies {\rm (\ref{eq0.6})}, {\rm (\ref{eq1})}. Then we write
\begeq
\label{eq23.1}
p(\xi_1,\xi_2)=p(f(\xi_2),\xi_2)+g(\xi_1,\xi_2)(\xi_1-f(\xi_2))^2,\quad f(0)=0,
\endeq
where $g(0,0)>0$. For each $N\in \nat$ there exists an elliptic Fourier integral operator
in the complex domain $F=F_{\eps}^{(N)}$ such that the symbol of $M^{-1}F^{-1}P_{\eps} FM$ is of the form
\begin{eqnarray}
\label{eq23.2}
P'_{\eps}(x_1,\xi_1+f(\xi_2),\xi_2;h)+{\cal O}\left(\eps^{N+1}+\eps \xi_1^N+h^{N+1}\right).
\end{eqnarray}
Here $M$ is the operator of multiplication by $e^{\frac{i}{h}f(\xi_2)x_1}$, and $P'_{\eps}(x_1,hD_{x_1},hD_{x_2};h)$
is defined in {\rm (\ref{eq8})}.

Furthermore, let $(\eps+h)^{1/2} \ll \mu \ll 1$, and let us view $\mu^{-2} P'_{\eps}(x_1,hD_{x_1}+f(\xi_2),\xi_2;h)$ as an
$\widetilde{h}$--pseudodifferential operator in $x_1$, with $\widetilde{h}=h/\mu$. There exists an elliptic
$\widetilde{h}$--Fourier integral operator $G$ in $x_1$,
defined in {\rm (\ref{eq15})}, microlocally in
$\abs{\widetilde{\xi}}\sim 1$, such that the full symbol of $G^{-1}\mu^{-2} P'_{\eps}(x_1,hD_{x_1}+f(\xi_2),\xi_2;h)G$
is independent of $x_1$. The operator $G$ quantizes a holomorphic canonical transformation whose ge\-ne\-rating function is of the form
$x_1\widetilde{\eta_1}+\varphi_{\mu}(x_1,\widetilde{\eta_1})$, where $\varphi_{\mu}$ is defined in
$1\ll \abs{\widetilde{\eta_1}}\ll \frac{1}{\mu}$ and satisfies there
\begeq
\label{eq23.3}
\varphi_{\mu}(x_1,\widetilde{\eta_1})={\cal O}\left(\frac{\eps+h}{\mu^2 \abs{\widetilde{\eta_1}}}\right).
\endeq
In this region we have, when $(\alpha_1,\beta_1)\in \nat^2$,
\begeq
\label{eq23.4}
\partial^{\alpha_1}_{x_1}\partial^{\beta_1}_{\widetilde{\eta_1}} \varphi_{\mu}=
{\cal O}_{\alpha_1 \beta_1}\left(\frac{\eps+h}{\mu^2 \abs{\widetilde{\eta_1}}^{1+\abs{\beta_1}}}\right).
\endeq
\end{prop}

\section{Microlocal Hilbert spaces near the rational torus}
\setcounter{equation}{0}
Let $P_{\eps}$ be as in section 1. In section 3, we have constructed a microlocal normal form for $P_{\eps}$
near the rational Lagrangian torus $\Lambda_{1,r}\subset p^{-1}(0)\cap
\real^4$, but away from an ${\cal O}((\eps+h)^{1/2})$--\neigh{} of this set --- see
(\ref{eq21}). The purpose of this section is to follow up the preceding formal
constructions with suitable function spaces and to construct a microlocal Hilbert space in a sufficiently small but fixed
\neigh{} of $\Lambda_{1,r}$, implementing the reduction scheme of Proposition 3.1.

\subsection{Microlocal Hilbert spaces outside of a tiny \neigh{} of $\Lambda_{1,r}$}
Let us consider an operator $P_{\eps}$, microlocally defined near $\xi=0$ in $T^*{\bf T}^2$, with the leading symbol given by
(\ref{eq2}), (\ref{eq2.1}). We shall work as much as possible with
functions on ${\bf T}^2$, and
with corresponding Fourier integral operators operating in 2
variables. Adopting this point of view, we see that the multiplication by
$e^{\frac{i}{h}f(\xi_2)x_1}$, introduced in (\ref{eq9}), can be viewed as the semiclassical
Fourier integral operator
\begin{eqnarray}
\label{eq24}
Mu(x)& = & \frac{1}{2\pi h}\int\!\!\!\!\int e^{\frac{i}{h}\left(f(\eta_2)x_1+(x_2-y_2)\eta_2\right)}u(x_1,y_2)
\,dy_2\,d\eta_2 \\ \nonumber
& = & \frac{1}{(2\pi h)^2} \int\!\!\!\!\int
e^{\frac{i}{h}\left(f(\eta_2)x_1+(x-y)\cdot
\eta\right)}u(y)\,dy\,d\eta,
\end{eqnarray}
associated to the \cantransf{}
\begeq
\label{eq25}
\kappa_M: (x_1,x_2+f'(\eta_2)x_1; \eta_1,\eta_2)\mapsto (x_1,x_2;
\eta_1+f(\eta_2),\eta_2).
\endeq
Let us recall now the operators $F$ and $G$, introduced in (\ref{eq7}) and (\ref{eq23}), respectively.
In the previous section we have obtained that formally,
\begeq
\label{eq26}
G^{-1}M^{-1}F^{-1}P_{\eps}FMG=P_{\eps}''(hD_x,h)+(M G)^{-1}R_{\eps} M G,
\endeq
with $P''_{\eps}$ and $R_{\eps}$ given in (\ref{eq21}) and (\ref{eq7.1}), respectively. The fact that the phase
$\varphi_{\mu}$ in (\ref{eq23}) (see also Proposition 3.1) is only defined for
$1 \ll \abs{\widetilde{\eta_1}} \ll \frac{1}{\mu}$, $(\eps+h)^{1/2}\ll \mu \ll 1$, is a difficulty that
we shall address later in this section. Ignoring that problem for a moment and still arguing formally, we would like
to consider $P_{\eps}''$ acting on the space $L^2_{\theta}({\bf
T}^2)$, microlocally defined near the zero section, but away from
the exceptional region
$$
\abs{\xi_1} = {\cal O}((\eps+h)^{1/2}).
$$
Consequently, the natural formal Hilbert space for considering $P_{\eps}$ should be given by
$FMG(L^2_{\theta}({\bf T}^2))$. When realizing the latter, it is going to be convenient to work on the FBI transform side.

\medskip
We shall work with the standard FBI--Bargmann transform,
\begeq
\label{eq27}
T u(x)= T_{h,h} u(x)=Ch^{-3/2} \int e^{-\frac{1}{2h}(x-y)^2}u(y)\,dy, \quad C>0,
\endeq
acting on $L^2_{\theta}({\bf T}^2)$, and mapping this space to a
weighted space of Floquet periodic holomorphic functions on $\comp^2$. Associated to $T$, there is a canonical transformation
\begeq
\label{eq27.01}
\kappa_{T_{h,h}}=\kappa_T: (y,\eta)\mapsto (x,\xi)=(y-i\eta,\eta),
\endeq
mapping the real phase space $T^*{\bf T}^2$ to the IR-manifold
\begeq
\label{eq27.01.1}
\Lambda_{\Phi_0}:\,\, \xi=\frac{2}{i}\frac{\partial \Phi_0}{\partial x}=-\Im x,\,\, \Phi_0(x)=\frac{1}{2}(\Im x)^2.
\endeq

\noindent
Let us also recall that the transformation
\begeq
\label{eq27.02}
T: L^2({\bf T}^2)\rightarrow H_{\Phi_0}(\comp^2/2\pi \z^2)
\endeq
is unitary, for a suitable choice of $C>0$ in (\ref{eq27}),
and it has been verified in section 3 of~\cite{MeSj2} that it remains unitary when acting on the Floquet space
$L^2_{\theta}({\bf T}^2)$. Here and in what follows, when $\Omega\subset \comp^2/2\pi \z^2={\bf T}^2+i\real^2$ is open
and $\Phi$ is a suitable strictly plurisubharmonic weight, close to $\Phi_0$ in (\ref{eq27.01.1}),
we shall let $H_{\Phi}(\Omega)$ stand for the closed subspace of $L^2(\Omega; e^{-\frac{2\Phi}{h}} L(dx))$,
consisting of functions that are holomorphic in $\Omega$ --- see also
the appendix.

Neglecting the Floquet conditions for the time being, we should have,
\begeq
\label{eq27.1}
T F M G(L^2({\bf T}^2))=H_{\Phi},
\endeq
where the weight $\Phi$ is such that
$$
\Lambda_{\Phi}:=\left\{\left(x,\frac{2}{i}\frac{\partial \Phi}{\partial x}\right)\right\}=\kappa_T\circ
\kappa_{\eps}\circ \kappa_M\circ \kappa_{\mu,\eps,h}(T^*{\bf T}^2).
$$
Here $\kappa_{\eps}$ and $\kappa_{\mu,\eps,h}$ are the canonical transformations corresponding to $F$ and $G$, and
introduced in (\ref{eq5}) and (\ref{eq13.5}), respectively. The weight $\Phi$ in (\ref{eq27.1}) should be a small
perturbation of $\Phi_0$ since $\kappa_{\eps}$, $\kappa_{\mu,\eps,h}$ are small perturbations of the identity,
and $\kappa_M$ in (\ref{eq25}) is a real \cantransf{}.

\medskip
We shall assume from now on that
\begeq
\label{eq27.2}
\eps \gg h,
\endeq
and abusing the previous notation slightly, we shall take
\begeq
\label{eq27.3}
\mu=\sqrt{\eps}.
\endeq
Because of the blow-up of the normal form construction in the region
where $\eta_1={\cal O}(\eps^{1/2})$ (see (\ref{eq19})), when realizing the formal space in (\ref{eq27.1}), we shall
have to make some modifications. First, the operator $G$ should be written as in (\ref{eq23}) with
$$
\widetilde{h}=\frac{h}{\sqrt{\eps}}\ll 1,
$$
and correspondingly, in order to define
a microlocal space corresponding to the formal space $G(L^2({\bf T}^2))$, we shall consider the mixed transform
\begeq
\label{eq28}
T_{\widetilde{h},h}u(x)=C \widetilde{h}^{-3/4}h^{-3/4} \int\!\!\!\!\int
e^{-\frac{1}{2\widetilde{h}}(x_1-y_1)^2-\frac{1}{2h}(x_2-y_2)^2} u(y_1,y_2)\,dy_1\,dy_2.
\endeq
Here $C>0$ is the same constant as in (\ref{eq27}).
For future reference, we notice that when viewed as an $h$--Fourier integral operator, the transform
$T_{\widetilde{h},h}$ is associated with the canonical transformation
\begeq
\label{eq28.5}
\kappa_{T_{\widetilde{h},h}}(y_1,\eta_1;y_2,\eta_2)=\left(y_1-i\frac{\eta_1}{\sqrt{\eps}},\eta_1;y_2-i\eta_2,\eta_2\right).
\endeq
Here we have written $(y_1,\eta_1;y_2,\eta_2)$ rather than $(y,\eta)$.

\medskip
We shall show that $T_{\widetilde{h},h}G(L^2({\bf T}^2))$ becomes a
well-defined exponentially weighted space of holomorphic functions
$u(x_1,x_2)$ in a region $1\ll \abs{\Im x_1} \ll \frac{1}{\mu}$,
$\abs{\Im x_2} \ll 1$. Once this has been done and the basic properties of the
weight have been investigated, we shall extend the definition of the weight
to the entire domain $\abs{\Im x_1} \ll
\frac{1}{\mu}$, $\abs{\Im x_2} \ll 1$ --- this will then lead to a definition of a microlocal Hilbert space corresponding
to a formal space $G(L^2({\bf T}^2))$, in a full \neigh{} of the
rational torus, and we shall be able to proceed as indicated above.

Let us compute $T_{\widetilde{h},h}Gu$, when $u\in L^2$. In doing
so, it will be convenient to do the computation first
in the $x_1$-variable alone, and as in (\ref{eq15}), we introduce, with $\mu=\sqrt{\eps}$,
\begin{eqnarray}
\label{eq29}
& & Gu(x_1,y_2,\eta_2) \\ \nonumber
& = & \frac{1}{2\pi \widetilde{h}} \int\!\!\!\!\int
e^{\frac{i}{\widetilde{h}}\left(\varphi_{\mu}\left(x_1,\widetilde{\eta_1},
\frac{h}{\eps},\eta_2\right)+(x_1-y_1)\widetilde{\eta_1}\right)}
a(x_1,\widetilde{\eta_1},\frac{h}{\eps}, \eta_2;\widetilde{h})
u(y_1,y_2)\,dy_1\,d\widetilde{\eta_1}.
\end{eqnarray}
Composing this expression with the one-variable transform
$T_{\widetilde{h}}$, we get
\begin{eqnarray}
\label{eq30}
& & T_{\widetilde{h}}G u(z_1,y_2,\eta_2) \\ \nonumber
& = & \frac{C \widetilde{h}^{-\frac{3}{4}}} {2\pi \widetilde{h}}
\int\!\!\!\!\int\!\!\!\!\int e^{\frac{i}{\widetilde{h}}\left(\frac{i}{2}(z_1-x_1)^2+\varphi_{\mu}(x_1,\widetilde{\eta_1}, \frac{h}{\eps},\eta_2)
+ (x_1-y_1)\widetilde{\eta_1}\right)} a(x_1, \widetilde{\eta_1}, \frac{h}{\eps}, \eta_2;\widetilde{h})u(y_1,y_2)\,
dy_1d\widetilde{\eta_1}dx_1 \\ \nonumber
& = & C_1 \widetilde{h}^{-\frac{3}{4}} \int e^{\frac{i}{\widetilde{h}}\left(\frac{i}{2}(z_1-y_1)^2+
\widetilde{\varphi}_{\mu}(z_1,y_1,\frac{h}{\eps}, \eta_2)\right)} b(z_1,y_1,\frac{h}{\eps},\eta_2;\widetilde{h})
u(y_1,y_2)\,dy_1, \quad C_1>0.
\end{eqnarray}
Here the last expression follows from the stationary phase method in
the variables $x_1$, $\widetilde{\eta_1}$~\cite{Sj82}, whereby we notice that the
critical point of the phase in (\ref{eq30}), $(x_1^c,\widetilde{\eta_1}^c)$, satisfies
$$
x_1^c=y_1+{\cal O}\left(\frac{1}{\abs{z_1-y_1}^2}\right),\quad \widetilde{\eta_1}^c=i(z_1-y_1)+
{\cal O}\left(\frac{1}{\abs{z_1-y_1}}\right).
$$
It follows that the phase $\widetilde{\varphi}_{\mu}$ in (\ref{eq30}) is a well-defined holomorphic
function of $z_1$ in a region $1\ll \abs{z_1-y_1}\ll 1/\eps^{1/2}$, $\abs{\Im z_1}\gg \abs{\Re z_1-y_1}$,
and enjoys the same estimates as $\varphi_{\mu}$ in (\ref{eq23.3}), (\ref{eq23.4}),
\begeq
\label{eq31}
\widetilde{\varphi}_{\mu}={\cal O}\left(\frac{1}{\abs{\Im z_1}}\right),
\quad
\partial^l_{y_1} \partial^m_{z_1} \widetilde{\varphi}_{\mu}=
{\cal O}_{lm}\left(\frac{1}{\abs{ \Im z_1}^{1+m}}\right).
\endeq
Here, as before, the estimates on the derivatives of $\widetilde{\varphi}_{\mu}$ follow from the
Cauchy inequali\-ties.

It follows from (\ref{eq30}) that
\begeq
\label{eq32}
T_{\widetilde{h}}G u(z_1,y_2,\eta_2)\in H_{\Phi_1(\cdot,\eta_2), \widetilde{h}},
\endeq
in the region $1 \ll \abs{\Im z_1} \ll \frac{1}{\eps^{1/2}}$, where
\begin{eqnarray}
\label{eq32.1}
\Phi_1(z_1,\eta_2) & = & \sup_{y_1\in {\bf R}} \left(-\frac{1}{2} \Re(z_1-y_1)^2-\Im \widetilde{\varphi}_{\mu}(z_1,
y_1, \frac{h}{\eps}, \eta_2)\right)\\ \nonumber
& = & \frac{1}{2}(\Im z_1)^2+\Phi_2(z_1,\eta_2),
\end{eqnarray}
and
\begeq
\label{eq33}
\Phi_2(z_1,\eta_2)={\cal O}\left(\frac{1}{\abs{\Im z_1}}\right).
\endeq
The critical point $y_1^c$ corresponding to the supremum in
(\ref{eq32.1}) satisfies
\begeq
\label{eq33.1}
y_1^c=\Re z_1+{\cal O}\left(\frac{1}{\abs{\Im z_1}}\right).
\endeq
Using (\ref{eq31}) together with a scaling argument very similar to
the one described in detail in the proof of Proposition 4.3 below, we see that the ${\cal
O}$--term in (\ref{eq33.1}) satisfies
$$
\partial^k_{{\rm Re} z_1}\partial^l_{{\rm Im} z_1}{\cal
O}\left(\frac{1}{\abs{\Im z_1}}\right)={\cal O}\left(\frac{1}{\abs{\Im z_1}^{1+l}}\right).
$$
It follows that for $\Phi_2(z_1,\eta_2)$ in (\ref{eq32.1}) we have
$$
\partial^k_{{\rm Re} z_1}\partial^l_{{\rm Im}
z_1}\Phi_2(z_1,\eta_2)={\cal O}_{k l}\left(\frac{1}{\abs{\Im
z_1}^{1+l}}\right).
$$

If we now let $G$ denote the full 2-variable operator in (\ref{eq23}), we get from (\ref{eq30}),
\begin{eqnarray}
\label{eq34}
& & T_{\widetilde{h},h} G u = C_2 \widetilde{h}^{-\frac{3}{4}} h^{-\frac{3}{4}-1}
\int\!\!\!\!\int\!\!\!\!\int\!\!\!\!\int\!\!\!
e^{\frac{i}{\widetilde{h}}
\left(\frac{i}{2}(z_1-y_1)^2+\widetilde{\varphi}_{\mu}(z_1,y_1,\frac{h}{\eps},\eta_2)\right)
+\frac{i}{h} \left(\frac{i}{2}(z_2-x_2)^2+(x_2-y_2)\eta_2\right)} \\ \nonumber
& \times & b(z_1,y_1,\frac{h}{\eps}, \eta_2;\widetilde{h}) u(y_1,y_2) \,dy_1\,dy_2\,d\eta_2\,dx_2 = C_3
\widetilde{h}^{-\frac{3}{4}} h^{-\frac{3}{4}}\times \\ \nonumber
& \times & \int\!\!\!\!\int
e^{\frac{i}{h}\left(\frac{i}{2}\sqrt{\eps}(z_1-y_1)^2+\frac{i}{2}(z_2-y_2)^2+
\sqrt{\eps}\widehat{\varphi}_{\mu}(z_1,y_1,z_2-y_2,\frac{h}{\eps})\right)} c(z_1,y_1, z_2-y_2;\widetilde{h},h)
u(y_1,y_2)\,dy_1\,dy_2,
\end{eqnarray}
where the last identity follows from stationary phase in $x_2$, $\eta_2$, and $\widehat{\varphi}_{\mu}$ satisfies the
same estimates as $\widetilde{\varphi}_{\mu}$,
\begeq
\label{eq35}
\partial^k_{z_2} \partial^l_{z_1} \widehat{\varphi}_{\mu}
={\cal O}_{k l}\left(\frac{1}{\abs{\Im z_1}^{1+l}}\right).
\endeq

It follows that
\begeq
\label{eq36}
T_{\widetilde{h},h} Gu\in H_{\Phi_{3,h}}, \quad \Phi_3(z_1,\Im z_2)=
\frac{\sqrt{\eps}}{2} (\Im z_1)^2+ \frac{1}{2}(\Im z_2)^2 +\Phi_4(z_1,\Im z_2),
\endeq
and
\begeq
\label{eq37}
\partial_{{\rm Re} z_1, {\rm Im} z_2}^k \partial_{{\rm Im} z_1}^l
\Phi_4(z_1,\Im z_2)={\cal O}_{k l}\left(\frac{\sqrt{\eps}}{\abs{\Im z_1}^{1+l}}\right).
\endeq
Here $\Phi_3$, $\Phi_4$ are independent of $\Re z_2$.

The discussion above is summarized in the following, somewhat informal, proposition.

\begin{prop}
Let us assume that $\eps\gg h$ and set $\widetilde{h}=h/\sqrt{\eps}$. Via the $(\widetilde{h},h)$--Bargmann transform $T_{\widetilde{h},h}$ defined in
{\rm (\ref{eq28})}, the formal space $G(L^2({\bf T}^2))$ corresponds to the weighted space of holomorphic functions
$H_{\Phi_3,h}$ in the region
$1\ll \abs{\Im z_1}\ll \frac{1}{\sqrt{\eps}}$, $\abs{\Im z_2}\ll 1$. The weight $\Phi_3=\Phi_3(z_1,\Im z_2)$ is such that
\begeq
\label{eq37.5}
\Phi_3(z_1,\Im z_2)=\frac{\sqrt{\eps}}{2}(\Im z_1)^2+\frac{1}{2}(\Im z_2)^2+\Phi_4(z_1,\Im z_2),
\endeq
where the perturbation $\Phi_4$ satisfies
\begeq
\label{eq37.6}
\partial_{{\rm Re} z_1, {\rm Im} z_2}^k \partial_{{\rm Im} z_1}^l
\Phi_4(z_1,\Im z_2)={\cal O}_{k l}\left(\frac{\sqrt{\eps}}{\abs{\Im z_1}^{1+l}}\right).
\endeq
The corresponding statement also holds when considering the formal space $G(L^2_{\theta}({\bf T}^2))$ of Floquet periodic
functions.
\end{prop}

\noindent
{\it Remark}. Let us remark that the cutoff and remainder errors not
written out explicitly in
the stationary phase expansions above
are all of the size ${\cal O}(1) \exp(-1/C\widetilde{h})={\cal O}(1)
\exp(-\frac{\sqrt{\eps}}{Ch})$~\cite{Sj82}, while the
deviation of the weight, due to $\Phi_4$, corresponds to an exponential factor
$$
\exp({\cal O}(1)\frac{\sqrt{\eps}}{h\abs{{\rm Im} z_1}})\ll
\exp(\frac{\sqrt{\eps}}{Ch}),
$$
since we work in a region where $\abs{\Im z_1} \gg 1$.

\noindent
{\it Remark}. Constructing and working with the $\widetilde{h}$--Fourier integral
operator $G$ in the domain where
$$
\frac{\sqrt{\eps}}{h^{\delta_1}}\leq \abs{\xi_1}\ll 1,
$$
for some $\delta_1>0$ small (see also (\ref{eq16.5})), we find that the formal space $G(L^2({\bf T}^2))$
corresponds, via the $(\widetilde{h},h)$--Bargmann transform, to the space $H_{\Phi_3,h}$,
as in Proposition 4.1, now viewed in the region $h^{-\delta_1}\leq
\abs{\Im z_1}\ll \frac{1}{\sqrt{\eps}}$, $\abs{\Im z_2}\ll 1$.

\subsection{Fourier series expansions in $H_{\Phi}$--spaces}

The purpose of this subsection is to obtain a relation between the 1-variable
weight $\Phi_1(z_1,\eta_2)$ and the 2-variable weight $\Phi_3(z_1,\Im
z_2)$, introduced in (\ref{eq32.1}) and (\ref{eq36}),
respectively. The starting point for us will be the following remark concerning Fourier series on the FBI transform side.
Let us rewrite (\ref{eq29}) with
slightly different notation, now for a function $u\in L^2$ of one variable only:
\begeq
\label{eq38}
G_{\eta_2} u(x_1) = \frac{1}{2\pi \widetilde{h}} \int\!\!\!\int
e^{\frac{i}{\widetilde{h}}
\left(\varphi_{\mu}(x_1,\widetilde{\eta_1},\frac{h}{\eps},\eta_2)+(x_1-y_1)\widetilde{\eta_1}\right)}
a(x_1,\widetilde{\eta_1}, \frac{h}{\eps}, \eta_2;\widetilde{h}) u(y_1)\,dy_1\,d\widetilde{\eta_1}.
\endeq
If $u=u(y_1,y_2)\in L^2({\bf T}^2)$ depends on 2 variables and we introduce the Fourier series expansion in $y_2$,
\begeq
\label{eq39}
u(y_1,y_2)=\sum_{k\in {\bf Z}} e^{\frac{i}{h}y_2 kh} \widehat{u}(y_1,kh),
\endeq
then
$$
G u(x_1,x_2) = \sum_{k\in {\bf Z}} e^{\frac{i}{h}x_2 kh} \left(G_{kh}\widehat{u}(\cdot,kh)\right)(x_1),
$$
and therefore, applying $T_{\widetilde{h},h}$ of (\ref{eq28}), we get
\begeq
\label{eq40}
T_{\widetilde{h},h}Gu(z_1,z_2)=\sum_{k\in {\bf Z}}T^{(2)}_h(e^{\frac{i}{h}x_2 k h})(z_2)
T_{\widetilde{h}}^{(1)}G_{kh} \widehat{u}(\cdot,k h)(z_1).
\endeq
Here the superscripts $(1)$, $(2)$ in (\ref{eq40}) indicate the
variable in which the corresponding operators are applied. A straightforward computation shows that
\begeq
\label{eq41}
T^{(2)}_h\left(e^{\frac{i}{h}(\cdot)\xi_2}\right)(z_2)=C h^{-\frac{1}{4}}
e^{-\frac{\xi_2^2}{2h}}e^{\frac{i}{h}z_2\xi_2}=: e_{\xi_2}(z_2),\quad
C>0,
\endeq
and clearly, as can also be verified directly, this function is normalized in $H_{\Phi_0}(\comp/2\pi \z)$,
$\Phi_0(z_2)=\frac{1}{2}(\Im z_2)^2$. The functions $e_{kh}$,
$k\in \z$, form an orthonormal basis in this space, and hence a general element of
$H_{\Phi_0}(\comp/2\pi \z)$ has an expansion
\begeq
\label{eq42}
v=\sum_{k\in {\bf Z}} \widetilde{v}_k e_{kh},\quad \widetilde{v}_k=(v|e_{kh})_{H_{\Phi_0}}.
\endeq

We shall now pause to review Fourier series expansions in $H_{\Phi}(\comp/2\pi
\z)$, where $\Phi=\Phi(\Im z)$ is a general smooth weight such that $t\mapsto \Phi(t)$ is strictly convex:
\begeq
\label{eq43}
v(z)=\sum_{k\in {\bf Z}} \widehat{v}_k e^{\frac{i}{h}z k h}.
\endeq
Here the scalar product
$$
\int_{{\bf C}/2\pi {\bf Z}} e^{\frac{i}{h}z kh}\overline{e^{\frac{i}{h} z \ell h}}
e^{-\frac{2\Phi({\rm Im} z)}{h}}\,L(dz),\quad
k,l\in \z
$$
vanishes for $\ell\neq k$ and for $k=\ell$ it is equal to
\begeq
\label{eq43.5}
2\pi \int e^{-\frac{2}{h}\left(\Phi({\rm Im} z)+kh {\rm Im} z\right)}\,d\Im z,
\endeq
which can be evaluated by the method of stationary phase. The critical
point $t=\Im z$ in (\ref{eq43.5}) is given by $\Phi'(t)+kh=0$ and
$\widetilde{\Phi}(kh):=\inf (kh t+\Phi(t))=-\sup ((-kh)t-\Phi(t))=-{\cal L}\Phi(-kh)$, where ${\cal L}$ is the Legendre
transformation. Notice also that the critical point $t$ can be characterized by
$$
\frac{2}{i} \frac{\partial \Phi}{\partial z}(x+it)=kh,\quad x\in \real,
$$
when identifying $\Phi(z)=\Phi(\Im z)$. Thus, by stationary phase (the Laplace method), we get
$$
\norm{e^{\frac{i}{h}(\cdot)kh}}_{H_{\Phi}}^2=h^{1/2} a_{\Phi}(kh;h) e^{\frac{2}{h}{\cal L}\Phi(-kh)},
$$
where $a_{\Phi}(t;h)\sim a_0(t)+h a_1(t)+\ldots $ is a positive elliptic symbol.

For the Fourier series expansion (\ref{eq43}) we therefore have the Parseval relation,
\begeq
\label{eq44}
\norm{v}_{H_{\Phi}}^2=\sum_{k\in {\bf Z}} h^{\frac{1}{2}} a_{\Phi}(kh;h) e^{\frac{2}{h}{\cal L} \Phi(-kh)}
\abs{\widehat{v}_k}^2,
\endeq
telling us to which weighted $l^2$-space the Fourier coefficients $\widehat{v}_k$ belong.

\bigskip
Applying (\ref{eq44}) to (\ref{eq40}),
(\ref{eq41}) viewed as a Fourier series in $z_2$ with $z_1$ as a parameter, we get
\begin{eqnarray}
\label{eq45}
& & \norm{T_{\widetilde{h},h} G u(z_1,\cdot)}^2_{H_{\Phi_3(z_1,\cdot),h}} \\ \nonumber
&=&\abs{C}^2 \sum_{k\in {\bf Z}}
a_{\Phi_3(z_1,\cdot)}(kh;h) e^{\frac{2}{h}{\cal L}\Phi_3(z_1,-kh)-\frac{(kh)^2}{h}}
\abs{T^{(1)}_{\widetilde{h}} G_{kh} \widehat{u}(\cdot,kh)(z_1)}^2.
\end{eqnarray}
Now recall that the weights $\Phi_3$ and $\Phi_1$ have been chosen so that
\begin{eqnarray}
\label{eq46}
& & \norm{T_{\widetilde{h},h}G u}^2_{H_{\Phi_{3,h}}}\sim \norm{u}^2_{L^2({\bf T}^2)}=
\sum_{k\in {\bf Z}} \norm{\widehat{u}_k}^2_{L^2({\bf T}^1)} \\ \nonumber
& = & \sum_{k \in {\bf Z}} \norm{T_{\widetilde{h}} G_{kh} \widehat{u}_k}^2_{H_{\Phi_1(\cdot,kh),\widetilde{h}}},
\end{eqnarray}
where as in (\ref{eq39}), $u(y_1,y_2)=\sum_{k\in {\bf Z}} e^{iky_2} \widehat{u}_{k}(y_1)$, and thus we want the
last member of (\ref{eq46}) to coincide with that of (\ref{eq45}) after an
integration with respect to $z_1$. This means that
$$
\frac{2}{h}({\cal L}\Phi_3)(z_1,-kh)-\frac{(kh)^2}{h}=-\frac{2\sqrt{\eps}}{h}\Phi_1(z_1,kh),
$$
so that
\begeq
\label{eq47}
({\cal L}\Phi_3)(z_1,-\eta_2)=\frac{1}{2}\eta_2^2-\sqrt{\eps} \Phi_1(z_1,\eta_2).
\endeq

When verifying (\ref{eq47}), we recall from (\ref{eq32.1}) that
\begeq
\label{eq48}
\Phi_1(z_1,\eta_2)=\sup_{y_1\in {\bf R}} \left(-\frac{1}{2}\Re (z_1-y_1)^2-
\Im \widetilde{\varphi}_{\mu}(z_1,y_1,\frac{h}{\eps},\eta_2)\right).
\endeq
We need a similar formula for $\Phi_3$. To that end, let us notice that in (\ref{eq34}) we can insert an
intermediate step, where we only integrate with respect to $x_2$, and exploiting that
$$
\int e^{\frac{i}{h}(x_2-y_2)\eta_2-\frac{1}{2h}(x_2-z_2)^2}\,dx_2=
\sqrt{2\pi h} e^{\frac{i}{h}(z_2-y_2)\eta_2-\frac{1}{2h}\eta_2^2},
$$
we get
\begin{eqnarray}
\label{eq49}
T_{\widetilde{h},h}G u(z) & = & C_3 \widetilde{h}^{-\frac{3}{4}}h^{-\frac{3}{4}-\frac{1}{2}}
\int\!\!\!\!\int\!\!\!\!\int\!\! e^{-\frac{1}{2\widetilde{h}}(z_1-y_1)^2+
\frac{i}{\widetilde{h}}\widetilde{\varphi}_{\mu}(z_1,y_1,\frac{h}{\eps},\eta_2)+
\frac{i}{h}(z_2-y_2)\eta_2-\frac{1}{2h}\eta_2^2} \\ \nonumber
& & \times\,\, b(z_1,y_1,\frac{h}{\eps},\eta_2;\widetilde{h}) u(y_1,y_2)\,dy_1\,dy_2\,d\eta_2.
\end{eqnarray}
The formula for $\Phi_3$ becomes
\begin{eqnarray}
\label{eq50}
& & \Phi_3(z_1,\Im z_2) \\ \nonumber
& = & \sup_{y_1\in {\bf R}} \sup_{y_2\in {\bf R}} {\rm vc}_{\eta_2} -\frac{\sqrt{\eps}}{2}\Re(z_1-y_1)^2
-\sqrt{\eps}\Im \widetilde{\varphi}_{\mu}-\Im (z_2-y_2)\eta_2-\frac{1}{2} \Re \eta_2^2.
\end{eqnarray}

For $y_1$ fixed, the $\sup_{y_2\in {\bf R}} {\rm vc}_{\eta_2}$ corresponds to taking the critical value with respect
to $y_2$, $\eta_2$ and the criticality with respect to $y_2$ requires $\eta_2$ to be real, making the right hand side
independent of $y_2$. Thus "${\rm vc}_{\eta_2}$" in (\ref{eq50}) can be replaced by "$\sup_{\eta_2\in {\bf R}}$"
and we get
\begin{eqnarray}
\label{eq51}
& & \Phi_3(z_1,\Im z_2) \\ \nonumber
& = & \sup_{y_1\in {\bf R}} \sup_{\eta_2\in {\bf R}} -\frac{\sqrt{\eps}}{2}\Re (z_1-y_1)^2-
\sqrt{\eps}\Im \widetilde{\varphi}_{\mu}(z_1,y_1,\frac{h}{\eps},\eta_2)-\frac{1}{2}\eta_2^2-\eta_2\Im z_2 \\ \nonumber
& = & \sup_{\eta_2\in {\bf R}} \sqrt{\eps} \Phi_1(z_1,\eta_2)-\frac{1}{2}\eta_2^2-\eta_2 \Im z_2 \\ \nonumber
& = & {\cal L}_{\eta_2\rightarrow {\rm Im} z_2} \left(\frac{1}{2}\eta_2^2-\sqrt{\eps} \Phi_1(z_1,\eta_2)\right)(-\Im z_2).
\end{eqnarray}
With $f(\eta_2)=\frac{1}{2}\eta_2^2-\sqrt{\eps} \Phi_1(z_1,\eta_2)$, $Ju(t)=u(-t)$, and $z_1$ treated as a parameter,
(\ref{eq47}) reads
$$
J{\cal L}\Phi_3=f,
$$
while (\ref{eq51}) tells us that $J{\cal L}f=\Phi_3$. Since $J^2={\cal L}^2=1$ and $J{\cal L}={\cal L}J$, we then see
that  (\ref{eq47}) follows from (\ref{eq51}).

\begin{prop}
Let the strictly plurisubharmonic weights $\Phi_1(z_1,\eta_2)$ and
$\Phi_3(z_1,\Im z_2)$ be defined in
{\rm (\ref{eq32.1})} and {\rm (\ref{eq50})}, respectively. Then we have the relation
\begeq
\label{eq51.1}
\left({\cal L}_{{\rm Im} z_2\rightarrow \eta_2} \Phi_3\right)(z_1,-\eta_2)=\frac{1}{2}\eta_2^2-\sqrt{\eps}\Phi_1(z_1,\eta_2).
\endeq
Here ${\cal L}f(\xi)=\sup_x (x\xi-f(x))$ is the Legendre transform of a strictly convex smooth function
$f:\real\rightarrow \real$.

It follows that if we have an expansion of $u\in H_{\Phi_3,h}$,
\begeq
\label{eq51.2}
u(z_1,z_2)=\sum_{k\in {\bf Z}} \widetilde{u}_k(z_1) e_{kh}(z_2),\quad
e_{kh}(z_2)=Ch^{-\frac{1}{4}} e^{-\frac{(kh)^2}{2h}} e^{\frac{i}{h}z_2 kh},\quad C>0,
\endeq
then
\begeq
\label{eq51.3}
\norm{u}_{H_{\Phi_3,h}}^2\sim \sum_{k\in {\bf Z}} \norm{\widetilde{u}_k}^2_{H_{\Phi_1(\cdot,kh),\widetilde{h}}}.
\endeq
\end{prop}

\subsection {Comparison with the ordinary transform away from $\Lambda_{1,r}$}

This subsection is a preparation for defining the global Hilbert space by
gluing together the local constructions near $\Lambda_{1,r}$ to the
weighted spaces that we used in~\cite{HiSjVu}. This discussion will
be continued in section 5.

In subsection 4.1 we have analyzed the space $T_{\widetilde{h},h}G(L^2({\bf T}^2))$, $\widetilde{h}=\frac{h}{\sqrt{\eps}}$,
and identified it with a weighted space $H_{\Phi_3,h}$ of holomorphic functions
defined in a region where $1\ll \abs{\Im z_1}\ll \frac{1}{\sqrt{\eps}}$,
$\abs{\Im z_2} \ll 1$. We shall now see that restricting the attention to a region where
$\abs{\Im z_1} \gg \eps^{-1/6}$, we can identify this space with a weighted space of holomorphic functions on the
$T_{h,h}$--transform side. Specifically, when studying $T_{h,h}G(L^2({\bf T}^2))$ as a weighted space, we shall
show that the region $\abs{\Im z_1} \gg \eps^{-\frac{1}{6}}$ on the $T_{\widetilde{h},h}$--side corresponds to a region
$\abs{\Im z_1} \gg \eps^{\frac{1}{3}}$ on the $T_{h,h}$--side.

All the work will concern the variables of index 1, and therefore we
shall restrict the attention to the
one-dimensional situation for a while and consider
(as appears also in the discussion of the second microlocalization in Chapter 16 of~\cite{Sj82}),
$$
T_h T^{-1}_{\widetilde{h}} u(x) = C h^{-\frac{3}{4}} \widetilde{h}^{-\frac{1}{4}} \int\!\!\!\int
e^{-\frac{1}{2h}(x-t)^2+\frac{1}{2\widetilde{h}}(t-y)^2} u(y)\,dy\,dt.
$$
Eliminating the $t$--integration by exact stationary phase, we get
\begeq
\label{eq52}
T_h T^{-1}_{\widetilde{h}} u(x) = C(1+{\cal O}(\sqrt{\eps}))h^{-\frac{1}{4}} \widetilde{h}^{-\frac{1}{4}}
\int e^{\frac{1}{2h}\frac{\sqrt{\eps}}{1-\sqrt{\eps}}(x-y)^2} u(y)\,dy.
\endeq
Let us consider first the operator $T_h T^{-1}_{\widetilde{h}}$ as a
map from $H_{\Phi_0,\widetilde{h}}$ to $H_{\Phi_0,h}$ with
$\Phi_0(x)=\frac{1}{2} (\Im x)^2$. Considering the reduced kernel of (\ref{eq52}), we then want to look at
\begin{eqnarray}
\label{eq53}
& &-\frac{1}{2} (\Im x)^2 + \Re \frac{1}{2} \frac{\sqrt{\eps}}{1-\sqrt{\eps}}(x-y)^2+\frac{\sqrt{\eps}}{2} (\Im y)^2
\\ \nonumber
& = & \frac{1}{2} \frac{\sqrt{\eps}}{1-\sqrt{\eps}} (\Re x-\Re y)^2-
\frac{1}{2}\left((\Im x)^2+\frac{\sqrt{\eps}}{1-\sqrt{\eps}}(\Im x-\Im y)^2-\sqrt{\eps} (\Im y)^2\right) \\ \nonumber
& = & \frac{1}{2}\frac{\sqrt{\eps}}{1-\sqrt{\eps}}\left(\Re x - \Re y\right)^2-
\frac{1}{2(1-\sqrt{\eps})} \left(\sqrt{\eps}\Im y-\Im x\right)^2.
\end{eqnarray}
This means that we can choose the integration contour $\Re y = \Re x$
in (\ref{eq52}), and using Schur's lemma we see that
\begin{eqnarray*}
\norm{T_h T_{\widetilde{h}}^{-1}}_{H_{\Phi_0,\widetilde{h}}\rightarrow H_{\Phi_0,h}} & = & {\cal O}(1) h^{-\frac{1}{4}} \widetilde{h}^{-\frac{1}{4}}
\left(\int e^{-\frac{1}{2h(1-\sqrt{\eps})}(\sqrt{\eps}{\rm Im} y-{\rm Im} x)^2}\,d\Im y\right)^{1/2} \\ \nonumber
& & \times \left(\int e^{-\frac{1}{2h(1-\sqrt{\eps})}(\sqrt{\eps}{\rm Im} y-{\rm Im} x)^2}\,d\Im x\right)^{1/2} \\ \nonumber
& = & {\cal O}(1) h^{-\frac{1}{4}} \widetilde{h}^{-\frac{1}{4}} \left(\frac{h}{\eps}\right)^{\frac{1}{4}}h^{\frac{1}{4}}
={\cal O}(1) \eps^{\frac{1}{8}-\frac{1}{4}}={\cal O}(1) \eps^{-\frac{1}{8}}.
\end{eqnarray*}
Here $H_{\Phi_0,\widetilde{h}}=H_{\sqrt{\eps}\Phi_0,h}$. We notice
that the factor $\eps^{-1/8}$ here represents a loss, since we know that
$T_h T^{-1}_{{\widetilde h}}: H_{\Phi_0,\widetilde{h}}\rightarrow
H_{\Phi_0,h}$ is unitary, modulo exponentially small errors. The loss
is due to the fact that here we are using contour integrals as a
preparation for the next case when the weights are no longer the standard quadratic ones.

\vskip 2mm
\noindent
Let us pass from $\sqrt{\eps}\Phi_0(y)=\frac{\sqrt{\eps}}{2}(\Im y)^2$ to
\begeq
\label{eq54}
\Phi_3(y)=\frac{\sqrt{\eps}}{2}(\Im y)^2+\Phi_4(y)
\endeq
in (\ref{eq36}), (\ref{eq37}), with
\begeq
\label{eq55}
\partial_{{\rm Re} y}^k \partial_{{\rm Im} y}^l \Phi_4(y)={\cal O}\left(\frac{\sqrt{\eps}}{\abs{\Im y}^{1+l}}\right),\quad
1\ll \abs{\Im y} \ll \frac{1}{\sqrt{\eps}}.
\endeq
(Here we continue to neglect the dependence on the variable of index 2.) When defining $T_h T^{-1}_{{\widetilde{h}}}$ on $H_{\Phi_3,h}$
we need to choose the integration contour in (\ref{eq52}) passing through the critical point of
\begin{eqnarray}
\label{eq55.1}
& & y\mapsto \Re \frac{1}{2}\frac{\sqrt{\eps}}{1-\sqrt{\eps}}(x-y)^2+\Phi_3(y)
\\ \nonumber
& & =\frac{1}{2}\frac{\sqrt{\eps}}{1-\sqrt{\eps}}(\Re x-\Re y)^2-\frac{\eps}{2(1-\sqrt{\eps})}
\left(\Im y-\frac{\Im x}{\sqrt{\eps}}\right)^2+\frac{1}{2}(\Im x)^2+\Phi_4(y).
\end{eqnarray}

We shall now discuss the estimates on the critical point $y(x)$ in (\ref{eq55.1}). Using (\ref{eq55}),
we see first that the criticality with respect to $\Im y$ means that
\begeq
\label{eq55.2}
\Im y(x)+{\cal O}\left(\frac{1}{\sqrt{\eps}(\Im y(x))^2}\right)=\frac{\Im x}{\sqrt{\eps}}.
\endeq
Working in a region where
\begeq
\label{eq55.3}
\frac{1}{\sqrt{\eps}\abs{\Im y}^2} \ll \abs{\Im y} \quad \wrtext{so that}\quad
\abs{\Im y} \gg \eps^{-\frac{1}{6}},
\endeq
we then see that
\begeq
\label{eq55.4}
\Im y(x) = \frac{\Im x}{\sqrt{\eps}}+{\cal O}\left(\frac{\sqrt{\eps}}{\abs{\Im x}^2}\right).
\endeq
Considering the $\Re y$-gradient of the phase in (\ref{eq55.1}), we get
\begeq
\label{eq55.5}
\Re y(x)=\Re x+{\cal O}\left(\frac{1}{\abs{\Im y(x)}}\right),
\endeq
and in view of (\ref{eq55.4}),
\begeq
\label{eq55.6}
\Re y(x)=\Re x+{\cal O}\left(\frac{\sqrt{\eps}}{\abs{\Im x}}\right).
\endeq

\begin{prop}
The critical point $y(x)$ in {\rm  (\ref{eq55.1})} satisfies
$$
\Re y = \Re x + {\cal O}\left(\frac{\sqrt{\eps}}{\abs{\Im
x}}\right),\quad \Im y = \frac{\Im x}{\sqrt{\eps}} + {\cal
O}\left(\frac{\sqrt{\eps}}{\abs{\Im x}^2}\right),
$$
where the remainders enjoy the following symbolic estimates: for each
$k$, $l\in \nat$, we have
\begeq
\label{eq55.65}
\partial_{{\rm Re} x}^k \partial_{{\rm Im} x}^l {\cal
O}\left(\frac{\sqrt{\eps}}{\abs{\Im x}}\right) = {\cal
O}\left(\frac{\sqrt{\eps}}{\abs{\Im x}^{1+l}}\right),\quad
\partial_{{\rm Re} x}^k \partial_{{\rm Im} x}^l {\cal
O}\left(\frac{\sqrt{\eps}}{\abs{\Im x}^2}\right) = {\cal
O}\left(\frac{\sqrt{\eps}}{\abs{\Im x}^{2+l}}\right).
\endeq
\end{prop}
\begin{proof}
The proof is a rescaling argument. With $\Re x= t$, $\frac{{\rm Im}
x}{\sqrt{\eps}}=s$, $\abs{s}\gg \eps^{-1/6}$, the equations
(\ref{eq55.2}) and (\ref{eq55.5}) become, if we write $\Re y=u$, $\Im
y=v$,
\begeq
\label{eq55.7}
\cases{ u=t+f(u,v) \cr
v=s+g(u,v),\cr}
\endeq
with
$$
\partial^k_u \partial^l_v f(u,v) = {\cal
O}\left(\frac{1}{\abs{v}^{1+l}}\right),\quad \partial^k_u \partial^l_v g(u,v)={\cal
O}\left(\frac{1}{\sqrt{\eps}\abs{v}^{2+l}}\right).
$$
Assume that $s\simeq s_0$, $\abs{s_0}\gg \eps^{-1/6}$, and write
$\tilde{v}=\frac{v}{s_0}$. Then (\ref{eq55.7}) gives
\begeq
\label{eq55.8}
\cases {u=t+f(u,s_0\tilde{v}) \cr
\tilde{v}=\tilde{s}+\frac{1}{s_0} g(u,s_0\tilde{v}).\cr}
\endeq
Here we have written $\tilde{s}=s/s_0$. Now
$$
\partial^k_u \partial^l_{\tilde{v}} f(u,s_0\tilde{v})={\cal
O}\left(\frac{1}{\abs{s_0}}\right)\ll 1,
$$
$$
\partial^k_u \partial^l_{\tilde{v}}
\frac{1}{s_0}g(u,s_0\tilde{v})={\cal
O}\left(\frac{1}{\sqrt{\eps}\abs{s_0}^3}\right)\ll 1,
$$
and we conclude that $u=u(t,\tilde{s})$ and
$\tilde{v}=\tilde{v}(t,\tilde{s})$ with
$\partial_t^k \partial^l_{\tilde{s}}u={\cal O}(1)$, $\partial_t^k
\partial_{\tilde{s}}^l \tilde{v}={\cal O}(1)$. Reinjecting this
information into (\ref{eq55.8}), we get
$$
u=t+a(t,\tilde{s}),\quad \tilde{v}=\tilde{s}+\tilde{b}(t,\tilde{s}),
$$
with
$$
\partial_t^k \partial_{\tilde{s}}^l a ={\cal
O}\left(\frac{1}{\abs{s_0}}\right),\quad \partial_t^k \partial_{\tilde{s}}^l
\tilde{b}={\cal O}\left(\frac{1}{\sqrt{\eps}\abs{s_0}^3}\right).
$$
Using that $\partial_{\tilde{s}}=s_0\partial_s$, we get
$$
u=t+{\cal O}\left(\frac{1}{\abs{s}}\right),\quad \partial_t^k \partial_s^l
{\cal O}\left(\frac{1}{\abs{s}}\right)={\cal
O}\left(\frac{1}{\abs{s}^{1+l}}\right),
$$
and
$$
v=s+{\cal O}\left(\frac{1}{\sqrt{\eps}\abs{s}^2}\right),\quad \partial_t^k
\partial_s^l {\cal O}\left(\frac{1}{\sqrt{\eps}\abs{s}^2}\right)={\cal
O}\left(\frac{1}{\sqrt{\eps}\abs{s}^{2+l}}\right).
$$
The symbolic estimates (\ref{eq55.65}) follow and this completes the proof.
\end{proof}

Choosing the integration contour $\{y; \Re y=\Re y(x)\}$ in
(\ref{eq52}) passing through the critical point $y(x)$, and noticing that the $y$-Hessian
of the phase occurring in (\ref{eq55.1}) along the contour is negative
definite, we obtain that $T_h T^{-1}_{\widetilde{h}}$ becomes
a well-defined operator of norm ${\cal O}(\eps^{-\frac{1}{8}})$ from $H_{\Phi_3,h}$ to $H_{\Phi_5,h}$, where
\begin{eqnarray}
\label{eq56}
& & \Phi_5(x) \\ \nonumber
& & = {\rm vc}_y \left(\frac{1}{2}(\Im x)^2+\frac{1}{2}\frac{\sqrt{\eps}}{1-\sqrt{\eps}}(\Re x-\Re y)^2
-\frac{\eps}{2(1-\sqrt{\eps})}\left(\Im y-\frac{\Im x}{\sqrt{\eps}}\right)^2+\Phi_4(y)\right).
\end{eqnarray}
Using the estimates (\ref{eq55.4}) and (\ref{eq55.5}), we see that
\begeq
\label{eq57}
\Phi_5(x)=\frac{1}{2}(\Im x)^2+{\cal O}\left(\frac{\eps}{\abs{\Im x}}+\frac{\eps^{3/2}}{\abs{\Im x}^2}+
\frac{\eps^2}{\abs{\Im x}^4}\right).
\endeq
In view of (\ref{eq55.3}), the strictly subharmonic function $\Phi_5(x)$ is naturally defined in a region
$\eps^{1/3}\ll \abs{\Im x} \ll 1$, and therefore we get
\begeq
\label{eq57.1}
\Phi_5(x)=\frac{1}{2}(\Im x)^2+{\cal O}\left(\frac{\eps}{\abs{\Im x}}\right).
\endeq
Estimating the derivatives of $\Phi_5$ using Proposition 4.3 and adding the dependence on the variable $x_2$, we get the
following result.

\begin{prop} Let us consider the weight $\Phi_3$, defined in {\rm (\ref{eq50})} and
satisfying {\rm (\ref{eq54})}, {\rm (\ref{eq55})}. Working in a region $\eps^{-1/6}\ll \abs{\Im y_1} \ll \eps^{-1/2}$,
$\abs{\Im y_2}\ll 1$, we have
\begeq
\label{eq57.2}
T_{h,h} T_{\widetilde{h},h}^{-1}={\cal O}(1) \eps^{-\frac{1}{8}}: H_{\Phi_3,h}\rightarrow H_{\Phi_5,h},
\endeq
where the strictly plurisubharmonic function $\Phi_5$ is given in {\rm (\ref{eq56})} and
is defined for $\eps^{1/3}\ll \abs{\Im x_1}\ll 1$, $\abs{\Im x_2}\ll 1$. We have
\begeq
\label{eq57.3}
\Phi_5(x)=\frac{1}{2} (\Im x)^2+\Phi_6(x_1,\Im x_2),
\endeq
where
\begeq
\label{eq57.35}
\partial_{{\rm Re} x_1, {\rm Im} x_2}^k \partial_{{\rm Im} x_1}^l \Phi_6(x_1,\Im x_2)=
{\cal O}_{kl}\left(\frac{\eps}{\abs{\Im x_1}^{l+1}}\right)
\endeq
\end{prop}

\noindent
{\it Remark}. Let us notice that we would also obtain the weighted space $H_{\Phi_5,h}$
more directly by studying $T_h Gu$ in the $x_1$-variable, with $G$ given in
(\ref{eq18}) and $u\in L^2$. Notice also that the canonical transformation
associated to $T_{h,h} T^{-1}_{\widetilde{h},h}$ in (\ref{eq57.2}) is
$\kappa_{T_{h,h}}\circ \kappa_{T_{\widetilde{h},h}}^{-1}$, and since
in view of (\ref{eq28.5}),
$$
\kappa_{T_{\widetilde{h},h}}^{-1}(x_1,\xi_1;x_2,\xi_2)=
\left(x_1+i\frac{\xi_1}{\sqrt{\eps}},\xi_1; x_2+i\xi_2,\xi_2\right),
$$
we get, using (\ref{eq27.01}),
\begeq
\label{eq57.4}
\kappa_{T_{h,h}}\circ \kappa_{T_{\widetilde{h},h}}^{-1}(x_1,\xi_1;x_2,\xi_2)=\left(x_1+i\xi_1
\frac{1-\sqrt{\eps}}{\sqrt{\eps}},\xi_1; x_2,\xi_2\right).
\endeq
The weights $\Phi_3$ in (\ref{eq54}) and $\Phi_5$ in (\ref{eq56}) are related through the formula
\begeq
\label{eq57.5}
\kappa_{T_{h,h}}\circ \kappa_{T_{\widetilde{h},h}}^{-1}(\Lambda_{\Phi_3})=\Lambda_{\Phi_5}.
\endeq


\subsection{Microlocal Hilbert space in a full \neigh{} of $\Lambda_{1,r}$}
Let us return to the situation discussed in section 3, and recall
from (\ref{eq7}) that the action of $P_{\eps}$ on $F(L^2_{\theta}({\bf T}^2))$ is, microlocally near $\xi=0$,
equivalent to the action of
$$
F^{-1} P_{\eps} F= P'_{\eps}(x_1,hD_x;h)+R_{\eps}(x,hD_x;h)
$$
on $L^2_{\theta}({\bf T}^2)$. Recall also from (\ref{eq26}) that
$$
G^{-1}M^{-1} P'_{\eps} MG=P_{\eps}''(hD_x;h),
$$
the operator $P_{\eps}''(hD_x;h)$ being defined in (\ref{eq21}).
With $\widetilde{h}=\frac{h}{\sqrt{\eps}}$, in Proposition 4.1 we have identified
$T_{\widetilde{h},h}G(L^2({\bf T}^2))$ with a space of holomorphic
functions $H_{\Phi_3,h}$ in the region $1\ll \abs{\Im
z_1} \ll \frac{1}{\sqrt{\eps}}$, $\abs{\Im z_2}\ll 1$, with $\Phi_3$ as
in (\ref{eq37.5}), (\ref{eq37.6}). Now let us extend $\Phi_3(\cdot,\Im z_2)$ from the region $R<\abs{\Im
z_1}\ll \frac{1}{\sqrt{\eps}}$ for $R\gg 1$, to the entire domain $\abs{\Im z_1}\ll \frac{1}{\sqrt{\eps}}$
so that we still have as in Proposition 4.1,
$$
\Phi_3(z_1,\Im z_2)=\frac{\sqrt{\eps}}{2}(\Im z_1)^2+\frac{1}{2} (\Im z_2)^2+\Phi_4(z_1,\Im z_2)
$$
with
\begeq
\label{eq60}
\partial^k_{{\rm Re} z_1, {\rm Im} z_2}\partial^l_{{\rm Im} z_1}\Phi_4(z_1,\Im z_2)
= {\cal O}_{k l}\left(\frac{\sqrt{\eps}}{(R+\abs{\Im
z_1})^{1+l}}\right),\quad R\gg 1,
\endeq
and with $\Phi_3$, $\Phi_4$ still independent of $\Re z_2$. It follows that
when
$$
(x_1,\widetilde{\xi_1},x_2,\xi_2)\in \Lambda_{\Phi_3}=\left\{\left(x,\frac{2}{i}\frac{\partial \Phi_3}{\partial x}\right)
\right\},
$$
then
$\Im \xi_2=-\frac{\partial \Phi_3}{\partial {\rm Re}\,x_2}=0$.

Recall that in subsection 4.1 we have defined a microlocal Hilbert space in a fixed
\neigh{} of
$$
\bigcup_{\abs{E}<\delta_0} \Lambda_{E,1,r}\quad 0<\delta_0 \ll 1,
$$
but away from a $\sqrt{\eps}$--\neigh{} of that set as $FMG(L^2_{\theta}({\bf
T}^2))$. Here the tori $\Lambda_{E,1,r}$ have been introduced in (\ref{eq1.2}).
Having extended $\Phi_3$, we now fill the gap by replacing
$G(L^2_{\theta}({\bf T}^2))$ by $T_{\widetilde{h},h}^{-1}H_{\Phi_3,h}$, and introduce a
microlocal Hilbert space defined in a full \neigh{} of $\Lambda_{1,r}$ and given by
\begeq
\label{eq61}
FM T^{-1}_{\widetilde{h},h} H_{\Phi_3,h}.
\endeq
Here it will be understood that the elements of $H_{\Phi_3,h}$ are Floquet periodic as in (\ref{2.2}). In what follows, in
order to simplify the presentation, we shall neglect the Floquet
conditions and work under the assumption that the
elements of the weighted space $H_{\Phi_3,h}$ are $2\pi \z^2$--periodic functions. It will be clear that the discussion below will
extend to the Floquet periodic case. Also, in (\ref{eq61}) we are identifying a \neigh{} of $\Lambda_{1,r}$ with a
\neigh{} of the zero section in $T^*{\bf T}^2$ by means of the canonical transformation $\kappa_0$ in (\ref{eq0.5}).

\medskip
Microlocally near $\Lambda_{1,r}$, the action of $P_{\eps}$ on the space (\ref{eq61}) can be identified with that of
\begeq
\label{eq62}
T_{\widetilde{h},h}M^{-1} P'_{\eps}M
T^{-1}_{\widetilde{h},h}+T_{\widetilde{h},h} M^{-1} R_{\eps} M
T_{\widetilde{h},h}^{-1}=:
\widetilde{P}_{\eps}+\widetilde{R}_{\eps}
\endeq
on $H_{\Phi_3,h}$, in view of (\ref{eq7}).
The operator $\frac{1}{\eps} M^{-1} P'_{\eps} M$ is given by
(\ref{eq11}) with $\xi_2$ replaced by $hD_{x_2}$. The operator
$\frac{1}{\eps}\widetilde{P}_{\eps}$ therefore becomes, with $\mu=\sqrt{\eps}$,
\begin{eqnarray}
\label{eq63}
\frac{1}{\eps}\widetilde{P}_{\eps} & = & \frac{1}{\eps}
p(f(\xi_2),\xi_2)+g(f(\xi_2)+\mu
\widetilde{h}D_{x_1},\xi_2)(\widetilde{h}D_{x_1})^2 \\ \nonumber
& & + \left(i\widetilde{\langle{q}\rangle_2}+\frac{{\cal
O}(h)}{\eps}+\eps r\right)(x_1+i\widetilde{h}D_{x_1}, f(\xi_2)+\mu
\widetilde{h}D_{x_1},\xi_2)\\ \nonumber
& & +{\cal O}(\widetilde{h}^2)(x_1+i\widetilde{h}D_{x_1},f(\xi_2)+\mu
\widetilde{h}D_{x_1},\xi_2),
\end{eqnarray}
where we replace $\xi_2$ by $hD_{x_2}$, since the
$\widetilde{h}$--Fourier integral operator
$T^{(1)}_{\widetilde{h}}$ is a convolution operator with the associated
\cantransf{} $(y_1,\eta_1)\mapsto (y_1-i\eta_1,\eta_1)$, and similarly for
$T^{(2)}_h$. From (\ref{eq7.1}) we also find that
\begeq
\label{eq63.5}
\frac{1}{\eps} \widetilde{R}_{\eps}(x_1,\widetilde{\xi}_1,x_2,\xi_2;h)={\cal O}\left(\eps^N+h^N+
\eps^{\frac{N}{2}}\widetilde{\xi_1}^N\right).
\endeq

To study the operator in (\ref{eq62}), we take a Fourier series expansion in $x_2$ of
a general element $u\in H_{\Phi_3,h}$,
\begeq
\label{eq64}
u(x_1,x_2)=\sum_{k\in {\bf Z}} \widetilde{u}_k(x_1) e_{kh}(x_2).
\endeq
Here the functions $e_{kh}(x_2)$ have been introduced in (\ref{eq51.2}). From Proposition 4.2 we recall that
\begeq
\label{eq65}
\norm{u}^2_{H_{\Phi_3,h}}\sim \sum_{k\in {\bf Z}}
\norm{\widetilde{u}_k}^2_{H_{\Phi_1(\cdot,kh),\widetilde{h}}},
\endeq
and correspondingly,
\begeq
\label{eq66}
\frac{1}{\eps} \widetilde{P}_{\eps} u = \sum_{k\in {\bf Z}}
\frac{1}{\eps}
\widetilde{P}_{\eps}\left(x_1,\widetilde{h}D_{x_1},kh;h\right)
\widetilde{u}_k(x_1) e_{kh}(x_2).
\endeq
Therefore we have to study
\begeq
\label{eq67}
\frac{1}{\eps} \widetilde{P}_{\eps}(x_1,\widetilde{h}D_{x_1},kh;h):
H_{\Phi_1(\cdot,kh),\widetilde{h}}\rightarrow
H_{\Phi_1(\cdot,kh),\widetilde{h}}.
\endeq
Here we recall from (\ref{eq32.1}) and (\ref{eq33}) that
\begeq
\label{eq67.1}
\Phi_1(x_1,\eta_2)=\frac{1}{2}(\Im x_1)^2+\Phi_2(x_1,\eta_2),\quad
\Phi_2(x_1,\eta_2)={\cal O}\left(\frac{1}{\abs{\Im x_1}}\right)
\endeq
is defined in the region $R< \abs{\Im x_1}\ll \frac{1}{\sqrt{\eps}}$,
and when extending the definition to the domain $\abs{\Im x_1}\leq R$,
we use Proposition 4.2. Using also (\ref{eq60}), we see that the
representation (\ref{eq67.1}) holds in the entire region $\abs{\Im
x_1}\ll \frac{1}{\sqrt{\eps}}$, with
\begeq
\label{eq67.10}
\partial_{{\rm Re} x_1,\eta_2}^k \partial_{{\rm Im} x_1}^l \Phi_2(x_1,\eta_2)=
{\cal O}\left(\frac{1}{(R+\abs{\Im x_1})^{l+1}}\right),\quad R\gg 1.
\endeq

\vskip 2mm
\noindent
In the region where $\abs{\Im x_1}\gg 1$, we have, from Proposition 3.1,
\begeq
\label{eq68}
\widetilde{P}_{\eps}=T_{\widetilde{h},h}G
P_{\eps}''G^{-1}T^{-1}_{\widetilde{h},h},
\endeq
and correspondingly for 1-variable pseudodifferential operators:
\begeq
\label{eq69}
\frac{1}{\eps} \widetilde{P}_{\eps}(x_1,\widetilde{h}D_{x_1},kh;h) = \frac{1}{\eps} T_{\widetilde{h}}
G_{kh}
P_{\eps}''(\widetilde{h}D_{x_1},kh;h)G^{-1}_{kh}T^{-1}_{\widetilde{h}}.
\endeq
Here $G_{kh}$ is defined in (\ref{eq38}) and $P''_{\eps}(\widetilde{h}D_{x_1},kh;h)$ is given in (\ref{eq21}):
\begin{eqnarray}
\label{eq70}
\frac{1}{\eps}
P_{\eps}''(\widetilde{h}D_{x_1},\xi_2;h) & = & \frac{1}{\eps}
p(f(\xi_2),\xi_2)+
g(f(\xi_2)+\sqrt{\eps}\widetilde{h}D_{x_1})(\widetilde{h}D_{x_1})^2\\
\nonumber
& & + \left(i\langle{\widetilde{\langle{q}\rangle_2}\rangle_1}+{\cal
O}\left(\frac{h}{\eps}\right)+\eps
\langle{r_{\eps}}\rangle_1\right)(f(\xi_2)+\sqrt{\eps}\widetilde{h}D_{x_1},\xi_2)\\
\nonumber
& & + {\rm Op}_{\widetilde{h}} \left({\cal
O}\left(\frac{1}{\widetilde{\xi_1}^2}\right)\right)+\widetilde{R}(\widetilde{h}D_{x_1},\xi_2,\eps;h),
\end{eqnarray}
where
$$
\widetilde{R}\sim \sum_{j=2}^{\infty}\widetilde{h}^j \eps^{-j/2}
\widetilde{R}_j,\quad \widetilde{R}_j={\cal
O}\left(\frac{1}{\abs{\widetilde{\xi_1}}^{2j-2}}\right).
$$
An application of Egorov's theorem then shows that in the region where
$\abs{\Im x_1}\gg 1$, the symbol of $\frac{1}{\eps}
\widetilde{P}_{\eps}(x_1,\widetilde{h}D_{x_1},kh;h)$ restricted to
$$
\Lambda_{\Phi_1(\cdot,kh)}=\left\{\left(x_1,\frac{\partial \Phi_1(x_1,kh)}{\partial x_1}\right)\right\},
$$
can be identified with the symbol of (\ref{eq70}) restricted to
$T^*{\bf T}^2$, modulo an error ${\cal O}(\widetilde{h})$. Let us notice also
that if $(x_1,\widetilde{\xi_1})\in \Lambda_{\Phi_1(\cdot,kh)}$ then
from (\ref{eq67.1}), (\ref{eq67.10}),
\begeq
\label{eq70.01}
\Re \widetilde{\xi}_1=-\frac{\partial \Phi_1}{\partial {\rm Im} x_1}(x_1,kh)=
-\Im x_1+{\cal O}\left(\frac{1}{(R+\abs{{\rm Im} x_1})^2}\right),
\endeq
and
\begeq
\label{eq70.02}
\Im \widetilde{\xi_1}=-\frac{\partial \Phi_1}{\partial {\rm
Re} x_1}(x_1,kh)={\cal O}\left(\frac{1}{R+\abs{\Im x_1}}\right),\quad
R\gg 1,
\endeq
so that the imaginary part of the term
$g(f(\xi_2)+\sqrt{\eps}\widetilde{\xi_1},\xi_2)\widetilde{\xi_1}^2$,
occurring in the symbol in (\ref{eq63}),
restricted to $\Lambda_{\Phi_1(\cdot,\eta_2=kh)}$, is small, when
$\abs{\Im x_1}={\cal O}(1)$.

\medskip
We shall finish this section by discussing the action of the remainder in (\ref{eq62}),
$\frac{1}{\eps}\widetilde{R}_{\eps}(x_1,\widetilde{h}D_{x_1},x_2,hD_{x_2};h)$, on $H_{\Phi_3,h}$. In doing so,
we shall work, as we may, with the classical rather than the Weyl quantization. We shall study the
scalar product $(\frac{1}{\eps} \widetilde{R}_{\eps} u_k|u_{\ell})_{H_{\Phi_3,h}}$, where
$$
u_k(x_1,x_2)=\widetilde{u}_k(x_1)e_{kh}(x_2) ,\quad u_{\ell}(x_1,x_2)=\widetilde{u}_{\ell}(x_1) e_{\ell h}(x_2),
\quad k,\ell\in \z,\,\, k\neq \ell.
$$
Here we have
\begeq
\label{eq70.03}
\abs{kh}\leq \frac{1}{\widetilde{C}},\quad \abs{\ell h}\leq \frac{1}{\widetilde{C}},
\endeq
for some $\widetilde{C}\gg 1$.
Let us consider first
\begeq
\label{eq70.04}
\frac{1}{2\pi}
\int_{{\bf C}/2\pi {\bf Z}} \frac{1}{\eps}\widetilde{R}_{\eps}(x_1,\widetilde{h}D_{x_1},x_2,kh;h)\widetilde{u}_k(x_1)
e^{\frac{i}{h}(kh)x_2} e^{-\frac{i}{h}(\ell h)\overline{x_2}} e^{-\frac{2\Phi_3(x_1,{\rm Im} x_2)}{h}}\,L(dx_2),
\endeq
which is equal to
\begeq
\label{eq70.05}
\int \widehat{\frac{1}{\eps}\widetilde{R}_{\eps}}(x_1,\widetilde{h}D_{x_1},\cdot+i {\rm Im}x_2,kh;h)\widetilde{u}_k(x_1)
(\ell-k) e^{-\frac{2}{h}(\Phi_3(x_1,{\rm Im} x_2)+\frac{(k+\ell)}{2}h{\rm Im}x_2)}\,d{\rm Im} x_2,
\endeq
where $\widehat{\frac{1}{\eps} \widetilde{R}_{\eps}}(x_1,\widetilde{h}D_{x_1},\cdot+i\Im x_2,kh;h)
\widetilde{u}_k(x_1)(\ell-k)$
is the Fourier coefficient of
$$
\real/2\pi \z \ni \Re x_2\mapsto \frac{1}{\eps} \widetilde{R}_{\eps}(x_1,\widetilde{h}D_{x_1},\Re x_2+i\Im x_2,kh;h)
\widetilde{u}_k(x_1)
$$
at the point $\ell-k$, and is therefore equal to the Fourier coefficient of
\begeq
\label{eq70.06}
\real/2\pi \z \ni \Re x_2\mapsto \frac{1}{\eps} \widetilde{R}_{\eps}(x_1,\widetilde{h}D_{x_1},\Re x_2,kh;h)
\widetilde{u}_k(x_1)
\endeq
at the same point times $e^{(k-\ell){\rm Im} x_2}$. It follows that (\ref{eq70.04})
is equal to
\begeq
\label{eq70.07}
\left(\widehat{\frac{1}{\eps}\widetilde{R}_{\eps}}(x_1,\widetilde{h}D_{x_1},\cdot,kh;h)
\widetilde{u}_k(x_1)\right)(\ell-k) \int e^{-\frac{2}{h}\left(\Phi_3(x_1,{\rm Im} x_2)+\ell h {\rm Im x_2}\right)}\,
d{\rm Im} x_2,
\endeq
and evaluating the integral in (\ref{eq70.07}) by the method of stationary phase, as in subsection 4.2, we get
\begeq
\label{eq70.08}
\left(\widehat{\frac{1}{\eps}\widetilde{R}_{\eps}}(x_1,\widetilde{h}D_{x_1},\cdot,kh;h)
\widetilde{u}_k(x_1)\right)(\ell-k) h^{\frac{1}{2}} a_{\Phi_3(x_1,\cdot)}(\ell h;h)
e^{\frac{2}{h}{\cal L}\Phi_3(x_1,-\ell h)}.
\endeq
Here the amplitude $a_{\Phi_3(x_1,\cdot)}$ is as in (\ref{eq45}).

When estimating the first factor in (\ref{eq70.08}), we recall that as
in~\cite{Sj82}, modulo an error that is ${\cal
O}(e^{-1/C\widetilde{h}})$, $C>0$, we may write
\begin{eqnarray}
\label{eq70.081}
& & \frac{1}{\eps}\widetilde{R}_{\eps}(x_1,\widetilde{h}D_{x_1},x_2,kh;h)\widetilde{u}_k(x_1)\\ \nonumber
= & & \frac{1}{2\pi \widetilde{h}}\int\!\!\!\int e^{\frac{i}{\widetilde{h}}(x_1-y_1)\widetilde{\xi}_1} \frac{1}{\eps}
\widetilde{R}_{\eps}(x_1,\widetilde{\xi}_1,x_2,kh;h)\chi(x_1-y_1)\widetilde{u}_k(y_1)\,dy_1\,d\widetilde{\xi}_1,
\end{eqnarray}
where $\chi$ is a suitable cutoff in a \neigh{} of $0$, and in (\ref{eq70.081}) we choose a good contour adapted to
the weight $\Phi_1(\cdot,\ell h)$ and given by
$$
\widetilde{\xi}_1=\frac{2}{i}\frac{\partial \Phi_1(x_1,\ell h)}{\partial x_1}+iC\overline{(x_1-y_1)},\quad C\gg 1.
$$
It follows, using also (\ref{eq63.5}) and (\ref{eq67.1}) that the absolute value of the kernel of
$$
e^{-\Phi_1(\cdot,\ell h)/\widetilde{h}}\frac{1}{\eps}\widetilde{R}_{\eps} e^{\Phi_1(\cdot,\ell h)/\widetilde{h}}
$$
does not exceed
$$
\frac{{\cal O}(1)}{\widetilde{h}} e^{-\abs{x_1-y_1}^2/\widetilde{h}}\left(\eps^N+h^N+\eps^{\frac{N}{2}}
\abs{\Im x_1}^N\right),
$$
and since $\Phi_1(x_1,\eta_2)$ is defined for $\abs{\Im x_1}\leq \frac{1}{R\sqrt{\eps}}$, $R\gg 1$, it follows
that the $H_{\Phi_1(\cdot,\ell h),\widetilde{h}}$--norm of (\ref{eq70.081}) does not exceed, uniformly in
$x_2$, $\abs{\Im x_2}\ll 1$,
$$
{\cal O}\left(\eps^N+h^N+\frac{1}{R^N}\right)\norm{\widetilde{u}_k}_{H_{\Phi_1(\cdot,\ell h),\widetilde{h}}}.
$$
Shifting also the contour of integration in $x_2$, we conclude that
the $H_{\Phi_1(\cdot,\ell h),\widetilde{h}}$--norm of
\begeq
\label{eq70.082}
x_1\mapsto \left(\widehat{\frac{1}{\eps}\widetilde{R}_{\eps}}(x_1,\widetilde{h}D_{x_1},\cdot,kh;h)
\widetilde{u}_k(x_1)\right)(\ell-k)
\endeq
can be estimated by
\begeq
\label{eq70.083}
{\cal O}\left(\eps^N+h^N+\frac{1}{R^N}\right) e^{-\abs{k-\ell}/{\cal O}(1)}
\norm{\widetilde{u}_k}_{H_{\Phi_1(\cdot,\ell h),\widetilde{h}}}.
\endeq
Combining (\ref{eq70.08}), (\ref{eq70.083}), and Proposition 4.2 we
see that the scalar product
\begeq
\label{eq70.09}
(\frac{1}{\eps} \widetilde{R}_{\eps} u_k|u_{\ell})_{H_{\Phi_3,h}}
\endeq
can be estimated by
\begeq
\label{eq70.095}
{\cal O}\left(\eps^N+h^N+\frac{1}{R^N}\right) e^{-\abs{k-\ell}/{\cal O}(1)}e^{\frac{h}{2}(\ell^2-k^2)}
\norm{\widetilde{u}_k}_{H_{\Phi_1,(\cdot,\ell h),\widetilde{h}}}
\norm{\widetilde{u}_l}_{H_{\Phi_1,(\cdot,\ell h),\widetilde{h}}}.
\endeq
Here when considering $\widetilde{u}_k$, we want to replace $\Phi_1(\cdot,\ell h)$ by $\Phi_1(\cdot,kh)$, and
according to (\ref{eq67.10}), we can do it at the expense of the exponential factor
$\exp({\cal O}(1)\frac{\sqrt{\eps}\abs{k-\ell}}{R+\abs{{\rm Im} x_1}})$,
which is permissible due to presence of the factor $\exp(-\abs{k-\ell}/{\cal O}(1))$ in (\ref{eq70.095}).
Taking into account also (\ref{eq70.03}) and (\ref{eq51.3}), we may summarize this discussion in the following result.

\begin{prop}
Assume that $k,\ell \in \z$ are such that {\rm (\ref{eq70.03})} holds, and make the assumption {\rm (\ref{eq63.5})}.
Then the scalar product
$$
\left(\frac{1}{\eps}\widetilde{R}_{\eps} u_k|u_{\ell}\right)_{H_{\Phi_3,h}},
$$
where
$$
u_{k}(x_1,x_2)=C h^{-\frac{1}{4}} e^{-\frac{(kh)^2}{2h}} e^{\frac{i}{h}(kh) x_2} \widetilde{u}_k(x_1),\quad
\widetilde{u}_k(x_1)\in H_{\Phi_1(\cdot,kh),\widetilde{h}},\,\,C>0,
$$
and $u_{\ell}$ is defined similarly, can be estimated by
\begeq
\label{eq70.096}
{\cal O}\left(\eps^N+h^N+\frac{1}{R^N}\right)e^{-\abs{k-\ell}/{\cal O}(1)}
\norm{u_k}_{H_{\Phi_3,h}} \norm{u_{\ell}}_{H_{\Phi_3,h}},\quad R\gg 1.
\endeq
\end{prop}

\section{Global Hilbert space and spectral asymptotics for $P_{\eps}$}
\setcounter{equation}{0}

\subsection{Behavior of the Diophantine weight near $\Lambda_{1,r}$}

Let us recall from the introduction that our spectral parameter $z$ varies in a rectangle of the form
$$
\abs{\Re z}<\frac{\eps}{{\cal O}(1)},\quad \abs{\frac{\Im z}{\eps}-F_0}<\frac{1}{{\cal O}(1)},
$$
where $F_0\in Q_{\infty}(\Lambda_{1,r})$ satisfies (\ref{eq0.41.1}), (\ref{eq0.41.4}),
(\ref{eq0.42}), and (\ref{eq0.44}). Recall also that we assume for simplicity that $L=2$ in (\ref{eq0.41.02})
and $L'=1$ in (\ref{eq0.41.3.5}).

In the absence of rational tori corresponding to the energy level
$(0,\eps F_0)$, the global weight that we used in~\cite{HiSjVu} when
away from a small but fixed \neigh{} of $\cup_{j=1}^2 \Lambda_{j,d}$,
was coming from an averaging procedure along the $H_p$--flow, and it
is the weight that we should use in the present case,
also when away from a \neigh{} of $\Lambda_{1,r}$. Following~\cite{HiSjVu}, we shall now recall the definition
of the weight in question.

Let $0\leq K\in C^{\infty}_0(\real)$ be even and such that $\int
K(t)\,dt=1$. When $T>0$, we introduce the smoothed out flow average of $q$,
\begeq
\label{eq71}
\langle{q}\rangle_{T,K}=\int K_T(t) q\circ \exp(tH_p)\,dt,\quad
K_T(t)=\frac{1}{T} K\left(\frac{t}{T}\right),
\endeq
the standard flow average in (\ref{eq016}) corresponding to taking $K=1_{[-1/2,1/2]}$.
Let $G_T$ be an analytic function defined near $p^{-1}(0)\cap \real^4$, such that
\begeq
\label{eq72}
H_p G_{T}=q-\langle{q}\rangle_{T,K}.
\endeq
As in~\cite{HiSjVu}, we solve (\ref{eq72}) by setting
\begeq
\label{eq73}
G_T=\int T J_T(-t) q\circ \exp(tH_p)\,dt,\quad J_T(t)=\frac{1}{T}
J\left(\frac{t}{T}\right),
\endeq
where the function $J$ is compactly supported, smooth away from 0, and with
\begeq
\label{eq74}
J'(t)=\delta(t)-K(t).
\endeq

The behavior of $G_T$ near the Diophantine tori $\Lambda_{j,d}$, $j=1$, $2$, as
$T\rightarrow \infty$, has been analyzed in~\cite{HiSjVu}.
We shall now consider the behavior of $G_T$ near $\Lambda_{1,r}$. Passing to the torus side by means
of the canonical transformation in (\ref{eq0.5}) and composing $p=p(\xi)$ in (\ref{eq2.1})
with $\kappa_M$ in (\ref{eq25}), we may reduce ourselves to the case when
\begeq
\label{eq74.5}
p(\xi_1,\xi_2)=p(f(\xi_2),\xi_2)+g(\xi_1+f(\xi_2),\xi_2)\xi_1^2,\quad f(0)=0,
\endeq
where $g(0,0)>0$. The expression (\ref{eq73}) gives
$$
G_T(x,\xi)=\int J\left(-\frac{t}{T}\right) q(x+tp'(\xi),\xi)\,dt,
$$
and expanding $q(\cdot,\xi)$ in a Fourier series, we get
\begeq
\label{eq75}
G_T(x,\xi) = \sum_{k=(k_1,k_2)\neq 0, k\in {\bf Z}^2} T \widehat{J}(T p'(\xi)\cdot
k)\widehat{q}(k,\xi) e^{ix\cdot k},
\endeq
since it follows from (\ref{eq74}) and the fact that $K$ is even that $\widehat{J}(0)=0$. Here $\widehat{q}(k,\xi)$ are the Fourier coefficients
of $q(x,\xi)$ and $\widehat{J}(\tau)=\int e^{-it\tau} J(t)\,dt$ is the Fourier transform of $J$.

We write
\begeq
\label{eq76}
G_T(x,\xi) = \sum_{k_2\neq 0}T \widehat{J}(T p'(\xi)\cdot
k)\widehat{q}(k,\xi) e^{ix\cdot k} + \sum_{k_2=0} T \widehat{J}(T p'(\xi)\cdot k)\widehat{q}(k,\xi) e^{ix\cdot k}
= {\rm I}+{\rm II},
\endeq
with the natural definitions of ${\rm I}$ and ${\rm II}$.
When estimating ${\rm I}$, we notice that when $k_2\neq 0$, $\abs{p'(\xi)\cdot k}\geq
\abs{p'_{\xi_2} k_2}-C\abs{\xi_1}\abs{k_1}\geq 1/2$, $C>0$, provided
that $2C\abs{\xi_1}\abs{k_1}\leq 1$. (Here for notational simplicity we assume that
$\abs{p'_{\xi_2}}\geq 1$.) Let now $0\leq \chi\in
C^{\infty}_0((-1,1))$ be such that $\chi=1$ on $[-1/2,1/2]$ and write,
using also (\ref{eq74}),
\begin{eqnarray}
\label{76.5}
{\rm I}& = & \sum_{k_2\neq 0} \chi(2C\abs{\xi_1}\abs{k_1}) T \widehat{J}(T
p'(\xi)\cdot k) \widehat{q}(k,\xi)e^{ix\cdot k} \\ \nonumber
& & + \sum_{k_2\neq 0} (1-\chi(2C\abs{\xi_1}\abs{k_1}))T \widehat{J}(T
p'(\xi)\cdot k) \widehat{q}(k,\xi)e^{ix\cdot k} \\ \nonumber
& = & \sum_{k_2\neq 0} \chi(2C\abs{\xi_1}\abs{k_1})
\frac{1-\widehat{K}(T p'(\xi)\cdot k)}{i p'(\xi)\cdot k}
\widehat{q}(k,\xi) e^{ix\cdot k} \\ \nonumber
& & +\sum_{k_2\neq 0} (1-\chi(2C\abs{\xi_1}\abs{k_1}))T \widehat{J}(T
p'(\xi)\cdot k) \widehat{q}(k,\xi)e^{ix\cdot k}.
\end{eqnarray}
It is easy to see that
\begeq
\label{eq76.6}
{\rm I}={\cal
O}\left(1+T\abs{\xi_1}^{\infty}\right),\quad T\geq 1.
\endeq
When considering the contribution coming from ${\rm II}$,  we notice that
\begin{eqnarray}
\label{eq78}
{\rm II} & = & \sum_{k_2=0, k_1\neq 0} T\widehat{J}(T p'_{\xi_1}k_1) e^{ix_1 k_1}
\widehat{q}(k,\xi) \\ \nonumber
& = & \sum_{k_2=0, k_1\neq 0}
\frac{1-\widehat{K}(Tp'_{\xi_1}k_1)}{ip'_{\xi_1}k_1} e^{ix_1 k_1} \widehat{q}(k,\xi),
\end{eqnarray}
and therefore, since $\abs{p'_{\xi_1}}\sim \abs{\xi_1}$, in view of
(\ref{eq74.5}), we get uniformly in $T\geq 1$,
\begeq
\label{eq79}
{\rm II}={\cal O}(1)\frac{1}{\abs{\xi_1}}.
\endeq
Combining (\ref{eq79}) with the bound ${\rm II}={\cal O}(T)$, we get
\begeq
\label{eq80}
{\rm II}={\cal O}(1)\frac{T}{T\abs{\xi_1}+1}.
\endeq

\begin{prop}
Let $G_T$ be defined in {\rm (\ref{eq73})}, {\rm (\ref{eq74})}, so that it satisfies {\rm
(\ref{eq72})}. Assume that near $\xi=0$ we have
$$
p(\xi_1,\xi_2)=p(f(\xi_2),\xi_2)+g(\xi_1+f(\xi_2),\xi_2)\xi_1^2,\quad f(0)=0,\quad g(0,0)>0.
$$
Then
\begeq
\label{eq81}
G_T(x,\xi)={\cal
O}\left(1+T\abs{\xi_1}^{\infty}+\frac{T}{T\abs{\xi_1}+1}\right),\quad T\geq 1.
\endeq
\end{prop}

\subsection{Global Hilbert space and the reference operators}
In the first part of this subsection, we shall construct
a global $h$-dependent Hilbert space where we shall study
resolvent bounds for $P_{\eps}$. The Hilbert space will be associated to a globally defined
IR-manifold $\Lambda_{\eps}\subset \comp^4$, which in a complex \neigh{} of $p^{-1}(0)\cap \real^4$,
away from a sufficiently small but fixed \neigh{} of
\begeq
\label{eq81.5}
\bigcup_{\abs{E}<\delta_0} \Lambda_{E,1,r} \quad 0<\delta_0\ll 1,
\endeq
and away from a small \neigh{} of $\bigcup_{j=1}^2 \Lambda_{j,d}$, will be given by
\begeq
\label{eq82}
\Lambda_{\eps}=\Lambda_{\eps G_T}:=\{ \exp(i\eps H_{G_T})(\rho);\, \rho\in
\real^4\}\subset \comp^4.
\endeq
Here the function $G_T$ has been defined in (\ref{eq73}). In view of the assumption
(\ref{eq0.44}) and Lemma 2.4 of~\cite{HiSjVu}, the imaginary part of $p_{\eps}$ in (\ref{eq002})
along $\Lambda_{\eps}$ in this region avoids the value $\eps F_0$, provided that $T$ is taken sufficiently large but fixed.


\medskip
When defining the global IR-manifold $\Lambda_{\eps}$ near the union of the Diophantine tori $\Lambda_{j,d}$,
$j=1,2$, we follow the procedure of~\cite{HiSjVu}, implementing a
Birkhoff normal form construction there. Therefore, it only remains to
discuss the definition of $\Lambda_{\eps}$ in a full \neigh{} of
$\Lambda_{1,r}$, and how to extend it further to $\Lambda_{\eps G_T}$ in
(\ref{eq82}).

\medskip
From the discussion in section 4, we know that near (\ref{eq81.5}), on the torus side,
$H(\Lambda_{\eps})$ should agree with the microlocal Hilbert space
\begeq
\label{eq82.1}
FMT^{-1}_{\widetilde{h},h} H_{\Phi_3,h},
\endeq
introduced in (\ref{eq61}). Now let us recall from Proposition 4.4 that in the region where
\begeq
\label{eq82.2}
\eps^{-1/6}\ll \abs{\Im x_1}\ll \eps^{-1/2},\quad \abs{\Im x_2}\ll 1,
\endeq
on the $T_{\widetilde{h},h}$--transform side, we have an identification
$T^{-1}_{\widetilde{h},h}H_{\Phi_3,h} \simeq T^{-1}_{h,h}
H_{\Phi_5,h}$, with the weight $\Phi_5$ having the properties described in
(\ref{eq57.3}), (\ref{eq57.35}). Moreover, on the $T_{h,h}$--transform side, the region in (\ref{eq82.2}) corresponds to
a region where $\eps^{1/3}\ll \abs{\Im x_1}\ll 1$, $\abs{\Im x_2}\ll 1$. In this region we
may therefore identify the microlocal Hilbert space in (\ref{eq82.1})
with
\begeq
\label{eq83}
F M T^{-1}_{h,h} H_{\Phi_5,h} = M T^{-1}_{h,h} H_{\Phi_7,h},
\endeq
where the smooth strictly plurisubharmonic function $\Phi_7(x)$ is
such that
$$
\kappa_{T_{h,h}}\circ \kappa_M^{-1}\circ \kappa_{\eps}\circ \kappa_M \circ
\kappa^{-1}_{T_{h,h}}\left(\Lambda_{\Phi_5}\right) = \Lambda_{\Phi_7}.
$$
Here $\kappa_{T_{h,h}}: (y,\eta)\mapsto (y-i\eta,\eta)$ is the canonical transformation associated to the
Bargmann transform $T_{h,h}$ on ${\bf T}^2$, given in
(\ref{eq27}). The transform $\kappa_{\eps}$ corresponding to the
operator $F$ has been introduced in (\ref{eq5}).

The transformation $\kappa_M^{-1}\circ \kappa_{\eps}\circ \kappa_M$ is
${\cal O}(\eps)$--close to the identity in the $C^{\infty}$--sense and
hence it follows from Proposition 4.4 that
\begeq
\label{eq84}
\Phi_7(x)=\Phi_0(x)+\Phi_8(x),\quad \Phi_0(x)=\frac{1}{2}(\Im x)^2,
\endeq
where the perturbation $\Phi_8(x)$ satisfies
\begeq
\label{eq84.5}
\partial^k_{{\rm Re} x_1, x_2} \partial_{{\rm Im} x_1}^l \Phi_8(x) =
{\cal O}_{kl}\left(\frac{\eps}{\abs{{\rm Im} x_1}^{1+l}}\right).
\endeq
In particular, the Hessian of $\Phi_7$ is uniformly bounded in a region
where $\eps^{1/3}\ll \abs{\Im x_1}\ll 1$, $\abs{\Im x_2}\ll 1$.

\medskip
We conclude that near (\ref{eq81.5}) but away from an ${\cal
O}(\eps^{1/3})$-\neigh{} of that set, we should choose
\begeq
\label{eq85}
\Lambda_{\eps}=\kappa_0^{-1}\circ \kappa_{\eps}\circ \kappa_M \circ
\kappa_{T_{h,h}}^{-1}(\Lambda_{\Phi_5}) = \kappa_0^{-1}\circ \kappa_M
\circ \kappa^{-1}_{T_{h,h}}\left(\Lambda_{\Phi_7}\right)
\endeq
where $\kappa_0$ is the action-angle transform defined in (\ref{eq0.5}).

\vskip 2mm
\noindent
We shall now glue the manifolds $\Lambda_{\eps G_T}$ in (\ref{eq82}) and $\Lambda_{\eps}$ in (\ref{eq85}). To that
end, from subsection 5.1 we recall that we have simplified the symbol
$p$ in (\ref{eq2.1}) by composing it with the transformation $\kappa_M$
in (\ref{eq25}). Hence
$$
\Lambda_{\eps G_T}=\kappa_0^{-1}\circ \kappa_M \left(\Lambda_{\eps
G_T\circ \kappa_0^{-1}\circ \kappa_M}\right),
$$
where $G_T:=G_T\circ \kappa_0^{-1}\circ \kappa_M$ is given in
Proposition 5.1. Recall next for example from~\cite{GeSj} that if $\Phi_{d}$ is such that
$\kappa_{T_{h,h}}(\Lambda_{\eps G_T})=\Lambda_{\Phi_{d}}$, then
\begeq
\label{eq87}
\Phi_{d}(x)=\Phi_0(x)+\eps G_T(\Re x,-\Im x)+{\cal O}(\eps^2
\abs{\nabla G_T}^2).
\endeq
Let $\chi=\chi(\Im x_1)\in C^{\infty}_0$, $0\leq \chi\leq 1$, be a standard cut-off function in a sufficiently small
but fixed \neigh{} of $0$, and consider
\begeq
\label{eq87.5}
\widetilde{\Phi}(x)=\chi(\Im x_1)\Phi_7(x)+(1-\chi(\Im x_1))\Phi_{d}(x).
\endeq
The function $\widetilde{\Phi}$ is strictly plurisubharmonic in a region
$\eps^{1/3}\ll \abs{\Im x_1}\leq \frac{1}{{\cal O}(1)}$, $\abs{\Im x_2}\leq \frac{1}{{\cal O}(1)}$. Moreover, it follows
from (\ref{eq84}), (\ref{eq84.5}), (\ref{eq87}), and Proposition 5.1 that
\begeq
\label{eq87.6}
\widetilde{\Phi}(x)=\Phi_0(x)+\Phi_9(x),
\endeq
where $\Phi_9$ and its derivatives satisfy the same estimates as
$\Phi_8(x)$ in (\ref{eq84.5}). It follows that in a fixed \neigh{} of
the set in (\ref{eq81.5}) but away from its $\eps^{1/3}$--\neigh{}, the IR-manifold $\Lambda_{\eps}$
is defined as
\begeq
\label{eq87.7}
\Lambda_{\eps}=\kappa_0^{-1}\circ \kappa_M \circ \kappa_{T_{h,h}}^{-1} \left(\Lambda_{\widetilde{\Phi}}\right),
\endeq
and we need to fill the remaining gap. To that end, it will be
convenient to go back to (\ref{eq82.1}) and to work on the
$T_{\widetilde{h},h}$-transform side. Let us recall the relation (\ref{eq57.5}) between
the weights $\Phi_3$ and $\Phi_5$,
$$
\kappa_{T_{h,h}}\circ \kappa^{-1}_{T_{\widetilde{h},h}}(\Lambda_{\Phi_3})=\Lambda_{\Phi_5},
$$
with the transform $\kappa_{T_{h,h}}\circ \kappa^{-1}_{T_{\widetilde{h},h}}$ defined in (\ref{eq57.4}). Corresponding
to the weight $\widetilde{\Phi}$ in (\ref{eq87.5}), on the $T_{h,h}$--transform side, we introduce a weight
$\widehat{\Phi}(x)$ on the $T_{\widetilde{h},h}$-transform side given by the analogous relation
\begeq
\label{eq87.8}
\kappa_{T_{h,h}}\circ \kappa_{T_{\widetilde{h},h}}^{-1}\left(\Lambda_{\widehat{\Phi}}\right)=\Lambda_{\widetilde{\Phi}}.
\endeq
We have
\begeq
\label{eq87.9}
\widehat{\Phi}(x)=\frac{\sqrt{\eps}}{2}(\Im x_1)^2+\frac{1}{2}(\Im x_2)^2+\Phi_{10}(x),
\endeq
where $\Phi_{10}$ and its derivatives satisfies the same estimates as
$\Phi_4$ in Proposition 4.1. Moreover, in a region where
$\eps^{-1/6}\ll \abs{\Im x_1}\ll \eps^{-1/2}$, $\abs{\Im x_2}\ll 1$,
the weight $\widehat{\Phi}$ is an ${\cal O}({\sqrt{\eps}})$-perturbation of
$\Phi_3$, and as such it extends to the entire region $\abs{\Im
x_1}\ll \eps^{-1/2}$, $\abs{\Im x_2}\ll 1$, in the same way as in subsection 4.4.

The definition of $\Lambda_{\eps}\subset \comp^4$ in a full \neigh{}
of $\Lambda_{1,r}$, including the gluing region, is then as follows,
\begeq
\label{eq89}
\Lambda_{\eps}=\kappa_0^{-1}\circ \kappa_M \circ \kappa_{T_{\widetilde{h},h}}^{-1}\left(\Lambda_{\widehat{\Phi}}\right).
\endeq
where the transform $\kappa_{T_{\widetilde{h},h}}$ has been defined in
(\ref{eq28.5}). Further away from $\Lambda_{1,r}$, we have
$\Lambda_{\eps}=\Lambda_{\eps G_T}$ in (\ref{eq82}),
and when approaching the Diophantine region $\Lambda_{1,d}\cup
\Lambda_{2,d}$, we define $\Lambda_{\eps}$ as in~\cite{HiSjVu}. This gives a
global definition of the IR--manifold $\Lambda_{\eps}\subset \comp^4$, which agrees with $\real^4$ outside a
bounded set.

\medskip
Let $T$ be the standard FBI--Bargmann transform, defined as in as in (\ref{eq27}), acting on
$L^2(\real^2)$, and with the associated canonical transformation
$\kappa_T: T^*\comp^2\rightarrow T^*\comp^2$, defined as in (\ref{eq27.01}). From~\cite{HiSjVu} we
know that away from a \neigh{} of the rational region, we have
\begeq
\label{eq90.1}
\kappa_T(\Lambda_{\eps})=\Lambda_{\Phi_{\eps}}:=\left\{(x,\xi)\in \comp^2\times \comp^2;\,\xi=\frac{2}{i}\frac{\partial
\Phi_{\eps}}{\partial x}\right\},
\endeq
where $\Phi_{\eps}$ is strictly plurisubharmonic with $\Phi_{\eps}-\Phi_0={\cal O}(\eps)$, $\nabla(\Phi_{\eps}-\Phi_0)={\cal O}(\eps)$,
$\Phi_0(x)=\frac{1}{2}\left(\Im x\right)^2$.
Associated to $\Lambda_{\eps}$, we then introduce a global $h$--dependent
Hilbert space $H(\Lambda_{\eps})$, which agrees with $L^2(\real^2)$ as
a set, and which is equipped with the norm
\begeq
\label{eq90.5}
\norm{u}:=\norm{T\left(1-\chi\right)u}_{H_{\Phi_{\eps}}}+\norm{T_{\widetilde{h},h}M^{-1}F^{-1}U^{-1}\chi
u}_{H_{\Phi_{3},h}}.
\endeq
Here $\chi\in C^{\infty}_0(\Lambda_{\eps})$ is a cut-off to a small
\neigh{} of the rational region, which we quantize as a Toeplitz operator on the FBI--Bargmann transform side --- see also
the following discussion in this section. The elliptic Fourier
integral operator $U$ quantizes the action-angle symplectomorphism $\kappa_0^{-1}$ in
(\ref{eq0.5}).


\bigskip
We shall now introduce a more precise description of the spectral window to which the spectral parameter $z$ is confined.
In doing so, let us recall the assumption (\ref{eq0.42}), and assume, in order to fix the ideas,
that $F_0 < \langle{q}\rangle(\Lambda_{1,r})$. Introduce a rectangle
\begeq
\label{eq90.6}
R_{\ell} = \left[-\frac{\eps}{C_0},\frac{\eps}{C_0}\right]+i\eps \left[F_0-\frac{1}{C_1},F_0+\frac{1}{C_2}\right],
\endeq
where $C_0>0$ is large enough. Moreover, we shall take $C_2>1$ so large that
\begeq
\label{eq90.7}
\frac{\Im z}{\eps} < \langle{q}\rangle(\Lambda_{1,r}),\quad z\in R_{\ell}.
\endeq
We further take $C_1>0$ so that
\begeq
\label{eq90.8}
F_0-\frac{1}{C_1} < \inf Q_{\infty}(\Lambda_{1,r}).
\endeq
Our goal now is to construct a trace class Toeplitz operator $K: H(\Lambda_{\eps})\rightarrow H(\Lambda_{\eps})$
such that the operator
$$
\frac{1}{\eps}\left(P_{\eps}+i\eps K-z\right)
$$
becomes elliptic, in the $\widetilde{h}$--pseudodifferential operator sense,
in a full \neigh{} of $\Lambda_{1,r}$, for $z$ varying in (\ref{eq90.6}). To this end, we shall
restrict the attention to the rational region.

When constructing the operator $K$, we recall that microlocally near $\Lambda_{1,r}$, the action of $P_{\eps}$ on
$H(\Lambda_{\eps})$ can be identified with the action of the operator in (\ref{eq62}) on the weighted space
$H_{\Phi_3,h}$. In what follows, as in (\ref{eq66}), (\ref{eq67}),
we shall consider the one-parameter family of operators
$\frac{1}{\eps} \widetilde{P}_{\eps}(x_1,\widetilde{h}D_{x_1},\xi_2;h)$ acting on
$H_{\Phi_1(\cdot,\xi_2),\widetilde{h}}$, where $\xi_2$ is given in (\ref{eq7.5}) and
$$
\abs{\xi_2}=\abs{h\left(k-\frac{k_0(\alpha_2)}{4}\right)-\frac{S_2}{2\pi}}\ll 1.
$$

We now claim that for $z\in \comp$ in the domain (\ref{eq90.6}) and in the region where $\abs{\xi_2}\gg \eps$, the elliptic
bound
\begeq
\label{eq91.1}
\abs{\frac{1}{\eps} \widetilde{P}_{\eps}(x_1,\widetilde{\xi}_1,\xi_2;h)-\frac{z}{\eps}}\geq \frac{1}{{\cal O}(1)}
\endeq
holds true. Here $\widetilde{\xi}_1=\frac{2}{i}\frac{\partial
\Phi_1}{\partial x_1}(x_1,\xi_2)$, so that $(x_1,\widetilde{\xi_1})\in
\Lambda_{\Phi_1(\cdot,\xi_2)}$. When verifying (\ref{eq91.1}), we
recall from subsection 4.4 that in the region where $\abs{\Im x_1}\gg 1$, the symbol of
$\frac{1}{\eps}\widetilde{P}_{\eps}(x_1,\widetilde{h}D_{x_1},\xi_2;h)$,
restricted to $\Lambda_{\Phi_1(\cdot,\xi_2)}$, is identified with the
symbol of (\ref{eq70}) restricted to $T^*{\bf T}^2$, modulo ${\cal
O}(\widetilde{h})$, and
(\ref{eq91.1}) follows by considering the imaginary part of
$\frac{1}{\eps}\left(\widetilde{P}_{\eps}(x_1,\widetilde{\xi}_1,\xi_2;h)-z\right)$, and using (\ref{eq90.7}).

It remains therefore to check (\ref{eq91.1}) in the region where
$\abs{\Im x_1}={\cal O}(1)$. Here it follows by considering the real part of
$\frac{1}{\eps}\widetilde{P}_{\eps}(x_1,\widetilde{\xi}_1,\xi_2;h)-\frac{z}{\eps}$ in (\ref{eq63}) and using that
$p(f(\xi_2),\xi_2)=h(\xi_2)\xi_2$, $h(\xi_2)>0$, and that $g(0,0)>0$,
together with (\ref{eq70.01}), (\ref{eq70.02}).

In what follows, when considering the one-parameter family
$\frac{1}{\eps}\widetilde{P}_{\eps}\left(x_1,\widetilde{h}D_{x_1},\xi_2;h\right)$,
we shall therefore restrict the attention to the quantum numbers $k\in \z$ given by the condition
\begeq
\label{eq91.2}
\xi_2=h\left(k-\frac{k_0(\alpha_2)}{4}\right)-\frac{S_2}{2\pi}={\cal O}(\eps).
\endeq
When $\widetilde{\xi}_1=\frac{2}{i}\frac{\partial \Phi_1}{\partial x_1}(x_1,\xi_2)$, using
(\ref{eq63}) together with (\ref{eq70.01}), (\ref{eq70.02}), we get,
for $\abs{\Im x_1}={\cal O}(1)$,
\begeq
\label{eq92}
\Im \frac{1}{\eps}\widetilde{P}_{\eps}(x_1,\widetilde{\xi}_1,\xi_2;h)=
\widetilde{\langle{q}\rangle_2}(\Re x_1,-\mu {\rm Im} x_1+f(\xi_2),\xi_2)
+{\cal O}\left(\frac{h}{\eps}+\eps+\frac{1}{R+\abs{\Im x_1}}\right).
\endeq
Here we recall that $\mu=\sqrt{\eps}$ and $R\gg 1$. Furthermore, as
already exploited above, in the region where $\abs{\Im x_1}\gg 1$,
the closure of the range of the imaginary part of the symbol of
$\frac{1}{\eps}\widetilde{P}_{\eps}(x_1,\widetilde{h}D_{x_1},\xi_2;h)$, restricted to $\Lambda_{\Phi_1(\cdot,\xi_2)}$,
avoids the value $F_0\in Q_{\infty}(\Lambda_{1,r})$.

For each $k\in \z$ satisfying (\ref{eq91.2}), let $0\leq r_{k}=r_{k}(\Im x_1)\in C_0^{\infty}(\real)$ be such that $r_{k}$
vanishes for $\abs{\Im x_1}\gg 1$ and such that the value $F_0$ is away from the closure of the range of
\begeq
\label{eq93}
\Im \frac{1}{\eps} \widetilde{P}_{\eps}\left(x_1,\frac{2}{i}\frac{\partial \Phi_1}{\partial x_1}(x_1,\xi_2),\xi_2;h\right)+
r_{k}(\Im x_1),
\endeq
when $\abs{\Im x_1}\leq R_1$, $R_1$ large enough. We notice that we can take $r_{k}$ to be a suitably large multiple of
some standard cutoff function. Associated with $r_{k}$ we then have a Toeplitz operator
\begeq
\label{eq94}
{\rm Top}(r_{k}): H_{\Phi_1(\cdot,\xi_2),\widetilde{h}}\rightarrow H_{\Phi_1(\cdot,\xi_2),\widetilde{h}},
\endeq
defined as in the appendix. Using the one-dimensional operators ${\rm Top}(r_k)$, we introduce
an operator ${\cal F}_{x_2}^{-1} {\rm Top}(r_{k}){\cal F}_{x_2}: H_{\Phi_3,h}\rightarrow H_{\Phi_3,h}$
given by
\begeq
\label{eq94.5}
{\cal F}_{x_2}^{-1} {\rm Top}(r_{k}){\cal F}_{x_2} u(x_1,x_2)=\sum_{\xi_2={\cal O}(\eps)}
\left({\rm Top}(r_{k})\widetilde{u}_k\right)(x_1)e_{\xi_2}(x_2),\quad u\in H_{\Phi_3,h},
\endeq
with $\xi_2$ as in (\ref{eq91.2}). Here, as in (\ref{eq64}), we have written
$$
u(x_1,x_2)=\sum_{k\in {\bf Z}} \widetilde{u}_k(x_1) e_{\xi_2}(x_2).
$$

Combining (\ref{eq91.1}) together with Proposition 4.5, and the
construction of ${\rm Top}(r_k)$, for $k\in \z$ satisfying
(\ref{eq91.2}), we conclude that for $z$ in the domain (\ref{eq90.6}), we have an elliptic estimate
\begeq
\label{eq94.6}
\left | \left | \, \left(\frac{1}{\eps}\widetilde{P}_{\eps}+\frac{1}{\eps}\widetilde{R}_{\eps}+
i {\cal F}_{x_2}^{-1}{\rm Top}(r_k){\cal
F}_{x_2}-\frac{z}{\eps}\right)u\,\right | \right |_{H_{\Phi_3,h}}\geq \frac{1}{{\cal O}(1)}
\norm{u}_{H_{\Phi_3,h}}.
\endeq
Here we are also using the basic formula relating quantization and
symbol multiplication on the FBI--Bargmann transform side, established
in Theorem 1.3 in~\cite{Sj90} (see also section 3 of~\cite{HerSjSto}).

\bigskip
\noindent
Back on the globally defined manifold $\Lambda_{\eps}$, we let now
$0\leq \chi_0\in C^{\infty}_0(\Lambda_{\eps})$
be such that $\chi_0=1$ near the rational torus and with $\supp\, \chi_0$ contained
in a small \neigh{} of the torus. We then take $0\leq {\chi}_1\in C^{\infty}_0(\Lambda_{\eps})$ supported near
$\Lambda_{1,r}$, such that $\chi_1=1$ in a \neigh{} of ${\rm supp}\,\chi_0$, and consider
\begeq
\label{eq95}
K:={\chi}_1U FM T^{-1}_{\widetilde{h},h} {\cal F}_{x_2}^{-1} {\rm Top}(r_{k}) {\cal F}_{x_2} T_{\widetilde{h},h}
M^{-1} F^{-1} U^{-1}\chi_0={\cal O}(1): H(\Lambda_{\eps})\rightarrow
H(\Lambda_{\eps}).
\endeq
Here, as in (\ref{eq90.5}), $U$ is a unitary Fourier integral operator quantizing the action-angle transformation
$\kappa_0^{-1}$ in (\ref{eq0.5}). When defining the operators corresponding to the
functions $\chi_0$ and $\chi_1$ in (\ref{eq95}), we identify $H(\Lambda_{\eps})$ with
$FMT^{-1}_{\widetilde{h},h}H_{\Phi_3,h}$ and use the Toeplitz
quantization on the FBI--Bargmann transform side.

Now it is clear that the operator in
(\ref{eq94.5}) is of trace class on $H_{\Phi_3,h}$, with its trace
class norm not exceeding
\begeq
\label{eq97}
{\cal O}\left(\frac{\eps}{h}\right)\sup_k \norm{{\rm Top}(r_{k})}_{{\rm tr}}\leq {\cal O}\left(\frac{\eps^{3/2}}{h^2}\right),
\endeq
since an application of Proposition A.1 shows that the trace class norm of the Toeplitz operator (\ref{eq94}) is
\begeq
\label{eq98}
{\cal O}\left(\frac{1}{\widetilde{h}}\right)={\cal O}\left(\frac{\sqrt{\eps}}{h}\right).
\endeq
It follows that $K$ in (\ref{eq95}) is of trace class on $H(\Lambda_{\eps})$,
its trace class norm not exceeding
$$
{\cal O}\left(\frac{\eps^{3/2}}{h^2}\right).
$$

\begin{prop}
Let us keep all the general assumptions from the introduction, and
assume that $F_0\in \cup_{\Lambda\in J} Q_{\infty}(\Lambda)$ satisfies the assumption
{\rm (\ref{eq0.41.1})}--{\rm (\ref{eq0.44})}.
Assume also that $h\ll \eps={\cal O}(h^{\delta})$, for some $\delta>0$.
Then there exists a globally defined IR-manifold $\Lambda_{\eps}\subset \comp^4$ and smooth Lagrangian tori
$\widehat{\Lambda}_{1,d}$, $\widehat{\Lambda}_{2,d}$, $\widehat{\Lambda}_{1,r}\subset \Lambda_{\eps}$ such that when
$\rho\in \Lambda_{\eps}$ is away from a small \neigh{} of $\widehat{\Lambda}_{1,d}\cup \widehat{\Lambda}_{2,d}\cup
\widehat{\Lambda}_{1,r}$ we have
\begeq
\label{eq98.1}
\abs{\Re P_{\eps}(\rho)}\geq \frac{1}{{\cal O}(1)}\quad \wrtext{or} \quad
\abs{\Im P_{\eps}(\rho)-\eps F_0}\geq \frac{\eps}{{\cal O}(1)}.
\endeq
The estimates {\rm (\ref{eq98.1})} remain valid for $\rho\in \Lambda_{\eps}$ near $\widehat{\Lambda}_{1,r}$ when away
from an ${\cal O}(\eps^{1/2})$--\neigh{} of this set. The manifold
$\Lambda_{\eps}$ is close to $\real^4$ and agrees with it outside a bounded set. We have
$$
P_{\eps}={\cal O}(1): H(\Lambda_{\eps})\rightarrow H(\Lambda_{\eps}).
$$
For $j=1,2$ there exists an elliptic Fourier integral operator
$$
U_j={\cal O}(1): H(\Lambda_{\eps})\rightarrow L^2_{\theta}({\bf T}^2)
$$
such that microlocally near $\widehat{\Lambda}_{j,d}$, $j=1,2$, we
have
$$
U_j P_{\eps}=\left(P^{(N)}_j(hD_x,\eps;h)+
R_{N+1,j}(x,hD_x,\eps;h)\right)U_j.
$$
Here $P_j^{(N)}(hD_x,\eps;h)+R_{N+1,j}(x,hD_x,\eps;h)$ is defined microlocally near $\xi=0$ in $T^*{\bf T}^2$,
the full symbol of $P_j^{(N)}(hD_x,\eps;h)$ is independent of $x$, and
$$
R_{N+1,j}(x,\xi,\eps;h)={\cal O}\left((\xi,\eps,h)^{N+1}\right).
$$
Here $N$ is arbitrarily large but
fixed. The leading symbol of $P^{(N)}_j(hD_x,\eps;h)$ is of the form
$$
p_j(\xi)+i\eps \langle{q_j}\rangle(\xi)+{\cal O}(\eps^2),
$$
with the differentials of $p_j$ and $\langle{q_j}\rangle$ being linearly
independent when $\xi=0$, $j=1,2$.

\medskip
\noindent
Furthermore, there exists a trace class Toeplitz operator
$$
K={\cal O}(1): H(\Lambda_{\eps})\rightarrow H(\Lambda_{\eps}),
$$
which has the following properties:
\begin{itemize}
\item $K$ is concentrated to the torus $\widehat{\Lambda}_{1,r}$ in the sense that if $\psi\in C^{\infty}_0(\Lambda_{\eps})$
is supported away from $\widehat{\Lambda}_{1,r}$ then
\begeq
\label{eq98.2}
\psi K=K \psi ={\cal O}(h^{\infty}): H(\Lambda_{\eps})\rightarrow H(\Lambda_{\eps}).
\endeq
\item The trace class norm of $K$ satisfies
$$
\norm{K}_{{\rm tr}}={\cal O}\left(\frac{\eps^{3/2}}{h^2}\right).
$$
\item For $\rho\in \Lambda_{\eps}$ near $\widehat{\Lambda}_{1,r}$, we have
$$
\abs{P_{\eps}(\rho)+i\eps K(\rho)-z}\geq \frac{\eps}{{\cal
O}(1)},
$$
provided that the spectral parameter $z\in \comp$ belongs to the
domain {\rm (\ref{eq90.6})}, assuming {\rm (\ref{eq90.7})}, {\rm (\ref{eq90.8})}.
\end{itemize}
\end{prop}

\noindent
{\it Remark}. It follows from the discussion preceding Proposition 5.2
that the operator $K$ enjoys better localization properties than
(\ref{eq98.2}), and is in fact concentrated to an ${\cal
O}(\eps^{1/2})$--\neigh{} of $\widehat{\Lambda}_{1,r}\subset \Lambda_{\eps}$.

\bigskip
\noindent
We shall now derive resolvent bounds for the perturbed operator
$P_{\eps}+i\eps K$ in the space $H(\Lambda_{\eps})$. To this end, let
us recall the set $E_d$, defined in Theorem 1.1, which consists of the quasi-eigenvalues $z(j,k)$,
$1\leq j\leq 2$, $k\in \z^2$, introduced in (\ref{BS1}). We introduce
an additional small parameter $0<\widetilde{\eps}={\cal
O}(h^{\delta})$ such that $\widetilde{\eps}\gg \eps^{1/2}$, $\widetilde{\eps}> h^{1/2-\delta}$.
Then it follows from Proposition 5.2 (see also
Proposition 5.1 in~\cite{HiSjVu}) that when $\rho\in \Lambda_{\eps}$
is away from an $\widetilde{\eps}$--\neigh{} of
$\widehat{\Lambda}_{1,d}\cup \widehat{\Lambda}_{2,d}\cup
\widehat{\Lambda}_{1,r}$, we have
\begeq
\label{eq98.3}
\abs{\Re P_{\eps}(\rho;h)}\geq \frac{\widetilde{\eps}}{{\cal
O}(1)}\quad \wrtext{or}\quad \abs{\Im P_{\eps}-\eps F_0}\geq
\frac{\eps\widetilde{\eps}}{{\cal O}(1)}.
\endeq

In what follows, we shall let $z\in \comp$ vary in the rectangle
\begeq
\label{eq98.4}
\left[-\frac{\eps}{C},\frac{\eps}{C}\right]+i\eps\left[F_0-\frac{\widetilde{\eps}}{C},F_0+\frac{\widetilde{\eps}}{C}\right],
\endeq
for some $C>0$ sufficiently large but fixed. Let $N_0\geq 1$ be
arbitrarily large but fixed. When $z$ in the rectangle
(\ref{eq98.4}) avoids the union of $\eps h^{N_0}/{\cal O}(1)$-neighborhoods of the $z(j,k)$'s, we would like to
show that $P_{\eps}+i\eps K-z$ is invertible and to estimate the inverse in
$H(\Lambda_{\eps})$. When doing so, to be able to exploit the Birkhoff normal form in the
Diophantine region, as in~\cite{HiSjVu}, we shall use a partition of unity involving
cutoff functions to small $h$--dependent neighborhoods of the Lagrangian tori.

In what follows we shall write that a function $a=a(\rho;h)\in
C^{\infty}(\Lambda_{\eps})$ is in the symbol class
$S_{\widetilde{\eps}}(1)$ if uniformly on $\Lambda_{\eps}$, we have
$$
\nabla^m a={\cal O}_m (\widetilde{\eps}^{-m}),\quad m\geq 0.
$$
We take a smooth partition of unity on the manifold $\Lambda_{\eps}$,
\begeq
\label{eq98.5}
1=\sum_{j=1}^2
\chi_j+\psi_{1,+}+\psi_{1,-}+\psi_{2,+}+\psi_{2,-}+\psi_3.
\endeq
Here $0\leq \chi_j\in C^{\infty}_0(\Lambda_{\eps})\cap
S_{\widetilde{\eps}}(1)$ is a cut-off function to an
$\widetilde{\eps}$--\neigh{} of $\widehat{\Lambda}_{j,d}$, $j=1$,$2$,
and as in~\cite{HiSjVu} we arrange so that
\begeq
\label{eq98.6}
[P_{\eps},\chi_j]={\cal O}(h^{(N+1)\delta}):
H(\Lambda_{\eps})\rightarrow H(\Lambda_{\eps}).
\endeq

The functions $0\leq \psi_{1,\pm}\in S_{\widetilde{\eps}}(1)$ are
such that $\pm \Re P_{\eps}\geq \widetilde{\eps}/{\cal O}(1)$ in the
support of $\psi_{1,\pm}$, respectively. Next, the functions
$0\leq \psi_{2,\pm}\in C^{\infty}_0(\Lambda_{\eps})\cap
S_{\widetilde{\eps}}(1)$ are supported in regions invariant under the
$H_p$--flow, where $\pm\left(\Im P_{\eps}-\eps F_0\right)\geq \eps
\widetilde{\eps}/{\cal O}(1)$, respectively. We also arrange so that
$\psi_{2,\pm}$ Poisson commute with $p$ on $\Lambda_{\eps}$. Here we
have written $p$ to denote the leading symbol of $P_{\eps=0}$ acting
on $H(\Lambda_{\eps})$. Finally, the function $0\leq \psi_3\in
C^{\infty}_0(\Lambda_{\eps})\cap S_{\widetilde{\eps}}(1)$ is a cut-off
to an $\widetilde{\eps}$--\neigh{} of $\widehat{\Lambda}_{1,r}$ such
that $H_p \psi_3=0$. Moreover, we can arrange that
$$
\psi_3 K=K+{\cal O}(h^{\infty}): H(\Lambda_{\eps})\rightarrow
H(\Lambda_{\eps}).
$$

At this point, we may follow the arguments of section 5
of~\cite{HiSjVu} (see also~\cite{HiSj1}) to prove, using
(\ref{eq98.3}) together with the sharp G\aa{}rding inequality, that when
\begeq
\label{eq98.7}
\left(P_{\eps}+i\eps K-z\right)u=v,\quad u\in H(\Lambda_{\eps}),
\endeq
with $z\in \comp$ varying in (\ref{eq98.4}), we have
\begeq
\label{eq98.8}
\norm{\left(1-\sum_{j=1}^2 \chi_j-\psi_3\right)u}\leq \frac{{\cal
O}(1)}{\eps \widetilde{\eps}}\norm{v}+{\cal O}(h^{\infty})\norm{u},
\endeq
provided that
\begeq
\label{eq98.81}
\frac{h}{\widetilde{\eps}^5}\leq h^{\delta}.
\endeq
Here $\norm{\cdot}$ is the norm in $H(\Lambda_{\eps})$.
Let us also remark that when establishing (\ref{eq98.8}),
following~\cite{HiSj1}, we use, in particular, that, on the operator
level,
$$
[P_{\eps},\psi_{2,\pm}]=[P_{\eps=0},\psi_{2,\pm}]+{\cal
O}\left(\frac{\eps h}{\widetilde{\eps}^2}\right)=
{\cal O}\left(\frac{h^2}{\widetilde{\eps}^4}\right)+{\cal
O}\left(\frac{\eps h}{\widetilde{\eps}^2}\right)
={\cal O}\left(\frac{\eps h}{\widetilde{\eps}^4}\right),
$$
since $h\leq \eps$. Furthermore, since $z$ belonging to (\ref{eq98.4})
is such that ${\rm
dist}(z,E_d)\geq \eps h^{N_0}/{\cal O}(1)$, directly from section 5 in~\cite{HiSjVu} we
see, using also (\ref{eq98.6}), that for $j=1$, $2$,
\begeq
\label{eq98.82}
\norm{\chi_j u}\leq \frac{{\cal O}(1)}{\eps h^{N_0}}\norm{v}+{\cal
O}(h^{(N+1)\delta-N_0-1})\norm{u},\quad (N+1)\delta-N_0-1\gg 1.
\endeq
Combining (\ref{eq98.8}) and (\ref{eq98.82}), we get
\begeq
\label{eq98.83}
\norm{(1-\psi_3)u}\leq \frac{{\cal O}(1)}{\eps h^{N_0}}\norm{v}+{\cal
O}(h^{(N+1)\delta-N_0-1})\norm{u}.
\endeq

It remains to derive an estimate for $\psi_3 u$. When doing so, we
write
\begeq
\label{eq98.84}
\left(P_{\eps}+i\eps K-z\right)\psi_3 u=\psi_3 v+[P_{\eps}+i\eps
K,\psi_3]u.
\endeq
Here
$$
[P_{\eps}+i\eps K,\psi_3]={\cal O}\left(\frac{\eps
h}{\widetilde{\eps}^4}\right): H(\Lambda_{\eps})\rightarrow
H(\Lambda_{\eps}),
$$
and using (\ref{eq98.83}) with a cut-off closer to $\widehat{\Lambda}_{1,r}$ we
see that the $H(\Lambda_{\eps})$--norm of the commutator term in the
right hand side of (\ref{eq98.84}) is controlled by
\begeq
\label{eq98.85}
{\cal O}\left(\frac{\eps h}{\widetilde{\eps}^4}\right)\frac{1}{\eps
h^{N_0}}\norm{v}+{\cal O}(h^{(N+1)\delta-N_0-1})\norm{u}=\frac{{\cal
O}(1)}{\widetilde{\eps}^4 h^{N_0-1}}\norm{v}+{\cal O}(h^{(N+1)\delta-N_0-1})\norm{u}.
\endeq
Using (\ref{eq94.6}) together with (\ref{eq98.84}) and (\ref{eq98.85}), we get
\begeq
\label{eq98.86}
\norm{\psi_3 u}\leq \frac{{\cal O}(1)}{\eps
\widetilde{\eps}^4 h^{N_0-1}}\norm{v}+{\cal O}(h^{(N+1)\delta-N_0-2})\norm{u}.
\endeq
Combining (\ref{eq98.83}) and (\ref{eq98.86}), and using also
(\ref{eq98.81}), we obtain the resolvent bounds, summarized in the following proposition.

\begin{prop}
Assume that $\widetilde{\eps}={\cal O}(h^{\delta})$, $\delta>0$, is such that $\widetilde{\eps}\gg \eps^{1/2}$ and that
{\rm (\ref{eq98.81})} holds. Let
\begeq
\label{eq98.88}
z\in \left[-\frac{\eps}{C},\frac{\eps}{C}\right]+
i\eps\left[F_0-\frac{\widetilde{\eps}}{C},F_0+\frac{\widetilde{\eps}}{C}\right],\quad C\gg 1,
\endeq
be such that ${\rm dist}(z,E_d)\geq \eps h^{N_0}/{\cal O}(1)$, for
some $N_0\geq 1$. Then, with the norm being the operator norm on
$H(\Lambda_{\eps})$, we have
\begeq
\label{eq99}
\norm{\left(P_{\eps}+i\eps K-z\right)^{-1}}\leq \frac{{\cal
O}(1)}{\eps h^{N_0}}.
\endeq
\end{prop}

\noindent
{\it Remark}. Continuing to argue as in~\cite{HiSjVu} and solving a suitable Grushin problem as in that paper,
we see that the eigenvalues of $P_{\eps}+i\eps K$ in the domain (\ref{eq98.88}) are given
by the elements of the set $E_d$ in (\ref{BS1}), modulo ${\cal O}(h^{\infty})$, their total number being
$\sim \eps \widetilde{\eps} /h^2$. We may think therefore of $P_{\eps}+i\eps K$ as a reference operator associated to the
Diophantine region, and hereafter we shall often write
\begeq
\label{eq99.1}
P_d=P_{\eps}+i\eps K.
\endeq


\bigskip
\noindent
{\it Remark}. Applying a simplified version of the argument above, we
see that in the absence of Diophantine tori corresponding to the level $(0,\eps F_0)$, the reference operator
$$
P_{\eps}+i\eps K-z:H(\Lambda_{\eps})\rightarrow H(\Lambda_{\eps})
$$
is globally invertible, with
\begeq
\label{eq99.2}
\left(P_{\eps}+i\eps K-z\right)^{-1}={\cal O}\left(\frac{1}{\eps}\right):
H(\Lambda_{\eps})\rightarrow H(\Lambda_{\eps}),
\endeq
for $z$ belonging to the rectangle (\ref{eq0.45}). Indeed, when
checking the injectivity, and hence the invertibility, of
$P_{\eps}+i\eps K-z$, together with (\ref{eq99.1}), we may use a
partition of unity of the form (\ref{eq98.5}), without the $\chi_j$'s,
with all the terms there being symbols of class $S_{\widetilde{\eps}=1}(1)$. The bound (\ref{eq99.2}) is
relevant for the proof of Theorem 1.2.

\bigskip
In the following discussion, we shall let $z\in \comp$ vary in the rectangle (\ref{eq90.6}),
so that in particular (\ref{eq90.7}) and (\ref{eq90.8}) hold.

We shall now introduce a reference operator associated to the rational region. In doing so,
we let $0\leq \widetilde{\chi}_d \in C_0^{\infty}(\Lambda_{\eps})$ be such that $\widetilde{\chi}_d=0$
in a small but fixed \neigh{} of $\widehat{\Lambda}_{1,r}$ while $\widetilde{\chi}_d=1$ away from a slightly
larger \neigh{} of this set, when restricting the attention to the region where $\abs{\Re P_{\eps}}\leq 1/{\cal O}(1)$.
Also, $\widetilde{\chi}_d$ vanishes outside of a slightly larger set of the form $\abs{\Re P_{\eps}}\leq 1/{\cal O}(1)$.
When $C>1$ is large enough, let us consider the operator
\begeq
\label{eq100}
P_r=P_{\eps}+i\eps C\widetilde{\chi}_d.
\endeq
We may view $P_r$ as a reference operator associated to the rational
region. Notice
that the trace class norm of the perturbation $P_r-P_{\eps}$ on $H(\Lambda_{\eps})$ is ${\cal O}(\eps h^{-2})$.


\bigskip
Our purpose is to study the spectrum of $P_{\eps}$ in the domain (\ref{eq98.4}) in terms of the
spectral information about the reference operators $P_d$ and $P_r$ in this
region. In parti\-cular, Proposition 5.3 gives a polynomial in $1/h$ control on the
resolvent of $P_d$, and we also know that the
eigenvalues of $P_d$ in (\ref{eq98.4}) are given,
modulo ${\cal O}(h^{\infty})$, by the elements of the set $E_d$. While
the spectral information available for the rational reference operator
$P_r$ is not going to be as precise, as a next step in our analysis,
we shall derive resolvent bounds on $P_r$ in $H(\Lambda_{\eps})$,
when $z$ in (\ref{eq90.6}) is not too close to the spectrum of this operator.

\medskip
\noindent
Using the same arguments as earlier and choosing $C>0$ sufficiently large, it is easily seen that the
operator $P_r+i\eps K-z = P_{\eps}+i\eps C \widetilde{\chi}_d+i\eps K-z$ is globally invertible on
$H(\Lambda_{\eps})$, with
\begeq
\label{eq107.5}
\left(P_{r}+i\eps K-z\right)^{-1}={\cal
O}\left(\frac{1}{\eps}\right): H(\Lambda_{\eps})\rightarrow H(\Lambda_{\eps})
\endeq
Here $z$ varies in the rectangle (\ref{eq90.6}). Write
\begeq
\label{eq108}
P_r-z=\left(P_r+i\eps K-z\right)\left(1-i\eps\left(P_r+i\eps K-z\right)^{-1}K\right).
\endeq
Proposition 5.2 together with (\ref{eq107.5}) implies that
$i\eps\left(P_r+i\eps K-z\right)^{-1}K$ is of trace
class on $H(\Lambda_{\eps})$, and the corresponding trace class norm satisfies
\begeq
\label{eq108.5}
\norm{i\eps\left(P_r+i\eps K-z\right)^{-1}K}_{{\rm tr}}=
{\cal O}\left(\frac{\eps^{3/2}}{h^2}\right).
\endeq
It follows from (\ref{eq108.5}) together with a basic estimate of~\cite{GoKr} that the holomorphic function
\begeq
\label{eq109}
D(z)={\rm det}\left(I-i\eps(P_r+i\eps K-z)^{-1}K\right),
\endeq
defined for $z$ in the rectangle (\ref{eq90.6}), satisfies
\begeq
\label{eq110}
\abs{D(z)}\leq \exp\left({\cal O}\left(\frac{\eps^{3/2}}{h^2}\right)\right).
\endeq
The zeros of the perturbation determinant $D(z)$ in the domain (\ref{eq90.6}) are precisely
the eigenvalues of $P_r$ in this region. To estimate the number of the zeros in such a
domain, with slightly increased values of $C_0$, $C_1$, and $C_2$ in (\ref{eq90.6}), it suffices,
in view of Jensen's  formula (see for example~\cite{Ma}) to establish a
lower bound on $D(z)$ at a single point $z=z_0$ in (\ref{eq90.6}). To this end we notice that the condition (\ref{eq90.8})
allows us to find $z_0$ in the domain (\ref{eq90.6}) such that
\begeq
\label{eq110.1}
\frac{\Im z_0}{\eps} < \inf Q_{\infty}(\Lambda_{1,r}).
\endeq
As before, it follows that $P_r-z_0$ is invertible with
\begeq
\label{eq111}
\left(P_r-z_0\right)^{-1}={\cal
O}\left(\frac{1}{\eps}\right): H(\Lambda_{\eps}) \rightarrow H(\Lambda_{\eps}).
\endeq
We get, using (\ref{eq108}),
\begin{eqnarray}
\label{eq112}
& &\left(I-i\eps\left(P_r+i\eps K-z_0\right)K\right)^{-1} \\ \nonumber
&=& \left(P_r-z_0\right)^{-1}\left(P_r+i\eps K-z_0\right) =
I+i\eps \left(P_r-z_0\right)^{-1}K,
\end{eqnarray}
and it follows, using (\ref{eq111}), that the absolute value of the determinant of the right hand side of
(\ref{eq112}) is ${\cal O}(\eps^{3/2}/h^2)$. Therefore,
\begeq
\label{eq113}
\abs{D(z_0)}\geq \exp\left(-{\cal O}\left(\frac{\eps^{3/2}}{h^2}\right)\right),
\endeq
and combining this bound together with (\ref{eq110}) and Jensen's formula, we conclude that the number of
eigenvalues of $P_r$ in the rectangle (\ref{eq90.6}), after an arbitrarily small decrease of the
constants $C_0$, $C_1$, and $C_2$, is
$$
{\cal O}\left(\frac{\eps^{3/2}}{h^2}\right).
$$
The proof of Theorem 1.2 is now complete, in view of the second remark following Proposition 5.3.

\medskip
Continuing with the proof of Theorem 1.1, we now come to derive resolvent estimates for the reference operator
$P_r$. Let $z$ in the rectangle (\ref{eq90.6}) be such that
\begeq
\label{eq114}
{\rm dist}(z,{\rm Spec}(P_r))\geq g(h)>0,\quad g(h)\ll \eps.
\endeq
An application of Theorem 5.1 from chapter 5 in~\cite{GoKr} together with (\ref{eq108.5}) shows that
\begeq
\label{eq115}
\norm{\left(1-i\eps\left(P_r+i\eps K-z\right)^{-1}K\right)^{-1}}\leq \frac{1}{\abs{D(z)}}
\exp\left({\cal O}\left(\frac{\eps^{3/2}}{h^2}\right)\right),
\endeq
and in view of (\ref{eq107.5}) and (\ref{eq108}), it suffices to estimate $\abs{D(z)}$
from below, away from its zeros. At this point, rather than recalling the details of the now well established
argument for that, based on Cartan's lemma (or, alternatively, on Lemma 4.3 in~\cite{Sj01}) and the Harnack
inequality together with the maximum principle, we shall merely refer to~\cite{Ma} and~\cite{Sj97},~\cite{Sj01}. We obtain
that if $z$ in the domain (\ref{eq90.6}), with increased values of the constants there, satisfies (\ref{eq114}), then
\begeq
\label{eq116}
\abs{D(z)} \geq \exp \left(-{\cal O}\left(\frac{\eps^{3/2}}{h^2}\log\frac{1}{g(h)}\right)\right).
\endeq
Combining (\ref{eq107.5}), (\ref{eq108}), (\ref{eq115}), and (\ref{eq116}), we get the following result.
\begin{prop}
Assume that $z\in \comp$ is such that
\begeq
\label{eq116.5}
\abs{\Re z}<\frac{\eps}{C},\quad \abs{\Im z-\eps F_0}< \frac{\eps}{C},\quad C\gg 1,
\endeq
with ${\rm dist}(z,{\rm Spec}(P_r))\geq g(h)$, $0<g(h)\ll \eps$. Then
\begeq
\label{eq117}
\norm{\left(P_r-z\right)^{-1}}\leq
\frac{{\cal O}(1)}{\eps}\exp\left({\cal O}\left(\frac{\eps^{3/2}}{h^2}\right)
\log \frac{1}{g(h)}\right),
\endeq
Here $P_r=P_{\eps}+i\eps C\widetilde{\chi}_d$, where $\widetilde{\chi}_d\in C^{\infty}_0(\Lambda_{\eps}; [0,1])$
is such that $\widetilde{\chi}_d=0$ near the rational torus
$\widehat{\Lambda}_{1,r}$ and $\widetilde{\chi}_d=1$ further away from this set, in the region where
$\abs{\Re P_{\eps}}\leq 1/{\cal O}(1)$.
\end{prop}

Relying upon the resolvent estimates for the reference operators $P_d$ and $P_r$, given in Propositions 5.3 and 5.4,
we shall next address the invertibility properties of $P_{\eps}-z$. This is the subject of the next subsection.

\subsection{Exponentially weighted estimates and bounds on spectral projections}
Let us recall the reference operators
\begeq
\label{eq117.5}
P_{d}=P_{\eps}+i\eps K\quad \wrtext{and}\,\,\, P_r=P_{\eps}+i\eps C\widetilde{\chi}_d,\quad C\gg 1,
\endeq
introduced in (\ref{eq99.1}) and (\ref{eq100}). In the present subsection, as in Proposition 5.3, we shall let
$z\in \comp$ vary in the domain
\begeq
\label{eq117.6}
R_s:=\left[-\frac{\eps}{C},\frac{\eps}{C}\right]+i\eps
\left[F_0-\frac{\widetilde{\eps}}{C},F_0+\frac{\widetilde{\eps}}{C}\right],\quad C\gg 1,
\endeq
where we recall that $\widetilde{\eps}$ here can be chosen so that $\widetilde{\eps} \sim \eps^{\delta}$
for any small $\delta>0$. Let us assume that
\begeq
\label{eq118}
{\rm dist}(z,{\rm Spec}(P_d)\cup {\rm Spec}(P_r))\geq g(h),
\endeq
where, as in Proposition 5.3, we take
$$
g(h)=\eps h^{N_0},
$$
for some arbitrarily large but fixed $N_0\geq 2$. We know from (\ref{eq99}) that $(P_d-z)^{-1}$ enjoys polynomial
upper bounds as a bounded operator on $H(\Lambda_{\eps})$, while Proposition 5.4 provides an
exponential estimate for the resolvent of $P_r$.

\medskip
\noindent
Next we shall introduce a reference operator associated with the elliptic region, where $\abs{\Re P_{\eps}}$ is bounded
away from zero. To this end, let $\psi\in C^{\infty}_0(\Lambda_{\eps}; [0,1])$ be such that $\psi=1$ in a region where
$\abs{\Re P_{\eps}}\leq 1/{\cal O}(1)$, and assume that $\psi$ vanishes outside of a slightly larger region of the same
form. When $C\gg 1$ and $z$ varies in the domain (\ref{eq117.6}), we see that the operator
\begeq
\label{eq118.0}
P_{\eps}+i\eps C\psi-z: H(\Lambda_{\eps})\rightarrow H(\Lambda_{\eps})
\endeq
is invertible, with
\begeq
\label{eq118.01}
\left(P_{\eps}+i\eps C\psi-z\right)^{-1}={\cal O}\left(\frac{1}{\eps}\right):
H(\Lambda_{\eps})\rightarrow H(\Lambda_{\eps}).
\endeq

\medskip
Let us consider a smooth partition of unity on the manifold $\Lambda_{\eps}$,
\begeq
\label{eq118.1}
1=\chi_r+\chi_d+\chi_0.
\endeq
Here $0\leq \chi_r\in C^{\infty}_0(\Lambda_{\eps}; [0,1])$ is $=1$ near
$\widehat{\Lambda}_{1,r}$ and $\supp\,\chi_r$ is contained in a small but fixed \neigh{} of this set.
The function $\chi_d\in C_0^{\infty}(\Lambda_{\eps}; [0,1])$ is $=1$ near $\supp\, \widetilde{\chi}_d$,
while $\chi_0\in C^{\infty}_b(\Lambda_{\eps}; [0,1])$ is such that $\supp\, \chi_0$ is contained in a region where
$\abs{\Re P_{\eps}}\geq 1/{\cal O}(1)$ and $\chi_0=1$ further away from the region where $\abs{\Re P_{\eps}}$ is small.
We furthermore arrange so that the functions $\psi$ and $\chi_0$ have disjoint supports.

Recall that $z\in R_s$ in (\ref{eq117.6}) satisfies (\ref{eq118}). As an approximation to the inverse of $P_{\eps}-z$, we consider
\begeq
\label{eq119}
R_0(z)=\left(P_r-z\right)^{-1} \chi_r+\left(P_d-z\right)^{-1} \chi_d+\left(P_{\eps}+i\eps C\psi-z\right)^{-1}\chi_0.
\endeq
Using the definitions (\ref{eq99.1}) and (\ref{eq100}) of the operators $P_d$ and $P_r$, we see that
\begeq
\label{eq120}
\left(P_{\eps}-z\right)R_0(z)=1+L,
\endeq
where
\begeq
\label{eq120.5}
L=-i\eps C\widetilde{\chi}_d\left(P_r-z\right)^{-1}\chi_r-i\eps K \left(P_d-z\right)^{-1} \chi_d-i\eps C\psi
\left(P_{\eps}+i\eps C\psi-z\right)^{-1}\chi_0.
\endeq
The key step will consist of establishing the following result.

\begin{prop}
\begeq
\label{eq120.6}
L={\cal O}(e^{-\frac{1}{\widetilde{C} h}}): H(\Lambda_{\eps})\rightarrow H(\Lambda_{\eps}),
\endeq
for some $\widetilde{C}>0$.
\end{prop}

\bigskip
\noindent
When proving Proposition 5.5, we shall introduce additional modifications of the
exponential weight corresponding to the IR-manifold $\Lambda_{\eps}$. The various modifications of the weight
will take place only in regions away from a small \neigh{} of rational torus ${\Lambda}_{1,r}$.

\medskip
We start by considering the term
\begeq
\label{eq120.7}
L_1 = -i\eps C\widetilde{\chi}_d\left(P_r-z\right)^{-1}\chi_r: H(\Lambda_{\eps})\rightarrow H(\Lambda_{\eps}),
\endeq
occurring in (\ref{eq120.5}), and notice that the compact sets $\supp\, \widetilde{\chi}_d$ and $\supp\, \chi_r$
are disjoint. From (\ref{eq90.1}) let us recall that away from $\widehat{\Lambda}_{1,r}\subset \Lambda_{\eps}$, we have
\begeq
\label{eq126.1}
\kappa_T \left(\Lambda_{\eps}\right)=\left\{(x,\xi)\in T^*\comp^2;
\xi=\frac{2}{i}\frac{\partial \Phi_{\eps}}{\partial x}\right\},
\endeq
with $\Phi_{\eps}-\Phi_0={\cal O}(\eps)$ and $\nabla\left(\Phi_{\eps}-\Phi_0\right)={\cal O}(\eps)$,
$\Phi_0(x)=(1/2) (\Im x)^2$. Here, as usual,
\begeq
\label{eq126.12}
\kappa_T(y,\eta)=(y-i\eta,\eta)=(x,\xi).
\endeq
In the following discussion, we shall often identify an open set $\Omega\subset \Lambda_{\eps}$ whose closure
is away from $\widehat{\Lambda}_{1,r}$, with $\pi_x\left(\kappa_T(\Omega)\right)\subset \comp^2$. Here $\pi_x: T^*\comp^2 \rightarrow \comp^2$ is the
natural projection given by $\pi_x(x,\xi)=x$. Correspondingly, a function $F: \Omega\rightarrow \comp$ may be identified with
$F\circ \left(\pi_x\circ \kappa_T\right)^{-1}: \comp^2\rightarrow \comp$.

\medskip
Let us recall from the introduction that we assume, for simplicity of the exposition only, that the tori
$\Lambda_{j,d}$, $j=1,2$ and $\Lambda_{1,r}$ belong to the same open edge of $J$ in (\ref{eq0}) so that (\ref{eq0.46.44})
holds. Let $\widetilde{\Lambda}_{j,d}\subset \Lambda_{\eps}$, $j=1$, $2$, be "intermediate"
Diophantine tori belonging to the same open edge of $J$ as $\Lambda_{j,d}$ and $\Lambda_{1,r}$, away from
$\supp\, \widetilde{\chi}_d$, with
$$
\Lambda_{1,d} < \widetilde{\Lambda}_{1,d} < \Lambda_{1,r},\quad \Lambda_{1,r} < \widetilde{\Lambda}_{2,d} < \Lambda_{2,d}.
$$
Here, in order to simplify the notation, we are identifying the real tori
$\widetilde{\Lambda}_{j,d}\subset p^{-1}(0)\cap \real^4$ with their images in $\Lambda_{\eps}$, by means of the
canonical transformation $\exp(i\eps H_{G_T}): \real^4\rightarrow \Lambda_{\eps}$ --- see also (\ref{eq82}).
We shall introduce a new weight $G\in C^{\infty}_0(\comp^2)$ supported in a region where
$\abs{\Re P_{\eps}}\leq 1/{\cal O}(1)$, such that $G=0$ in a fixed \neigh{} of $\supp\, \chi_r$, while
$G=-\eta<0$ in a fixed \neigh{} of $\supp\, \widetilde{\chi}_d$. Here $\eta>0$ is very small but fixed and we shall have
$\abs{\nabla G}\ll 1$, $\abs{\nabla^2 G}\ll 1$ everywhere.
Moreover, $G$ will be chosen so that, when restricting the attention to the region where
$\abs{\Re P_{\eps}}\leq 1/{\cal O}(1)$, the support of $\nabla G$ is contained in a sufficiently small but fixed \neigh{}
of $\widetilde{\Lambda}_{1,d}\cup \widetilde{\Lambda}_{2,d}$.

\medskip
\noindent
We shall now define $G$ near $\widetilde{\Lambda}_{j,d}$, say, when $j=1$. When doing
so, take a smooth canonical diffeomorphism
\begeq
\label{eq126.2}
\widetilde{\kappa}: {\rm neigh}(\widetilde{\Lambda}_{1,d}, \Lambda_{\eps})\rightarrow {\rm neigh}(\xi=0, T^*{\bf T}^2),
\endeq
mapping $\widetilde{\Lambda}_{1,d}$ to the zero section in $T^*{\bf T}^2$ and obtained by composing the
action-angle canonical transformation near the real torus with the holomorphic
transformation $\exp(-i\eps H_{G_T})$.
Composing $p_{\eps}$ in (\ref{eq0.47}) with $\widetilde{\kappa}^{-1}$, we
obtain a new symbol, still denoted by $p_{\eps}$, defined near
the zero section $\xi=0$ in $T^*{\bf T}^2$, which is of the form
\begeq
\label{eq126.3}
p_{\eps}(x,\xi)=p(\xi)+i\eps \langle{q}\rangle_T(x,\xi)+{\cal O}_T(\eps^2).
\endeq
Here, as already observed in the beginning of subsection 5.2, we take $T>0$ sufficiently large but fixed,
so that $\langle{q}\rangle_T(x,\xi)$ avoids the value $F_0$ in this
region. In view of the implicit function theorem, we may assume that the energy surface $p^{-1}(0)$ is given by an equation
\begeq
\label{eq126.35}
\xi_2=f(\xi_1),\quad \abs{\xi_1}\leq a,\quad 0<a\ll 1,
\endeq
where the analytic function $f$ satisfies $f(0)=0$, $f'(0)\neq 0$.

When defining the weight $G$ near $\xi=0$ we shall require that it
should be constant on each invariant torus $\xi={\rm Const}$. In doing
so, we shall first define $G$ on $p^{-1}(0)$, and to that end we introduce the tori
$\Lambda_{\mu}\subset p^{-1}(0)$, $\abs{\mu}\leq a$, given by $\xi_1=\mu$, $\xi_2=f(\mu)$. In order to fix the ideas,
let us assume that when $\mu<0$, then the tori $\widetilde{\kappa}^{-1}\left(\Lambda_{\mu}\right)$ satisfy
$$
\widetilde{\Lambda}_{1,d} < \widetilde{\kappa}^{-1}\left(\Lambda_{\mu}\right) < \widehat{\Lambda}_{1,r},
$$
and for $\mu>0$, we have
$$
\widehat{\Lambda}_{1,d} < \widetilde{\kappa}^{-1}\left(\Lambda_{\mu}\right) < \widetilde{\Lambda}_{1,d}.
$$
When $\delta>0$ is very small but fixed, we then let $G_0=G_0(\xi_1)\in C^{\infty}([-a,a];[0,\delta])$ be increasing and such
that $G_0=0$ near $-a$, $G_0=\delta$ near $a$, and with $G'_0$ having a compact support in a small \neigh{} of $\xi_1=0$.
Taking $\delta>0$ small enough, we achieve that $\abs{G_0'}\ll 1$ and $\abs{G''_0}\ll 1$. Setting
$G(\xi_1,f(\xi_1))=-G_0(\xi_1)$, we see that we have defined $G$ on $p^{-1}(0)$. We then extend $G$ suitably to a full
\neigh{} of $\xi=0$ in $\real^2$ so that it still depends on $\xi$ only and $\abs{\nabla G}\ll 1$ is different from zero
only in a small \neigh{} of $\xi=0$.

\medskip
Introduce next the IR-manifold
\begeq
\label{eq126.4}
\left(T^*{\bf T}^2\right)_G=\left\{(x+iG'_{\xi}(\xi),\xi);\,\, (x,\xi)\in T^*{\bf T}^2\right\},
\endeq
defined in a complex \neigh{} of the zero section $\xi=0$. Then the imaginary part of the symbol of
$p_{\eps}$ in (\ref{eq126.3}), along $\left(T^*{\bf T}^2\right)_G$, still avoids the value $\eps F_0$.

Similarly, working in the action-angle variables, we define $G=G(\xi)$ in a \neigh{} of the Diophantine torus
$\widetilde{\Lambda}_{2,d}$. It is then clear that we can define the new global
IR--manifold $\widetilde{\Lambda}_{\eps}\subset \comp^4$ so that near $\widetilde{\Lambda}_{1,d}$, it is given by
$\widetilde{\kappa}^{-1}\left(\left(T^*{\bf T}^2\right)_G\right)$, and away from
$\widetilde{\Lambda}_{1,d}\cup \widetilde{\Lambda}_{2,d}\cup \widehat{\Lambda}_{1,r}$, we define
$\widetilde{\Lambda}_{\eps}$ so that the representation
$$
\kappa_T\left(\widetilde{\Lambda}_{\eps}\right)=\Lambda_{\widetilde{\Phi}_{\eps}}
$$
holds true. Here $\widetilde{\Phi}_{\eps}-\Phi_{\eps}\in C^{\infty}_0(\comp^2)$ and its gradient is supported
in a small \neigh{} of $\widetilde{\Lambda}_{1,d}\cup \widetilde{\Lambda}_{2,d}$, when
restricting the attention to the region where $\abs{\Re P_{\eps}}\leq 1/{\cal O}(1)$.
We have $\widetilde{\Phi}_{\eps}=\Phi_{\eps}-\eta$ in a fixed
\neigh{} of $\supp\, \widetilde{\chi}_d$, while
$\widetilde{\Lambda}_{\eps}=\Lambda_{\eps}$ near $\widehat{\Lambda}_{1,r}$.

\medskip
The discussion above is summarized in the following proposition. \

\begin{lemma} There exists an IR-manifold $\widetilde{\Lambda}_{\eps}\subset T^*\comp^2$ which coincides with
$\Lambda_{\eps}$ near $\widehat{\Lambda}_{1,r}$, such that away from $\widehat{\Lambda}_{1,r}$, after applying
the canonical transformation $\kappa_T$, defined in {\rm (\ref{eq126.12})}, so that $\Lambda_{\eps}$ becomes
$\Lambda_{\Phi_{\eps}}$, with
$$
\Phi_{\eps}=\Phi_0+{\cal O}(\eps),\quad \Phi_0(x)=\frac{(\Im x)^2}{2},
$$
$\widetilde{\Lambda}_{\eps}$ becomes $\Lambda_{\widetilde{\Phi}_{\eps}}$, where $\widetilde{\Phi}_{\eps}-\Phi_{\eps}$
is compactly supported and for some $\eta>0$ small enough but fixed we have
$\widetilde{\Phi}_{\eps}=\Phi_{\eps}-\eta$ near
$\pi_x\left(\kappa_T\left(\supp\, \widetilde{\chi}_d\right)\right)$. Furthermore,
$$
P_{\eps}={\cal O}(1): H(\widetilde{\Lambda}_{\eps})\rightarrow H(\widetilde{\Lambda}_{\eps}),
$$
and the resolvent bounds {\rm (\ref{eq99})} and {\rm (\ref{eq117})} hold true in the sense of bounded linear operators
on $H(\widetilde{\Lambda}_{\eps})$.
\end{lemma}


\medskip
\noindent
It is now easy to estimate the norm of the term (\ref{eq120.7}) as a
bounded operator on $H(\Lambda_{\eps})$. First notice that
$$
\chi_r={\cal O}(1): H(\Lambda_{\eps})\rightarrow H(\widetilde{\Lambda}_{\eps}),
$$
where we use, as before, the Toeplitz quantization of $\chi_r$ on the FBI--Bargmann side.
Combining this with Proposition 5.4 and Lemma 5.6, we get
\begeq
\label{eq126.5}
\left(P_r-z\right)^{-1} \chi_r=\frac{{\cal O}(1)}{\eps}
\exp\left({\cal O}\left(\frac{\eps^{3/2}}{h^2}\right)\log\frac{1}{h}\right):
H({\Lambda}_{\eps})\rightarrow H(\widetilde{\Lambda}_{\eps}).
\endeq
Now $\supp\, \widetilde{\chi}_d$ is contained in a region where
$\widetilde{\Phi}_{\eps}-\Phi_{\eps}=-\eta<0$ and hence,
\begeq
\label{eq126.6}
\widetilde{\chi}_d={\cal O}\left(e^{-\frac{\eta}{h}}\right):
H(\widetilde{\Lambda}_{\eps})\rightarrow H(\Lambda_{\eps}).
\endeq
Here $\widetilde{\chi}_d$ is quantized as a Toeplitz operator in the weighted space $H_{\Phi_{\eps}}$, by working on
the transform side. Using (\ref{eq126.5}) and (\ref{eq126.6}) together with the upper
bound $\eps={\cal O}(h^{2/3+\delta})$, $\delta>0$, we conclude that
\begeq
\label{eq127}
L_1=-i\eps C \widetilde{\chi}_d\left(P_r-z\right)^{-1}\chi_r={\cal O}\left(e^{-1/{\widetilde{C}h}}\right):
H(\Lambda_{\eps})\rightarrow H(\Lambda_{\eps}),
\quad \widetilde{C}>0.
\endeq

\bigskip
\noindent
When estimating the operator norm of the expression
\begeq
\label{eq127.1}
L_2=-i\eps K\left(P_d-z\right)^{-1}\chi_d:
H(\Lambda_{\eps})\rightarrow H(\Lambda_{\eps}),
\endeq
we argue similarly and introduce a weak but $h$-independent weight, supported in a region where $\abs{\Re
P_{\eps}}\leq 1/{\cal O}(1)$, which is equal to a very small strictly positive constant $\eta>0$ in a fixed \neigh{}
of $\supp\, \chi_d$. We then obtain a new microlocally weighted space
$H(\widehat{\Lambda}_{\eps})$ associated to an
IR--manifold $\widehat{\Lambda}_{\eps}$ defined similarly to $\widetilde{\Lambda}_{\eps}$, such that if
$\kappa_T(\widehat{\Lambda}_{\eps})=\Lambda_{\widehat{\Phi}_{\eps}}$, then $\widehat{\Phi}_{\eps} = \Phi_{\eps}+\eta$,
$0<\eta \ll 1$, in a fixed \neigh{} of $\supp\, \chi_d$. Then
\begeq
\label{eq128}
\chi_d={\cal O}\left(e^{-\frac{\eta}{h}}\right):
H(\Lambda_{\eps})\rightarrow H(\widehat{\Lambda}_{\eps}),
\endeq
and combining this estimate together with the fact that
$K={\cal O}(1): H(\widehat{\Lambda}_{\eps})\rightarrow
H(\Lambda_{\eps})$ and with Lemma 5.6, we infer that
\begeq
\label{eq129}
L_2=-i\eps K\left(P_d-z\right)^{-1} \chi_d = {\cal O}\left(e^{-\frac{1}{C h}}\right): H(\Lambda_{\eps})
\rightarrow H(\Lambda_{\eps}),\quad C>0.
\endeq

\medskip
To finish the proof of Proposition 5.5, we only need to estimate the norm of the operator
\begeq
\label{eq129.1}
L_3=i\eps C \psi \left(P_{\eps}+i\eps C \psi-z\right)^{-1}\chi_0:
H(\Lambda_{\eps})\rightarrow H(\Lambda_{\eps}),
\endeq
and this requires an introduction of a new weight on the FBI--Bargmann
transform side, that we shall still denote by $G$. We take $G\in
C^{\infty}_b(\comp^2)$ with $\abs{\nabla G}\ll 1$, $\abs{\nabla^2
G}\ll 1$, such that $G=0$ in a fixed \neigh{} of $\supp\, \psi$. We
shall furthermore choose $G$ so that it is equal to a very small but strictly
positive constant in a fixed \neigh{} of $\supp\, \chi_0$, and hence
in a \neigh{} of infinity. Here we may recall that the support of
$\chi_0$ does not intersect the compact set $\supp\, \psi$. We
also choose $G$ so that the support of $\nabla G$ is contained in a thin
domain included in a region where $\abs{\Re P_{\eps}}\geq 1/{\cal
O}(1)$. It is then easy to see that
\begeq
\label{eq129.2}
L_3 = {\cal O}\left(e^{-1/Ch}\right): H(\Lambda_{\eps})\rightarrow
H(\Lambda_{\eps}),\quad C>0,
\endeq
and combining this estimate together with (\ref{eq127}),
(\ref{eq129}), and (\ref{eq120.5}), we complete the proof of
Proposition 5.5.

\bigskip
Combining Proposition 5.5 with (\ref{eq120}) we see
that for $z$ satisfying (\ref{eq118}), the ope\-rator
$P_{\eps}-z: H(\Lambda_{\eps})\rightarrow H(\Lambda_{\eps})$ is invertible, with
\begeq
\label{eq131}
\left(P_{\eps}-z\right)^{-1}=R_0(z)\left(1+L\right)^{-1}.
\endeq
Writing $(1+L)^{-1}=1-(1+L)^{-1}L$ we get
\begeq
\label{eq132}
\left(P_{\eps}-z\right)^{-1}=R_0(z)-R_0(z)\left(1+L\right)^{-1}L.
\endeq
Let now $\gamma$ be a simple positively oriented closed $C^1$--contour contained in the domain (\ref{eq117.6}),
of length ${\cal O}(\eps)$, such that (\ref{eq118}) holds for each $z$ along $\gamma$.
Let
\begeq
\label{eq133}
\Pi=-\frac{1}{2\pi i}\int_{\gamma} \left(P_{\eps}-z\right)^{-1}\,dz
\endeq
be the spectral projection of $P_{\eps}$ associated to the spectrum of $P_{\eps}$ inside $\gamma$.
The finite-dimensional space $\Pi (H(\Lambda_{\eps}))$ is spanned by the generalized eigenfunctions of $P_{\eps}$
corresponding to the eigenvalues of $P_{\eps}$ in the interior of $\gamma$. Define also
\begeq
\label{eq134}
\Pi_0=-\frac{1}{2\pi i} \int_{\gamma} R_0(z)\,dz,
\endeq
and notice that the last term in the right hand side of (\ref{eq119}) does not contribute to the integral in
(\ref{eq134}), since $(P_{\eps}+i\eps C\psi-z)^{-1}$ is holomorphic in $z\in R_s$. Let us also introduce
the finite-dimensional space $E\subset H(\Lambda_{\eps})$ spanned by the generalized eigenfunctions of the operators
$P_d$ and $P_r$, corresponding to their spectra inside $\gamma$. Notice that the range
of $\Pi_0$ in (\ref{eq134}) is contained in $E$.

Now (\ref{eq132}) gives that
\begeq
\label{eq135}
\Pi=\Pi_0+\frac{1}{2\pi i}\int_{\gamma} R_0(z)\left(1+L\right)^{-1}L\,dz,
\endeq
and combining Proposition 5.3, Proposition 5.4 and Proposition 5.5
together with the fact that $\eps={\cal O}\left(h^{2/3+\delta}\right)$, $\delta>0$, we see that the operator
norm of the contour integral in the right hand side of (\ref{eq135}), is ${\cal O}(\exp(-1/\widehat{C} h))$, for
some $\widehat{C}>0$. In particular, if $u\in H(\Lambda_{\eps})$, $\norm{u}=1$, belongs to the range of $\Pi$, then
$$
\Pi_0 u=u+{\cal O}(e^{-1/\widehat{C} h}).
$$
Using the basic properties of the non--symmetric distance between two closed subspaces of a Hilbert space,
introduced and studied in~\cite{HeSj1} (see also~\cite{DiSj}), we conclude that
\begeq
\label{eq136}
{\rm dim}\, \Pi(H(\Lambda_{\eps})) \leq {\rm dim}\,E.
\endeq

\medskip
When proving the opposite inequality, we write, using (\ref{eq119}),
$$
\Pi_0=\Pi_r \chi_r + \Pi_d \chi_d,
$$
where
\begeq
\label{eq137}
\Pi_d =-\frac{1}{2\pi i}\int_{\gamma} \left(P_d-z\right)^{-1}\,dz={\cal O}\left(\frac{1}{h^{N_0}}\right):
H(\Lambda_{\eps})\rightarrow H(\Lambda_{\eps}),
\endeq
and
\begeq
\label{eq138}
\Pi_r =-\frac{1}{2\pi i}\int_{\gamma}
\left(P_r-z\right)^{-1}\,dz
\endeq
satisfies
\begeq
\label{eq139}
\Pi_r=\exp\left({\cal O}\left(\frac{\eps^{3/2}}{h^2}\right)
\log\frac{1}{h}\right): H(\Lambda_{\eps})\rightarrow
H(\Lambda_{\eps}).
\endeq
Here we have also used Propositions 5.3 and 5.4.

\medskip
Let $u\in E$ be a normalized generalized eigenfunction of, say,
$P_d$, corresponding to an eigenvalue of this operator inside $\gamma$. Then using exponentially weighted estimates,
in the same way as in
the proof of Proposition 5.5, together with (\ref{eq139}) and the upper bound $\eps={\cal O}(h^{2/3+\delta})$,
we see that
$$
\Pi_r \chi_r u={\cal O}(e^{-1/C h}),\quad C>0.
$$
Similarly, we find that $\Pi_d \chi_d u=u+{\cal O}(e^{-1/C h})$, and
therefore,
\begeq
\label{eq140}
\Pi_0 u=u+{\cal O}(e^{-1/C h}).
\endeq
We get the same conclusion also when $u\in E$ is a normalized generalized
eigenfunction of $P_r$.

\medskip
Let now $u\in E$ be such that $\norm{u}=1$. Using (\ref{eq135}) and (\ref{eq140}),
we infer that
$$
\Pi u = u+{\cal O}(e^{-1/C h}),
$$
and it follows that the dimension of $E$ does not exceed that of
$\Pi(H(\Lambda_{\eps}))$.
This together with (\ref{eq136}) implies that the spaces $\Pi(H(\Lambda_{\eps}))$ and $E$ have the same dimension,
and from here it is easy to see how to get the full statement of Theorem 1.1.

\section{An application to surfaces of revolution}
\setcounter{equation}{0}
The purpose of this section is to illustrate how Theorem 1.1 applies
to the case when $M$ is an analytic surface of revolution in
$\real^3$, and
\begeq
\label{eq141}
P_{\eps}=-h^2\Delta+i\eps q,
\endeq
where $\Delta$ is the
Laplace-Beltrami operator and $q$ is an analytic function on $M$. We
shall consider the same class of surfaces of revolution as
in~\cite{HiSjVu}, and begin by recalling the assumptions made on $M$
in that paper.

Let us normalize $M$ so that the $x_3$-axis is
its axis of revolution, and parametrize it by the cylinder
$[0,L]\times S^1$, $L>0$,
\begeq
\label{eq142}
[0,L]\times S^1 \ni (s,\theta)\mapsto (u(s)\cos\theta, u(s)\sin
\theta, v(s)),
\endeq
assuming, as we may, that the parameter $s$ is the arclength
along the meridians, so that $(u'(s))^2+(v'(s))^2=1$. In the coordinates $(s,\theta)$, the
Euclidean metric on $M$ takes the form
\begin{equation}
\label{eq143}
 g=ds^2+u^2(s)d\theta^2.
\end{equation}
The functions $u$ and $v$ are assumed to be real analytic on
$[0,L]$, and we shall assume that for each $k\in \nat$,
$$
u^{(2k)}(0)=u^{(2k)}(L)=0,
$$
and that $u'(0)=1$, $u'(L)=-1$. As we recalled in~\cite{HiSjVu}, these
assumptions guarantee the regularity of $M$ at the poles.

Assume furthermore that $M$ is a simple surface of revolution, in the sense that
$0\leq u(s)$ has precisely one critical point $s_0\in (0,L)$, and
that this critical point is a non-degenerate maximum,
$u''(s_0)<0$. To fix the ideas, we shall assume that
$f(s_0)=1$. Notice that $s_0$ corresponds
to the equatorial geodesic $\gamma_E\subset M$ given by $s=s_0$,
$\theta\in S^1$. This is an elliptic orbit.

\medskip
Writing
$$
T^*\left(M\backslash \{(0,0,v(0)),(0,0,v(L))\}\right)\simeq
T^*\left((0,L)\times S^1\right),
$$
and using (\ref{eq143}) we see that the leading symbol of $P_0=-h^2\Delta$ on $M$ is given by
\begin{equation}
\label{eq144}
p(s,\theta,\sigma,\theta^*)=\sigma^2+\frac{(\theta^*)^2}{f^2(s)}.
\end{equation}
Here $\sigma$ and $\theta^*$ are the dual variables to $s$ and
$\theta$, respectively. Since the function $p$ in (\ref{eq144})
does not depend on $\theta$, it follows that $\{p,\theta^*\}=0$, and we recover the well-known fact that the
geodesic flow on $M$ is completely integrable.

\medskip
Let $E>0$ and $\abs{F}< E^{1/2}$, $F\neq 0$. Then the set
$$
\Lambda_{E,F}: p=E,\,\, \theta^*=F,
$$
is an analytic Lagrangian torus contained inside the energy surface
$p^{-1}(E)$. Geometrically, the torus $\Lambda_{E,F}$ consists of
geodesics contained between and intersecting tangentially the
parallels $s_{\pm}(E,F)$ on $M$ defined by the equation
$$
u(s_{\pm}(E,F))=\frac{\abs{F}}{E^{1/2}}.
$$
For $F=0$, the parallels reduce to the two poles and we obtain a
torus consisting of a family of meridians. The case
$\abs{F}=E^{1/2}$ is degenerate and corresponds to the equator
$s=s_0$, traversed with the two different orientations. Writing
$\Lambda_a:=\Lambda_{1,a}$, we get a decomposition as in (\ref{eq0}),
$$
p^{-1}(1)=\bigcup_{a\in J} \Lambda_a,
$$
with $J = [-1,1]$, $S=\{\pm 1\}$.

\medskip
In~\cite{HiSjVu}, we have derived an explicit expression for the
rotation number $\omega(\Lambda_a)$ of the torus $\Lambda_{a}$, $0\neq a\in
(-1,1)$,
\begeq
\label{eq145}
\omega(\Lambda_a)=\frac{a}{\pi}\int_{s_-(a)}^{s_+(a)}
\frac{1}{u^2(s)}\left(1-\frac{a^2}{u^2(s)}\right)^{-1/2}\,ds,\quad
u(s_{\pm}(a))=\abs{a}.
\endeq
We are going to assume that the analytic function $(-1,1)\ni a\mapsto \omega(\Lambda_a)$
is not identically constant.

Let $\alpha>0$, $d>0$. In what follows we shall say that a torus $\Lambda_a\subset
p^{-1}(1)$, $a\in (-1,1)$, is $(\alpha,d)$--Diophantine if the rotation number $\omega(\Lambda_a)$
satisfies
\begeq
\abs{\omega(\Lambda_a)-\frac{p}{q}}\geq \frac{\alpha}{q^{2+d}},\quad
p\in\z,\,\, q\in \nat.
\endeq
From the introduction let us also recall that if a torus $\Lambda_a\subset p^{-1}(1)$ is rational, so that
$\omega(\Lambda_a)=\frac{m}{n}$, with $m\in \z$ and $n\in \nat$ relatively prime and $m={\cal O}(n)$, then
we define the height of $\omega(\Lambda_a)$ as $k(\omega(\Lambda_a))=\abs{m}+\abs{n}$.

\medskip
\noindent
Let $q=q(s,\theta)$ be a real-valued analytic function on $M$ which we shall view as
a function on $T^*M$. Associated to each $a\in J$, we introduce the
compact interval $Q_{\infty}(\Lambda_a)\subset \real$ defined as in (\ref{eq0.1}). We also define
an analytic function
$$
(-1,1)\ni a\mapsto \langle{q}\rangle(\Lambda_a),
$$
obtained by averaging $q$ over the invariant tori $\Lambda_a$. Assume
that $a\mapsto \langle{q}\rangle(\Lambda_a)$ is not identically constant.
From the introduction, let us recall that as
$a\rightarrow a_0\in S$, the set of the accumulation points of $\langle{q}\rangle(\Lambda_a)$ is contained in
$Q_{\infty}(\Lambda_{a_0})$.

\medskip
\noindent
Following~\cite{HiSjVu}, we now come to introduce uniformly good values in $\real$, for which the conclusion of Theorem 1.1
will be valid uniformly. In doing so, let us notice that the following discussion is not restricted to the case
of surfaces of revolution.

Let $d>0$ be fixed. Given $\alpha$, $\beta$, $\gamma>0$ we say that $F_0\in \real$ is $(\alpha,\beta,\gamma)$--good if the following conditions hold:
\begin{itemize}
\item $F_0$ is not in the union of all $Q_{\infty}(\Lambda_a)$ with ${\rm dist}(\Lambda_a,S)\leq \alpha$.
\item If $F_0\in Q_{\infty}(\Lambda_a)$ and $\omega(\Lambda_a)\notin {\bf Q}$ then
$\Lambda_a$ is $(\alpha,d)$--Diophantine and $\abs{d_a\langle{q}\rangle(\Lambda_a)}\geq \alpha$.
\item If $F_0\in Q_{\infty}(\Lambda_a)$ and $\omega(\Lambda_a)\in {\bf Q}$ then
$k(\omega(\Lambda_a))={\cal O}(\frac{1}{\alpha})$,
$\abs{d_a \omega(\Lambda_a)}\geq \alpha$, and $\abs{F_0-\langle{q}\rangle(\Lambda_a)}\geq \alpha$.
\item Let $\langle{q}\rangle^{-1}(F_0)=\{ \Lambda_{a_1,d},\ldots\, \Lambda_{a_L,d}\}$,
$\omega(\Lambda_{a_j,d})\notin {\bf Q}$, $1\leq j \leq L$,
and $F_0\in Q_{\infty}(\Lambda_{a_j,r})$, $\omega(\Lambda_{a_j,r})\in {\bf Q}$, $j=1,\ldots L'$. Then the distance in
$\real$ from $F_0$ to the union
$$
\bigcup_{\Lambda_a\in J;\, {\rm dist}_{J}\left(\Lambda_a, \left(\cup_{j=1}^{L} \Lambda_{a_j,d}\right)\cup
\left(\cup_{k=1}^{L'}\Lambda_{a_k,r}\right)\right)> \beta} Q_{\infty}(\Lambda_a)
$$
is $> \gamma$.
\end{itemize}

\medskip
\noindent
{\it Remark}. This definition of an $(\alpha,\beta,\gamma)$--good value is less restrictive than in our previous
work~\cite{HiSjVu}, since we now allow such a value $F_0$ to belong to an interval $Q_{\infty}(\Lambda_a)$ corresponding
to a rational torus satisfying the isoenergetic condition, provided that $F_0$ is not too close to the torus average
$\langle{q}\rangle(\Lambda_a)$.

\medskip
In the following proposition we shall make use of the fact, observed in the introduction, that in the
case when the subprincipal symbol of $P_{\eps=0}$ in (\ref{eq0012}) vanishes, the
validity of Theorem 1.1 extends to the range $h^2\ll \eps={\cal
O}(h^{2/3+\delta})$.

\medskip
\noindent
\begin{prop}
Assume that $M$ is a simple analytic surface of revolution with a parametrization {\rm (\ref{eq142})}, for which
the rotation number $\omega(\Lambda_a)$ defined in {\rm (\ref{eq145})} is not identically
constant. Consider an operator of the form $P_{\eps}=-h^2\Delta+i\eps
q$, where $q$ is a real valued analytic function on $M$, such that the
torus averages function $a\mapsto \langle{q}\rangle(\Lambda_a)$ is not identically constant. Let $\alpha$, $\beta$,
$\gamma>0$, and fix $0<\delta \ll 1$. There exists $C>0$ such that if $F_0$ is $(\alpha,\beta,\gamma)$--good,
$0<h \leq \frac{1}{C}$, and $h^2/C\leq \eps\leq h^{2/3+\delta}$, then Theorem {\rm 1.1} applies uniformly to describe
the spectrum of $P_{\eps}$ in the rectangle
$$
\left[-\frac{\eps}{C},\frac{\eps}{C}\right]+i\eps
\left[F_0-\frac{\eps^{\delta}}{C},
F_0+\frac{\eps^{\delta}}{C}\right].
$$
\end{prop}

\medskip
\noindent
{\it Remark}. If $\eps=h$, then the operator $P_{\eps}$ in
Proposition 6.1 is a semiclassical version of the stationary
damped wave operator~\cite{Le},~\cite{Sj00},~\cite{Hi02}.

\medskip
\noindent
{\it Remark}. In the corresponding discussion in subsection 7.2
of~\cite{HiSjVu}, it has been assumed that the complex perturbation
$q$ in Proposition 6.1 is close to a rotationally symmetric one. This additional
assumption has now been removed, thanks to Theorem 1.1, at the expense
of weakening the final result and restricting the bounds on the strength $\eps$ of the
non-selfadjoint perturbation.

\appendix
\section {Trace class estimates for Toeplitz operators}
\setcounter{equation}{0}

The purpose of this appendix is to derive a simple estimate on the trace class norm of a Toeplitz operator with a
compactly supported smooth symbol acting in a weighted $L^2$--space of holomorphic functions on $\comp^n$.
Indeed, the result will be seen to be a straightforward consequence of the analysis of~\cite{MeSj}.

\medskip
Let $\Phi_0(x)$ be a real quadratic form on $\comp^n$ and assume that $\Phi_0$ is strictly
plurisubharmonic. (In what follows we may think of the special case when $\Phi_0(x)=\frac{1}{2}(\Im x)^2$.) Let
\begeq
\label{a1}
H_{\Phi_0}:={\rm Hol}(\comp^n)\cap L^2(\comp^n; e^{-\frac{2\Phi_0}{h}}L(dx)),
\endeq
where $L(dx)$ is the Lebesgue measure on $\comp^n=\real^{2n}$ and ${\rm Hol}(\comp^n)$ is the
space of entire holomorphic functions on $\comp^n$. Then $H_{\Phi_0}$ is a closed subspace of the space
$L^2_{\Phi_0}:=L^2(\comp^n; e^{-\frac{2\Phi_0}{h}}L(dx))$, and from~\cite{MeSj}we recall the following expression for
the orthogonal projection $\Pi_{\Phi_0}: L^2_{\Phi_0}\rightarrow H_{\Phi_0}$,
\begeq
\label{a2}
\Pi_{\Phi_0} u(x)=\frac{C}{h^n} \int e^{\frac{2}{h} \psi_0(x,y)} u(y) e^{-\frac{2}{h} \Phi_0(y)}\,L(dy),
\endeq
where the constant $C$ is real and $\psi_0(x,y)$ is the unique quadratic form on $\comp^n_x\times \comp^n_y$ which
is holomorphic in $x$, anti-holomorphic in $y$, and satisfies
\begeq
\label{a3}
\psi_0(x,x)=\Phi_0(x).
\endeq
In the case when $\Phi_0(x)=\frac{1}{2}(\Im x)^2$, we have $\psi_0(x,y)=-\frac{1}{8}(x-\overline{y})^2$.

Now let $\Phi\in C^{\infty}(\comp^n;\real)$ be such that $\Phi-\Phi_0$ is bounded and
$\sup\abs{\frac{\partial \Phi}{\partial x}-\frac{\partial \Phi_0}{\partial x}}$ small enough. Assume also that
$\nabla^k \Phi$ is bounded for each $k\geq 2$ and that $\Phi$ is uniformly strictly plurisubharmonic, so that the
set
\begeq
\label{a4}
\Lambda_{\Phi}=\left\{\left(x,\frac{2}{i}\frac{\partial \Phi}{\partial x}(x)\right); x\in \comp^n\right \}
\endeq
is an IR-manifold. Associated with the weight $\Phi$ we have the orthogonal projection
\begeq
\label{a5}
\Pi_{\Phi}: L^2_{\Phi}\rightarrow H_{\Phi},
\endeq
where $L^2_{\Phi}=L^2(\comp^n; e^{-\frac{2\Phi}{h}}L(dx))$ and $H_{\Phi}={\rm Hol}(\comp^n)\cap L^2_{\Phi}$.
If now $p\in C^{\infty}_0(\comp^n)$, we introduce the corresponding Toeplitz operator
\begeq
\label{a6}
{\rm Top}(p)=\Pi_{\Phi} p \Pi_{\Phi}={\cal O}(1): H_{\Phi}\rightarrow H_{\Phi}.
\endeq
Our goal is to show that ${\rm Top}(p)$ is of trace class as an operator on $H_{\Phi}$ and to estimate its trace class
norm. In doing so, it is convenient to recall from~\cite{MeSj} the asymptotic description of the Bergman projection
$\Pi_{\Phi}$, as $h\rightarrow 0$.

Let $\psi(x,y)\in C^{\infty}(\comp^n_x\times \comp^n_y)$ be al\-most ho\-lo\-mor\-phic in $x$ and
al\-most anti-ho\-lo\-mor\-phic in $y$ at the diagonal ${\rm diag}(\comp^n_x\times \comp^n_y)$, such that
$\nabla^k \psi$ is bounded on $\comp^{2n}$ for each $k\geq 2$ and with
\begeq
\label{a7}
\psi(x,x)=\Phi(x).
\endeq
Then we know that
\begeq
\label{a8}
\Phi(x)+\Phi(y)-2\Re \psi(x,y)\sim \abs{x-y}^2,
\endeq
uniformly for $\abs{x-y}\leq 1/C$, for $C>0$ large enough.

It follows from~\cite{MeSj} that there exists
$f(x,y;h)\sim \sum_{j=0}^{\infty} f_j(x,y) h^j$ in $C^{\infty}_{{b}}(\comp^{2n})$, with $\supp f\subset
\{(x,y); \abs{x-y}\leq 1/C\}$, $C\gg 1$, with $f(x,x;h)$ real, $1/C\leq f_0(x,x)\leq C$, and with
\begeq
\label{a9}
\partial_{\overline{x},y} f={\cal O}\left(\abs{x-y}^{\infty}+h^{\infty}\right),
\endeq
such that if
\begeq
\label{a10}
\widetilde{\Pi}_{\Phi} u(x)=\frac{1}{(\pi h)^n} \int e^{\frac{2}{h}(\psi(x,y)-\Phi(y))} f(x,y;h) u(y)\,L(dy),
\endeq
then
\begeq
\label{a11}
\Pi_{\Phi}=\widetilde{\Pi}_{\Phi}+R.
\endeq
Here
\begeq
\label{a12}
R=e^{\frac{\Phi}{h}} \widetilde{R} e^{-\frac{\Phi}{h}},
\endeq
where $\widetilde{R}$ is a negligible integral operator in the sense of section 3 of~\cite{MeSj}. In particular,
it follows from~\cite{MeSj} that $R={\cal O}(h^{\infty}): L^2_{\Phi} \rightarrow L^2_{\Phi}$. It follows
furthermore from the results of Section 3 of~\cite{MeSj} that the operators $R p$ and $p R$ are of
trace class as operators $H_{\Phi}\rightarrow L^2_{\Phi}$, with the trace class norm ${\cal O}(h^{\infty})$.

When estimating the trace class norm of ${\rm Top}(p)$ on $H_{\Phi}$, we may therefore replace $\Pi_{\Phi}$ by
$\widetilde{\Pi}_{\Phi}$, and consider the corresponding operator
$\widetilde{{\rm Top}}(p)=\widetilde{\Pi}_{\Phi} p \widetilde{\Pi}_{\Phi}$. The factorization
\begeq
\label{a13}
\widetilde{{\rm Top}}(p)=\left(\widetilde{\Pi}_{\Phi} p_1\right) \left(p_2 \widetilde{\Pi}_{\Phi}\right),
\endeq
where $p=p_1 p_2$ and $p_2=\abs{p}^{1/2}$, shows that it suffices to
prove that
the operators $\widetilde{\Pi}_{\Phi} p_1:L^2_{\Phi}\rightarrow H_{\Phi}$ and
$p_2\Pi_{\Phi}: H_{\Phi}\rightarrow L^2_{\Phi}$ are of Hilbert-Schmidt class.
Now the reduced kernel of $\widetilde{\Pi}_{\Phi} p_1$, in view of (\ref{a10}), is equal to
\begeq
\label{a14}
\frac{1}{(\pi h)^n} f(x,y;h) p_1(y) e^{-\frac{\Phi(x)}{h}} e^{\frac{2}{h}\left(\psi(x,y)-\Phi(y)\right)}
e^{\frac{\Phi(y)}{h}},
\endeq
and using also (\ref{a8}) we immediately see that the square of its $L^2$--norm over $\comp^{2n}$ is bounded by
\begeq
\label{a15}
\frac{{\cal O}(1)}{h^{2n}} \int\!\!\!\int \abs{p(y)} e^{-c\abs{x-y}^2/h}\,L(dy)\,L(dx)=
\frac{{\cal O}(1)}{h^n} \norm{p}_{L^1},\quad c>0.
\endeq
It follows that $\widetilde{\Pi}_{\Phi} p_1: L^2_{\Phi} \rightarrow H_{\Phi}$ is of Hilbert-Schmidt class with
\begeq
\label{a16}
\norm{\widetilde{\Pi}_{\Phi}p_1}_{{\rm HS}} = \frac{{\cal O}(1)}{h^{n/2}}\norm{p}_{L^1}^{1/2}.
\endeq
Since a similar argument applies to $p_2\widetilde{\Pi}_{\Phi}$, we get the following result.

\begin{prop}
When $\Phi\in C^{\infty}(\comp^n;\real)$ is a strictly plurisubharmonic function sa\-tisfying the general assumptions
of the beginning of this section, let
$\Pi_{\Phi}: L^2_{\Phi}\rightarrow H_{\Phi}$ be the orthogonal projection. If $p\in C^{\infty}_0(\comp^n)$,
then the Toeplitz operator ${\rm Top}(p)=\Pi_{\Phi} p \Pi_{\Phi}: H_{\Phi}\rightarrow H_{\Phi}$
is of trace class and we have
\begeq
\label{a17}
\norm{{\rm Top}(p)}_{{\rm tr}} \leq \frac{{\cal O}(1)}{h^n} \norm{p}_{L^1}+{\cal O}(h^{\infty}),
\endeq
where the implicit constant in ${\cal O}(h^{\infty})$ is a continuous seminorm of $p$ on the Schwartz space
${\cal S}(\comp^n)$.
\end{prop}

\noindent
{\it Remark}. Rather than working on all of $\comp^n$, we could also consider an open domain
$\Omega\subset \comp^n$, with $p\in C^{\infty}_0(\Omega)$. Then Proposition A.1 remains valid if we replace
$\comp^n$ by $\Omega$ in (\ref{a1}), with $\Pi_{\Phi}$ still being the orthogonal projection on all of $\comp^n$.
In the main text, we work on $H_{\Phi}(\Omega)$, where $\Omega\subset \comp^n/2\pi \z^n$ is open, and Proposition A.1
then still holds.

\end{document}